\documentclass[11pt]{amsart}
\usepackage{Macros}

\title[Tableau Stabilization]{Tableau Stabilization and Rectangular Tableaux Fixed by Promotion Powers}
\author{Connor Ahlbach}
\thanks{Department of Mathematics, Texas State University, San Marcos, TX 78666, USA; c\_a518@txstate.edu}
\date{\today}

\begin{document}

\begin{abstract} We introduce tableau stabilization, a new phenomenon and statistic on Young tableaux based on jeu de taquin. We investigate bounds for tableau stabilization, the shape of stabilized tableaux, and tableau stabilization as a permutation statistic. We apply tableau stabilization to construct the sufficiently large rectangular tableaux fixed by powers of promotion, which were counted by Brendon Rhoades via the cyclic sieving phenomenon \cite[Theorem~1.3]{MR2557880}. 

\end{abstract}
\maketitle

\section{Introduction} 

In this paper, we introduce tableau stabilization, a new phenomenon we found in order to construct sufficiently large rectangular tableaux that are fixed by promotion powers. Central to defining and investigating tableau stabilization are Sch\"utenberger's jeu de taquin and the rectification operator, which are already well-established algorithms in the theory of Young tableaux. Tableau stabilization is the phenomenon that if we attach sufficiently many shifted copies of a skew tableau to its right and then rectify, some copy and all of those to its right only experience horizontal slides. We will investigate when the vertical slides stop, i.e. when the skew tableau stabilizes. We will also determine the shape of the stabilized tableau if the initial skew tableau has same size rows. The case where all rows have the same size is used to construct sufficiently large rectangular tableaux that are fixed by promotion powers. We leave the reader with open problems on tableau stabilization, most notably extending our bound from equal row sizes and finding its distribution as a permutation statistic. See \Cref{sec:bkgd} for definitions and further background.

\begin{Definition} \label{def:stab} For any standard skew tableau $ S $, let $ S^{(k)} $ denote the result of attaching $ (k - 1) $ shifted copies of $ S $ to the right of $ S $ so that the result is a standard skew tableau. Let $ m $ denote the size of $ S $ and $ k $ be a positive integer. We say $ S $ \emph{stabilizes} at $ k $ if the entries in $ [(k - 1)m + 1, km] $ lie in the same rows in $ \Rect(S^{(k)}) $ and $ S^{(k)} $. Let $ \stab(S) $ denote the minimum value at which $ S $ stabilizes.

\end{Definition}

\begin{Example} \label{ex:stab} Let $ \Rect $ denote the rectification operator. Consider
\begin{align*}
	S & = \byt \none & \none & 1 & 3 \\ \none & \none & 5 & 6 \\ 2 & 4 \eyt \, \qquad \qquad \qquad \quad \; \;
	T = \byt \none & \none & 1 & 6 \\ \none & 2 & 5 \\ 3 & 4 \eyt \\
	S^{(3)} & = \byt \none & \none & 1 & 3 &  *(yellow) 7 &  *(yellow) 9 & *(green) 13 & *(green) 15 \\ \none & \none & 5 & 6 &  *(yellow) 11 &  *(yellow) 12 & *(green) 17 & *(green) 18 \\ 2 & 4 &  *(yellow) 8 &  *(yellow) 10 & *(green) 14 & *(green) 16 \eyt \, , \quad 
	T^{(3)} = \byt \none & \none & 1 & 6 & *(yellow) 7 & *(yellow) 12 & *(green) 13 & *(green) 18 \\ \none & 2 & 5 & *(yellow) 8 & *(yellow) 11 & *(green) 14 & *(green) 17 \\ 3 & 4 & *(yellow) 9 & *(yellow) 10 & *(green) 15 & *(green) 16 \eyt \, , \\
	\Rect(S^{(3)}) & =  \byt 1 & 3 & 5 & 6 &  *(yellow) 7 &  *(yellow) 9 & *(green) 13 & *(green) 15 \\ 2 & 4 &  *(yellow) 11 &  *(yellow) 12 & *(green) 17 & *(green) 18 \\  *(yellow) 8 &  *(yellow) 10 & *(green) 14 & *(green) 16 \eyt \, , \quad 			
	\Rect(T^{(3)}) = \byt 1 & 4 & 5 & 6 & *(yellow) 7 & *(yellow) 12 & *(green) 13 & *(green) 18 \\ 2 & *(yellow) 8 & *(yellow) 10 & *(yellow) 11 & *(green) 14 & *(green) 17 \\ 3 & *(yellow) 9 & *(green) 15 & *(green) 16 \eyt \, . 
\end{align*}
Note that $ 2 $ does not lie in the same row in $ S^{(3)} $ and $ \Rect(S^{(3)}) $, but $ 7, 8, \dots, 12 $ do. Hence, $ \stab(S) = 2 $. As $ 13, 14, \dots, 18 $ also stay in the same row in $ S^{(3)} $ and $ \Rect(S^{(3)}) $, $ S $ stabilizes at 3 as well. Note $ 10 $ does not lie in the same row in $ T^{(3)} $ and $ \Rect(T^{(3)}) $, but $ 13, 14, \dots, 18 $ do. Hence, $ \stab(T) = 3 $. 

Also consider
\begin{align*}
	U & = \byt \none & \none & \none & 4 & 5 & 6 \\ \none & \none & 3 & 7 \\ 1 & 2 \eyt \, , \\
	U^{(3)} & = \byt \none & \none & \none & 4 & 5 & 6 & *(yellow) 11 & *(yellow) 12 & *(yellow) 13 & *(green) 18 & *(green) 19 & *(green) 20 \\ \none & \none & 3 & 7 & *(yellow) 10 & *(yellow) 14 & *(green) 17 & *(green) 21 \\ 1 & 2 & *(yellow) 8 & *(yellow)  9 & *(green) 15 & *(green) 16 \eyt \, , \\
	\Rect(U^{(3)}) & = \byt 1 & 2 & 3 & 4 & 5 & 6 & *(yellow) 11 & *(yellow) 12 & *(yellow) 13 & *(green) 18 & *(green) 19 & *(green) 20 \\ 7 & *(yellow) 9 & *(yellow) 10 & *(yellow) 14 & *(green) 17 & *(green) 21 \\ *(yellow) 8 & *(green) 15 & *(green) 16 \eyt \, .
\end{align*}
Note that 9 does not in the same row in $ U^{(3)} $ and $ \Rect(U^{(3)}) $, but $ 15, 16, \dots, 21 $ do, so $ \stab(U) = 3 $.

\end{Example}

Tableau stabilization is defined on skew tableaux whose row sizes weakly decrease from top to bottom. Otherwise, $ S^{(k)} $ need not be a standard skew tableau, see \Cref{rem:rowsizes}. Two of the most basic facts about tableau stabilization are that once a skew tableau stabilizes, it continues to stabilize, and any skew tableau must stabilize eventually, see \Cref{lem:stabcontinues}. A more interesting property of stabilization is that it is constant on dual equivalence classes, see \Cref{thm:dualequivstab}.

The case where all rows have the same size is of particular interest to us because we use it to construct sufficiently large rectangular tableaux that are fixed by certain promotion powers. We also have better results in this case.

It is natural to ask what bounds there are for when a tableau stabilizes. In \Cref{ex:stab}, $ S, T, U $ all have 3 rows and stabilized at 3. Note that $ S $ in \Cref{ex:stab} actually stabilizes earlier. Does every skew tableau with $ b $ rows whose sizes weakly decrease from top to bottom stabilize at $ b $? We have verified that the answer is yes for all standard skew tableaux of size at most 7 and random searching on larger tableaux has failed to produce a counterexample (computations in Sage, see \cite{Sage}).

\begin{Conjecture} \label{conj:stab} Any standard skew tableau with $ b $ rows and decreasing row sizes stabilizes at $ b $.

\end{Conjecture}

In the case where all rows have the same size, \Cref{conj:stab} is true, see \Cref{thm:stab}. \Cref{thm:stab} will give us explicit bounds on the dimensions of the rectangular tableaux we construct to be fixed by various promotion powers. Although we have not proven \Cref{conj:stab}, we have deduced a weaker, but still linear, bound for skew tableaux with $ b $ rows of decreasing size, see \Cref{thm:stabweaklydecrows}. This weaker bound shows \Cref{conj:stab} is true when $ b = 2 $. The case $ b = 3 $ is open.

\begin{Thm} \label{thm:stab} Any standard skew tableau with $ b $ rows of the same size stabilizes at $ b $.

\end{Thm}

\begin{Thm} \label{thm:stabweaklydecrows} Any standard skew tableau with $ b \ge 2 $ rows and decreasing row sizes stabilizes at $ 2b - 2 $.

\end{Thm}

We have a way to determine the shape of $ \Rect(S^{(k)}) $ if $ S $ has $ b $ rows of size $ r $, and $k \ge b - 1 $, \Cref{thm:stabshape}. This approach is instrumental in proving \Cref{thm:stab}. Both \Cref{thm:stab} and \Cref{thm:stabshape} are essential in our proof of \Cref{thm:prfixedtabs}. See \Cref{sec:bkgd} and \Cref{sec:stab} for missing definitions. Purbhoo and Rhee effectively proved the case $ r = 1 $ in \cite[Lemma 11(ii)]{MR3625918}. 

\begin{Thm} \label{thm:stabshape} ($ r = 1 $, see \cite[Lemma 11(ii)]{MR3625918}) Suppose $ S $ is a standard skew tableau with $ b $ rows of size $ r $. Let $ w_1, \dots, w_b $ denote the entries in each row read from left to right, starting from the bottom. For $ k \ge b -1 $, $ \Rect(S^{(k)}) $ has shape $ (\lam_1, \dots, \lam_b) $, where
\begin{align}
\label{eq:stabshape}
	\lam_j = kr + \sum_{i = 1}^{b - j} c_i - \sum_{i = 1}^{j - 1} c_i \quad \tx{ for all } j = 1, \dots, b,
\end{align}
and
\[
	c_i = (\tx{the length of the first row of } P(w_i w_{i + 1})) - r \quad \tx{ for all } i = 1, \dots, b - 1.
\]
\end{Thm}

We defined tableau stabilization in order to construct the sufficiently large rectangular tableaux fixed by powers of promotion. Dennis White conjectured and Brendon Rhoades proved a very interesting cyclic sieving phenomenon, see \cite{MR2087303}, regarding the action of promotion on rectangular standard Young tableaux, see \Cref{thm:rhoadesCSP} \cite[Theorem~1.3]{MR2557880}. This result shows that the number of rectangular standard Young tableaux of shape $ (a^b) $ fixed by $ d $ promotions equals the number of standard $ \f{ab}{d} $-ribbon tableaux of shape $ (a^b) $, see \Cref{cor:prfixedcount} \cite[Corollary 9.1]{MR2557880}, which is the same as the number of standard tableaux of the shape associated to the $ \f{ab}{d} $-quotient of $ (a^b) $, see \cite{MR1399504}, \Cref{cor:prfixedribbontab}.

\begin{Definition} \label{def:fpsets} For any set $ W $ and map $ g: W \to W $, define
\[
	W^g \coloneqq \{ w \in W : g(w) = w \}.
\]
\end{Definition}

\begin{Definition}
  \label{def:csp}
  Suppose $C_n$ is a cyclic group of order $n$ generated
  by $\sigma_n$, $W$ is a finite set on which $C_n$ acts,
  and $f(q) \in \bZ_{\geq 0}[q]$. We say the triple
  $(W, C_n, f(q))$ exhibits the
  \textit{cyclic sieving phenomenon (CSP)} if for all
  $k \in \bZ$,
  \begin{align}\label{eq:CSP_eval}
      \# W^{\si_n^k} = f(\omega_n^k),
  \end{align}
  where $ \omega_n $ is a fixed primitive $ n $-th root of unity.
\end{Definition}

For $ a, b \in \bZ_{\ge 1} $, let
\[
	(a^b) \coloneqq \underbrace{(a, \dots, a)}_{b \tx{ times}},
\]
and let $ \SYT(a^b) $ denote the set of standard Young tableaux of shape $ (a^b) $.

\begin{Thm} \cite[Theorem~4.4]{MR1158783} \label{thm:prnid} Sch\"utzenberger's promotion operator $ \pr: \SYT(a^b) \to \SYT(a^b) $ has order $ ab $.

\end{Thm}

\noindent By \Cref{thm:prnid}, the cyclic group $ \la \pr \ra $ generated by $ \pr: \SYT(a^b) \to \SYT(a^b) $ has order $ ab $. Following a conjecture by Dennis White, Rhoades proved that using the $ q $-analog of the hook length formula on shape $ (a^b) $ gives rise to a CSP for the action of $ \la \pr \ra $ \cite[Theorem~1.3]{MR2557880}. Let $ h_c $ denote the hook length of cell $ c $ in the Young diagram of $ (a^b) $.

\begin{Thm}{\cite[Theorem~1.3]{MR2557880}} (Conjectured by White, 2007) \label{thm:rhoadesCSP} For $ a, b \in \bZ_{\ge 1} $,
\[
	\lp \SYT(a^b), \la \pr \ra, \f{[ab]_q!}{\prod_{c \in (a^b)} [h_c]_q} \rp
\]
exhibits the CSP.

\end{Thm}

The polynomial $ \f{[ab]_q!}{\prod_{c \in (a^b)} [h_c]_q} $ in \Cref{thm:rhoadesCSP} has notable connections to tableau statistics and representation theory \cite{MR2557880}. First, the hook length formula \cite{MR0062127} tells us
\[
	\# \SYT(a^b) = \f{(ab)!}{\prod_{c \in (a^b)} h_c},
\]
making $ \f{[ab]_q!}{\prod_{c \in (a^b)} [h_c]_q} $ the $ q $-analog of the hook length formula for the shape $ (a^b) $. Secondly, it is a $ q $-shift of the major index generating function on Standard Young tableau of shape $ \lam $:
\[
	 \f{[ab]_q!}{\prod_{c \in (a^b)} [h_c]_q} = q^{- a \ch{b}{2}} \SYT(\lam)^{\maj}(q).
\]
Thirdly, $ \SYT(\lam)^{\maj}(q) $ corresponds to the graded multiplicities of the Specht module $ S^{\lam} $ in the coinvariant algebra \cite[Corollary 7.21.5]{MR1676282}. Namely, for $ \lam \vdash n $,
\[
	\SYT(\lam)^{\maj}(q) = \sum_{k \ge 0} \la R_{n}^k, S^{\lam} \ra q^k, 
\]
where $ R_n^k $ is the component of degree $ k $ of the coinvariant algebra $ R_n $ of $ S_n $. 

Rhoades proves \Cref{thm:rhoadesCSP} by finding a basis of the Specht module $ S^{(a^b)} $ on which the long cycle $ \si_{ab} $ acts by promotion up to a sign, see \cite[Proposition~3.5]{MR2557880}. He uses the Kazhdan--Lusztig construction of the irreducible $ S_n $-representations, which is governed by descent sets of tableaux and the leading coefficients of the Kazhdan--Lusztig polynomials. Identifying permutations with their insertion tableaux, he shows that the symmetrized Kazhdan--Lusztig $ \mu $ function on rectangular tableaux is invariant under simultaneous promotion. Moreover, he defines a cyclic descent set on rectangular tableaux which promotion cycles. 

Because $ \si_{ab} $ acts by promotion up to a sign on the Specht module indexed by $ (a^b) $, we have, by the definition of character,
\begin{align}
\label{eq:FPchar}
	\# \SYT(a^b)^{\pr^d} = | \chi^{(a^b)}(\si_{ab}^d) |,
\end{align}
where $ \chi^\lam $ is the character of the Specht module indexed by $ \lam $. The following Corollary then follows from \eqref{eq:FPchar}, the Murnagnan--Nakayama Rule, and the fact that all $ r $-ribbon tableaux of a given shape have the same height parity \cite[2.7.26]{MR644144}. 

\begin{Corollary}{\cite[Corollary~9.1]{MR2557880}} \label{cor:prfixedcount} For $ a, b \in \bZ_{\ge 1} $, and $ d \mid ab $,
\begin{align*}	
	\# \SYT(a^b)^{\pr^d} = \# \tx{ standard } \lp \f{ab}{d} \rp \tx{$-$ribbon tableaux of shape } (a^b).	
\end{align*}
\end{Corollary}

Kevin Purbhoo also gave an alternate proof of \Cref{cor:prfixedcount} using the Wronksi map \cite[Theorem 1.5]{MR3009651}. The Wronksi map takes a $ b $-dimensional subspace $ X $ of polynomials with degree up to $ a + b - 1 $ and outputs the determinant of the Jacobian matrix of a basis for $ X $, which is well-defined up to scalar multiplication. The generic fibers of the Wronski map are in bijection with rectangular standard Young tableaux of shape $ (a^b) $. If we restrict to points in the fiber of the Wronski map that are fixed by a certain $ C_d $-action, Purbhoo shows that the generic number is both $ \# \SYT(a^b)^{\pr^d} $ and the number of $ \f{ab}{d} $-ribbon tableaux of shape $ (a^b) $, proving \Cref{cor:prfixedcount}.

Now, $ r $-ribbon tableaux of shape $ \lam $ only exist when $ \lam $ has empty $ r $-core \cite{MR1399504}. Moreover, when $ \lam $ has empty $ r $-core, $ r $-ribbon tableaux of shape $ \lam $ bijectively correspond to standard fillings of the $ r $-quotient of $ \lam $, \cite{MR1399504} or \cite[Lemma~2.1]{SXP}. With these two facts, we can rephrase \Cref{cor:prfixedcount} as follows.

\begin{Corollary}{\cite[Corollary~9.1]{MR2557880}} \label{cor:prfixedribbontab} For $ a, b \in \bZ_{\ge 1} $ and $ d \mid ab $,
\begin{equation} 
\label{eq:FPcount}
\begin{aligned}	
	\# \SYT(a^b)^{\pr^d} = \begin{cases} \# \SYT(Q_{\f{ab}{d}}(\lam)), & \quad \tx{if $ (a^b) $ has empty $ \f{ab}{d} $-core, } \\
							      0, & \quad \tx{else}
					\end{cases} 	
\end{aligned}
\end{equation}
where $ Q_{r}(\lam) $ is the $ r $-quotient of $ \lam $ combined anti-diagonally into a single skew shape.
\end{Corollary}

However, neither Rhoades's nor Purbhoo's proof describes which rectangular tableaux are fixed by $ d $ promotions.  The problem of constructing these fixed points is still open. We make substantial progress on this problem by characterizing all of the sufficiently large tableaux fixed by a given power of promotion. By \Cref{thm:prnid}, $ \SYT(a^b)^{\pr^d} $ is nonempty only when $ d  \mid ab $. Furthermore, we show in \Cref{sec:bkgd} that all nonempty cases are of the form $ \SYT((ar)^b)^{\pr^{br}} $ up to conjugation. Thus, it suffices to answer \Cref{qu:prfixedtabs}.

\begin{Question} \label{qu:prfixedtabs} For $ a, b, r \in \bZ_{\ge 1} $, which tableaux lie in $ \SYT((ar)^b)^{\pr^{br}} $?
\end{Question}

Some cases of \Cref{qu:prfixedtabs} have already been answered. The case $ a = 1 $ is a trivial consequence of  \Cref{thm:prnid}: $ \SYT(r^b)^{\pr^{br}} = \SYT(r^b) $. The case $ r = 1, a \ge b $ was solved by Kevin Purbhoo and Donguk Rhee \cite{MR3625918} in 2017. Their construction uses an algorithm similar to tableau stabilization. For each $ w = w_1 \dots w_b \in S_b $, they put $ w_1, \dots, w_b $ in cells placed anti-diagonally. Then they perform rectification while refilling the anti-diagonal cells with $ n $ plus whatever entry just left it. This algorithm agrees with stabilizing the anti-diagonal tableau and then restricting to cells left of this anti-diagonal. Finally, they perform the analogous algorithm with outer slides toward the southwest corner of $ (a^b) $ and attach the two results along the anti-diagonal to get a rectangular tableau of shape $ (a^b) $ fixed by $ \pr^b $.

The case $ a = 2 $ was solved by Dennis White, \cite{Whiteprom}, and independently by Donguk Rhee in his Master's thesis, \cite{RheeThesis} using the same construction. For the sake of completeness and so there is a record of this construction in the literature, we present their construction in \Cref{sec:otherprfixed}.

We answer \Cref{qu:prfixedtabs} for all $ b, r \in \bZ_{\ge 1} $ and $ a \ge 2b - 1 $, see \Cref{thm:prfixedtabs}. In addition, the tableaux we construct in $ \SYT((ar)^b)^{\pr^{br}} $ are in natural bijection with $ \SYT(Q_a((ar)^b)) $, which are in bijection with standard $ \f{ab}{d} $-ribbon tableaux of shape $ (a^b) $, as in \Cref{cor:prfixedribbontab}, solving \cite[Problem~9.4]{MR2557880} for $ a \ge 2b - 1 $. Describing the tableaux in $ \SYT((ar)^b)^{\pr^{br}} $ for $ a \in [3, 2b - 2] $ remains open, see \Cref{prob:otherprfixed}.

Let $ \Rect $ denote the rectification operator and $ \Rect^\ast $ denote the anti-rectification operator, which slides a skew tableau until it is right-justified, see \Cref{def:antirect}. For partitions $ \lam, \mu $, let $ \lam \cup \mu $ denote the result of combining $ \lam $ and $ \mu $ anti-diagonally into a single skew shape, see \Cref{ex:quotients}. By \Cref{cor:prfixedribbontab}, the tableaux in $ \SYT((ar)^b)^{\pr^{br}} $ are in bijection with $ \SYT( (r) \cup \dots \cup (r) ) $, with $ b $ pieces. For any tableau $ S \in \SYT( (r) \cup \dots \cup (r) ) $ and integer $ a \ge 2b - 1 $, let $ R_a(S) $ be the rectangular tableau formed by row-concatenating $ \Rect(S^{(b - 1)}), S^{(a - 2b + 2)} + (b -1)r $, and $  \Rect^\ast(S^{(b - 1)}) + (a - b + 1)r $ together from left to right.

\begin{Thm} \label{thm:prfixedtabs} For any $ b,r \in \bZ_{\ge 1} $ and integer $ a \ge 2b - 1 $, the tableaux in $ \SYT((ar)^b)^{\pr^{br}} $ are all constructed as follows:
\[
	\SYT((ar)^b)^{\pr^{br}} = \lb R_a(S) : S \in \SYT( \underbrace{(r) \cup \dots \cup (r)}_{b \tx{ times}} ) \rb.
\]
In addition,
\begin{align}
\label{eq:SYTquotientcount}
	\# \SYT((ar)^b)^{\pr^{br}} = \ch{br}{r, r, \dots, r}.
\end{align}

\end{Thm}

Moreover, since promotion commutes with itself, $ \SYT((ar)^b)^{\pr^{br}} $ is closed under the promotion operator. Using the extension of promotion to skew shapes, we will show that promotion commutes with the $ R_a $ operator.

\begin{Corollary} \label{cor:prfixedpr} For all $ a \ge 2b - 1 $ and $ S \in \SYT( \underbrace{(r) \cup \dots \cup (r)}_{b \tx{ times}}) $,
\[
	\pr(R_a(S)) = R_a(\pr(S)).
\] 
\end{Corollary}

The case $ a = 2 $ uses a similar construction, but it does not require tableau stabilization. Using \Cref{cor:prfixedribbontab}, the tableaux in $ \SYT((2r)^b)^{\pr^{br}} $ are in bijection with $ \SYT ((r^{\left \lceil \f{b}{2} \right \rceil}) \cup (r^{\left \lfloor \f{b}{2} \right \rfloor} )) $. For $ S \in \SYT ((r^{\left \lceil \f{b}{2} \right \rceil}) \cup (r^{\left \lfloor \f{b}{2} \right \rfloor} )) $. Let $ R_2(S) $ be formed by attaching $ \Rect(S) $ and $ \Rect^\ast(S) + br $ together from left to right. We show promotion commutes with the $ R_2 $ operator as well.

\begin{Thm} \cite{Whiteprom} \cite{RheeThesis} \label{thm:fixedbyab/2prom} For $ b,r \in \bZ_{\ge 1} $,
\[
	\SYT((2r)^b)^{\pr^{br}} = \lb R_2(S) : S \in \SYT( (r^{\left \lceil \f{b}{2} \right \rceil}) \cup (r^{\left \lfloor \f{b}{2} \right \rfloor} )) \rb.
\]
In addition,
\[
	\# \SYT((2r)^b)^{\pr^{br}} = \ch{br}{\left \lfloor \f{b}{2} \right \rfloor r} \# \SYT( r^{\left \lceil \f{b}{2} \right \rceil} ) \dd \# \SYT( r^{\left \lfloor \f{b}{2} \right \rfloor}).
\]
\end{Thm}

\begin{Corollary} \cite{Whiteprom} \cite{RheeThesis} \label{cor:prfixedpr2} For all $ S \in \SYT ((r^{\left \lceil \f{b}{2} \right \rceil}) \cup (r^{\left \lfloor \f{b}{2} \right \rfloor} )) $,
\[
	\pr(R_2(S)) = R_2(\pr(S)).
\] 
\end{Corollary}

In \Cref{sec:bkgd}, we review the necessary background. In \Cref{sec:stab}, we introduce tableau stabilization and prove some of its general properties, notably \Cref{thm:stabweaklydecrows}. In \Cref{sec:samesizerows}, we restrict our attention to stabilization on tableaux with rows of the same size and prove Theorems \ref{thm:stab} and \ref{thm:stabshape}. In \Cref{sec:antistab}, we introduce anti-stabilization and prove that it coincides with stabilization when both are defined. In \Cref{sec:prfixedtabs}, we employ tableau stabilization to explicitly construct $ \SYT((ar)^b)^{\pr^{br}} $ for $ a \ge 2b - 1 $, proving \Cref{thm:prfixedtabs}, and describe promotion on $ \SYT((ar)^b)^{\pr^{br}} $ for $ a \ge 2b - 1 $. In \Cref{sec:otherprfixed}, we explicitly construct $ \SYT((2r)^b)^{\pr^{br}} $, proving \Cref{thm:fixedbyab/2prom}, and describe the action of promotion on $ \SYT((2r)^b)^{\pr^{br}} $. In \Cref{sec:stabonperms}, we study stabilization as a permutation statistic. In \Cref{sec:open}, we present open problems.

\section{Background}
\label{sec:bkgd}

\subsection{Words and Tableaux}

A \emph{word} is a finite sequence of letters in $ \bZ_{\ge 1} $, the set of positive integers. For any sequence $ w $, let $ w_j $ refer to the $ j $-th letter of $ w $ and $ \ell(w) $ refer to the length of $ w $. A word $ w_1 w_2 \dots w_n $ is \emph{increasing} if $ w_1 \le w_2 \le \dots \le w_n $ and \emph{decreasing} if $ w_1 \ge w_2 \ge \dots \ge w_n $. A \emph{subsequence} of a word $ w_1 w_2 \dots w_n $ is a word of the form $ w_{i_1} w_{i_2} \dots w_{i_k} $ for some $ 1 \le i_1 < i_2 < \dots < i_k \le n $. Let $ [n] \coloneqq \{ 1, 2, \dots, n \} $ and $ [a,b] \coloneqq \{ x \in \bZ : a \le x \le b \} $. The \emph{descent set} of $ w = w_1 \dots w_n $ is 
\[
	\Des(w) = \{ i \in [n - 1] : w_i > w_{ i + 1} \}.
\]

A \emph{partition} $ \lam = (\lam_1, \dots, \lam_k) $ is a weakly decreasing sequence of positive integers. We say $ \lam \vdash n $ or $ |\lam| = n $ if $ \lam_1 + \dots + \lam_k = n $.  We view $ \lam $ as a left-justified diagram with $ \lam_1, \dots, \lam_k $ cells in its rows from top to bottom, using English notation. If $ \lam, \mu $ are partitions so that all of the cells in $ \mu $ are in $ \lam $, the \emph{skew shape} $ \lam/\mu $ is the set of cells in $ \lam $ but not in $ \mu $. 

A (semistandard) \emph{tableau} is a filling of a partition with entries from $ \bZ_{\ge 1} $ so that the rows are weakly increasing and the columns are strictly increasing. A (semistandard) skew tableau is such a filling of a skew shape. For a skew tableau $ S $, let $ |S| $ denote the number of cells in $ S $. A skew tableau $ S $ is \emph{standard} if it uses $ 1, \dots, |S| $ exactly once. For partitions $ \lam, \mu $, let $ \SYT( \lam ), \SYT(\lam / \mu) $ denote the set of standard skew tableaux of shape $ \lam, \lam/\mu $, respectively. We say a skew tableau has \emph{straight shape} if its shape is a partition. For any skew tableau $ S $, let $ \rv(S) $, the \emph{row vector} of $ S $, denote the sequence of row sizes from top to bottom so that $ \rv(S)_j $ denotes the number of cells in row $ j $ of $ S $. The row vector of $ S $ coincides with the shape of $ S $ when $ S $ has straight shape. The \emph{reading word} of $ S $ is the sequence of entries in $ S $ read left to right along the rows, from bottom to top. Let $ \lam', S' $ denote the conjugates of $ \lam, S $, respectively, obtained by interchanging the rows and columns. For a standard skew tableau $ S $, its descent set is 
\[
	\Des(S)= \{ i : (i + 1) \tx{ lies strictly below } i \tx{ in } S \}. 
\]

\begin{Example} \label{ex:tabstuff}
\[
	\Yvcentermath1 S = \begin{ytableau} \none & \none & \none & \none & \none & 1 & 7 \\ \none & \none & \none & \none & \none \\ \none & 2 & 4 & 6 & 9 \\ 3 & 5 & 8 \end{ytableau}   
\]
is a standard skew tableau with skew shape $ (7,5,5,3)/(5,5,1,0), |S| = 9, \rv(S) = (2,0,4,3), \rv(S)_3 = 4 $, reading word $ 358246917 $, and $ \Des(S) = \{ 1, 2, 4, 7 \} $.
\end{Example}

\subsection{Row Insertion and Jeu de Taquin}

This paper uses row insertion, the Robinson--Schensted--Knuth correspondence (RSK), and jeu de taquin heavily. We will state the known properties we will use here, but we assume the reader is already familiar with row insertion, RSK, and jeu de taquin. We recommend \cite[Chapter~3]{MR1824028} or \cite[Chapter~1]{MR1464693} for background on these algorithms and their properties.

\begin{Definition} For any tableau $ T $ and $ x \in \bZ_{\ge 1} $, let $ T \leftarrow x $ be the tableau obtained by row inserting $ x $ into $ T $. For a word $ w = w_1 w_2 \dots w_n $, let
\[
	T \leftarrow w \coloneqq ( \dots ((T \leftarrow w_1 ) \leftarrow w_2) \dots ) \leftarrow w_n.
\]
Let $ P(w) \coloneqq \varnothing \leftarrow w $ denote the insertion tableau of $ w $ and $ Q(w) $ denote the recording tableau of $ w $, which records the order the cells are created in $ \varnothing \leftarrow w $. Then, 
\begin{align}
\label{eq:DesQ}
	\Des(w) = \Des(Q(w)),
\end{align}
which is a consequence of \Cref{lem:bumpingcomparisonorig}.
\end{Definition}

\begin{Definition} The \emph{bumping chain} for the row insertion $ T \leftarrow x $ is the set of cells in $ T \leftarrow x $ which are bumped into or created while row inserting $ x $ into $ T $. 
\end{Definition}

\begin{Example} We highlight the bumping chain of this insertion in yellow:
\[
	\Yvcentermath1 \begin{ytableau} 1 & 3 & 5 & 7 & 8 \\ 2 & 9 & 10 & 11 & 14 \\ 4 & 12 & 15 \\ 13 \end{ytableau} \leftarrow 6 = \begin{ytableau} 1 & 3 & 5 & *(yellow) 6 & 8 \\ 2 & *(yellow) 7 & 10 & 11 & 14 \\ 4 & *(yellow) 9 & 15 \\ *(yellow) 12 \\ *(yellow) 13 \end{ytableau} \, . 
\]
\end{Example}


\begin{Lemma} \cite[\S 1.1]{MR1464693} \label{lem:bumpingcomparisonorig} Let $ T $ be any tableau, and consider the sequential row insertions $ (T \leftarrow x) \leftarrow y $.
\begin{enumerate}[(i)]
\item If $ x \le y $, the bumping chain of $ y $ is strictly to the right of the bumping chain of $ x $. 
\item If $ x > y $, the bumping chain of $ y $ is weakly to the left of the bumping chain of $ x $.
\end{enumerate}
\end{Lemma}

We will need the following generalization of \Cref{lem:bumpingcomparisonorig}(ii) that allows for intermediate insertions between $ x $ and $ y $. Because insertions can only decrease the value of each cell, \Cref{lem:bumpingcomparisonorig}(ii) still holds if the intermediate bumping chains avoid the bumping chain of $ x $. 

\begin{Lemma} \label{lem:bumpingcomparison} Let $ T $ be any tableau, and consider the sequential row insertions $ ( \dots ( (T \leftarrow x) \leftarrow x_1) \dots \leftarrow x_k) \leftarrow y $. If the bumping chains of $ x_1, \dots, x_k $ are disjoint from the bumping chain of $ x $ and $ x > y $, then the bumping chain of $ y $ is weakly to the left of the bumping chain of $ x $.
\end{Lemma}

%

The insertion tableau was originally defined to study increasing subsequences of words by Schensted, who proved that the length of the first row of $ P(w) $ is the maximum length of an increasing subsequence of $ w $ \cite{MR0121305}. Later, Greene generalized this result to describe $ \rv(P(w)) $ in terms of increasing subsequences of $ w $ \cite{MR0354395}. Since reversing a word $ w $  transposes the shape of $ P(w) $, analogous results hold for decreasing subsequences and the conjugate shape.  


\begin{Thm} \cite{MR0354395} \label{thm:Greene} Suppose $ w $ is a word and let $ \lam = (\lam_1, \dots, \lam_b) = \rv(P(w)) $. Then, for any $ k \le b $, $ \lam_1 + \dots + \lam_k $ is the maximum total length of $ k $ disjoint increasing subsequences of $ w $. If $ \lam' = (\lam_1', \dots, \lam_m') $, then for any $ k \le m $, $ \lam_1' + \dots + \lam_k' $ is the maximum total length of $ k $ disjoint decreasing subsequences of $ w $.

\end{Thm}

\begin{Definition} For any skew tableau $ S $, an \emph{inner slide} on $ S $ is the act of sliding into an outer corner of the inner shape of $ S $ and continuing to slide into the created hole so that increasing rows and columns are preserved until we reach an outer corner of $ S $. Let $ \Rect(S) $ denote the \emph{rectification} of $ S $, obtained by continually performing inner slides until straight shape is achieved. $ \Rect(S) $ is well-defined in that it is independent of the order of the slides. See \Cref{ex:rect}.

\end{Definition}

\begin{Example} \label{ex:rect} Consider
\[
	\Yvcentermath1 S = \begin{ytableau} \none & \none & 1 & 6 \\ \none & 3 & 4 & 9 \\ 2 & 7 & 8 & 11 \\ 5 & 10 & 12 & 13 \end{ytableau} \, .
\] 
We continually perform inner slides to the yellow $ \ast $ cell as follows. The green cells indicate which cells were just moved by the previous inner slide.
\begin{align*}
	\begin{ytableau} \none & *(yellow) \ast  & 1 & 6 \\ \none & 3 & 4 & 9 \\ 2 & 7 & 8 & 11 \\ 5 & 10 & 12 & 13 \end{ytableau}
\rightarrow \begin{ytableau} \none & *(green) 1 & *(green) 4 & 6 \\ *(yellow) \ast & 3 & *(green) 8 & 9 \\ 2 & 7 & *(green) 11 & *(green) 13 \\ 5 & 10 & 12 \end{ytableau}
\rightarrow	 \begin{ytableau} *(yellow) \ast & 1 & 4 & 6 \\ *(green) 2 & 3 & 8 & 9 \\ *(green) 5 & 7 & 11 & 13 \\ *(green) 10 & *(green) 12 \end{ytableau}
\rightarrow	 \begin{ytableau} *(green) 1 & *(green) 3 & 4 & 6 \\ 2 & *(green) 7 & 8 & 9 \\ 5 & *(green) 11 & *(green) 13 \\ 10 & 12 \end{ytableau} \, .
\end{align*}
Therefore,
\[
	\Yvcentermath1 \Rect(S) =   \begin{ytableau} 1 & 3 & 4 & 6 \\ 2 & 7 & 8 & 9 \\ 5 & 11 & 13 \\ 10 & 12 \end{ytableau} \, .
\]
\end{Example}

Rectification has many important roles in algebraic combinatorics. It is related to the RSK correspondence as in \Cref{lem:rectP}. It is involved in one of the many ways to define the Littlewood--Richardson coefficients, see \Cref{def:LRCs}. Inner and outer slides are also essential in defining promotion, demotion, and evacuation, see Definitions \ref{def:pr} and \ref{def:evac}.

\begin{Lemma} \cite[\S 2.1, Corollary 2]{MR1464693} \label{lem:rectP} For any skew tableau $ S $, 
\begin{align}
	\Rect(S) = P(w),
\end{align}
where $ w $ is the reading word of $ S $.
\end{Lemma}

\begin{Definition} \cite[Chapter 5.1]{MR1464693} \label{def:LRCs} Given any skew shape $ \lam/\mu $ and $ \nu \vdash |\lam| - |\mu| $, the corresponding \emph{Littlewood--Richardson coefficient} is 
\[
	c_{\mu, \nu}^{\lam} \coloneqq \# \{ S \in \SYT(\lam/ \mu) : \Rect(S) = T_0 \}
\]
for any $ T_0 \in \SYT(\nu) $. The number $ c_{\mu, \nu}^{\lam} $ is independent of the choice of $ T_0 \in \SYT(\nu) $ \cite[Corollary 5.1.1]{MR1464693}. 

\end{Definition}




Sch\"utzenberger's promotion operator on standard tableaux, see \Cref{def:pr}, appears in the context of representation theory. Promotion also generalizes to linear extensions of finite posets, including skew tableaux.


\begin{Definition} \label{def:pr} For any skew shape $ \lam/\mu $, \emph{promotion} $ \pr: \SYT(\lam/\mu) \to \SYT(\lam/\mu) $ is given by erasing the label $ 1 $, sliding that cell until it hits an outer corner of $ \lam/\mu $, filling that outer corner with $ |\lam| - |\mu| + 1 $, and finally decrementing all of the entries by 1. \emph{Demotion} is given by erasing the label $ n $, sliding  that cell until it hits an inner corner of $ \lam $, filling that inner corner with $ 0 $, and finally incrementing all of the entries by 1. Demotion and promotion are inverses when $ \lam/\mu $ is rectangular, but not in general. The names ``promotion" and ``demotion" for these operations are often interchanged, such as in \cite{MR2557880}.
\end{Definition} 

\begin{Example}
\[
	\Yvcentermath1 \pr: \begin{ytableau} 1 & 3 & 5 & 6 \\ 2 & 4 & 8 & 11 \\ 7 & 9 & 10 & 12 \end{ytableau} \mapsto \begin{ytableau} *(yellow) \ast & 3 & 5 & 6 \\ 2 & 4 & 8 & 11 \\ 7 & 9 & 10 & 12 \end{ytableau} \mapsto \begin{ytableau} 2 & 3 & 5 & 6 \\ 4 & 8 & 10 & 11 \\ 7 & 9 & 12 & *(yellow) \ast \end{ytableau} \mapsto \begin{ytableau} 1 & 2 & 4 & 5 \\ 3 & 7 & 9 & 10 \\ 6 & 8 & 11 & 12 \end{ytableau} \, .
\]
For an example on a skew shape,
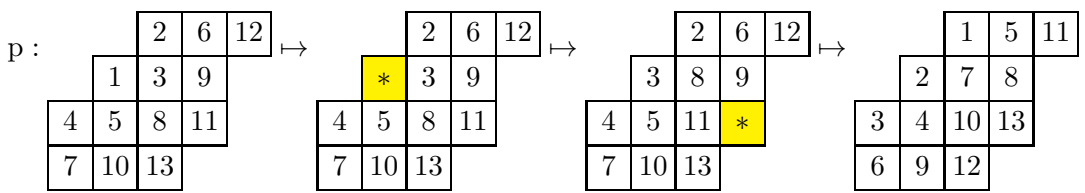
\[
	\Yvcentermath1 \pr: \begin{ytableau} \none & \none & 2 & 6 & 12 \\ \none & 1 & 3 & 9 \\ 4 & 5 & 8 & 11 \\ 7 & 10 & 13 \end{ytableau} \mapsto \begin{ytableau} \none & \none & 2 & 6 & 12 \\ \none & *(yellow) \ast & 3 & 9 \\ 4 & 5 & 8 & 11 \\ 7 & 10 & 13 \end{ytableau} \mapsto \begin{ytableau} \none & \none & 2 & 6 & 12 \\ \none & 3 & 8 & 9 \\ 4 & 5 & 11 & *(yellow) \ast \\ 7 & 10 & 13 \end{ytableau} \mapsto \begin{ytableau} \none & \none & 1 & 5 & 11 \\ \none & 2 & 7 & 8 \\ 3 & 4 & 10 & 13 \\ 6 & 9 & 12 \end{ytableau} \, .
\]
We will only use demotion on rectangular shapes, when it is the inverse of promotion.
\[
	\Yvcentermath1 \pr^{-1}: \begin{ytableau} 1 & 3 & 5 & 6 \\ 2 & 4 & 8 & 11 \\ 7 & 9 & 10 & 12 \end{ytableau} \mapsto \begin{ytableau} 1 & 3 & 5 & 6 \\ 2 & 4 & 8 & 11 \\ 7 & 9 & 10 & *(yellow) \ast \end{ytableau} \mapsto \begin{ytableau} *(yellow) \ast & 1 & 3 & 6 \\ 2 & 4 & 5 & 8 \\ 7 & 9 & 10 & 11 \end{ytableau} \mapsto \begin{ytableau} 1 & 2 & 4 & 7 \\ 3 & 5 & 6 & 9 \\ 8 & 10 & 11 & 12 \end{ytableau} \, .	
\]

\end{Example}

\bs


\begin{Definition} \label{def:evac} Suppose $ \lam \vdash n $. For $ T \in \SYT(\lam) $, the \emph{evacuation} of $ T $, $ e(T) \in \SYT(\lam) $ is the standard tableau that records the reverse order in which the cells are vacated as the smallest entry is repeatedly removed and the remaining skew tableau rectified.

\end{Definition}

\begin{Example} \label{ex:evac} Consider
\[
	T = \begin{ytableau} 1 & 4 & 5 \\ 2 & 6 \\ 3 \end{ytableau} \, .
\]
Repeatedly removing the smallest entry and rectifying, we get the sequence
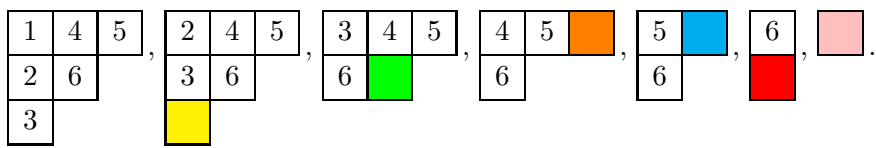
\[
	\begin{ytableau} 1 & 4 & 5 \\ 2 & 6 \\ 3 \end{ytableau} \, , \, \begin{ytableau} 2 & 4 & 5 \\ 3 & 6 \\ *(yellow) \,  \end{ytableau} \, , \, \begin{ytableau} 3 & 4 & 5 \\ 6 & *(green) \, \end{ytableau} \, , \, \begin{ytableau} 4 & 5 & *(orange) \\ 6 \end{ytableau} \, , \, \begin{ytableau} 5 & *(cyan) \\ 6 \end{ytableau} \, , \, \begin{ytableau} 6 \\ *(red) \,  \end{ytableau} \, , \,  \begin{ytableau} *(pink) \, \end{ytableau} \, .
\]
Recording the reverse order in which the cells are vacated gives
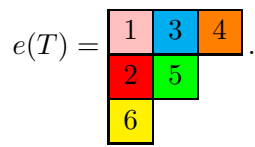
\[
	e(T) = \begin{ytableau} *(pink) 1 & *(cyan) 3 &  *(orange) 4 \\ *(red) 2 & *(green) 5 \\  *(yellow) 6 \end{ytableau} \, .
\]

\end{Example}

\begin{Definition} \label{def:dualequiv} Two standard skew tableaux $ S, T $ are \emph{dual equivalent} if they differ by a sequence of \emph{elementary dual equivalence moves}. The elementary dual equivalence move $ d_i $ swaps the cells of $ i \pm 1 $ and $ i $ if $ i \mp 1 $ appears between them in reading order:
\[
	d_i: \begin{matrix} & & & & \boxed{i} \\
			& & &  \reflectbox{$\ddots$} & \\
			 & & \boxed{i \mp 1} & & \\
			& \reflectbox{$\ddots$} & & &  \\
			\boxed{i \pm 1} & &&&
\end{matrix}
\leftrightarrow
\begin{matrix}  & & & & \boxed{i \pm 1} \\
			& & & \reflectbox{$\ddots$} & \\
			& & \boxed{i \mp 1} & & \\
			& \reflectbox{$\ddots$} & & &   \\
			\boxed{i} & & & & 
\end{matrix} \, .
\]
\end{Definition}

Haiman studies dual equivalence in \cite{MR1158783}. He defines dual equivalence as the property that the same sequence of slides produces the same shape and then proves that this property is characterized by \Cref{def:dualequiv}, see \cite[Theorem~2.6]{MR1158783}. Moreover, elementary dual equivalences moves commute with slides, see \cite[Lemma~2.3]{MR1158783}, so:
\begin{align}
\label{eq:DEcomm}
	d_i ( \Rect(S) ) = \Rect(d_i(S)), 
\end{align}
for all integers $ i $ and standard skew tableaux $ S $. Haiman uses dual equivalence to prove \Cref{thm:prnid}. 
Two permutations are dual equivalent if they differ by moves of the form
\[
	\dots (i \pm 1) \dots (i \mp 1) \dots i \dots \leftrightarrow \dots i \dots (i \mp 1) \dots (i \pm 1) \dots,
\]
for any integer $ i $. By \Cref{def:dualequiv}, dual equivalence of standard tableaux of the same skew shape corresponds to dual equivalence of their reading words \cite[Lemma~2.11]{MR1158783}. Dual equivalence classes of permutations are indexed by their common recording tableau: for all $ v, w \in S_n $,
\begin{align}
\label{eq:DEQ}
	Q(v) = Q(w) \iff \tx{ $ v $ and $ w $ are dual equivalent}. 
\end{align}

\ms

\subsection{Quotients and Cores of Partitions}

We review the definition of quotients and cores of partitions. We present a version due to James and Kerber \cite[Section 2.7]{MR644144} by viewing the boundary path of a partition as a binary sequence. Many other equivalent variations for the quotients and cores of partitions are known, such as \cite[I, Exercise 1.8]{MR1354144}. In our case, it will prove convenient to view the $ r $-quotient as a skew shape by combining its pieces anti-diagonally.

\begin{Definition} \label{def:quotients} Consider the map
\[
	\va: \{ \tx{Partitions} \} \to \{ \tx{Infinite binary strings with initial 0's and terminal 1's} \}
\] 
given by tracing the boundary path of the partition from southwest to northeast and writing 0 for each up step and 1 for each right step. The initial 0's represent the up steps along the negative $ y $-axis, and the terminal 1's represent the right steps along the positive $ x $-axis. We can thus view these sequences as finite binary strings where we can remove initial 0's and terminal 1's. 

For any partition $ \lam $ and $ r \in \bZ_{\ge 1} $, we form the \emph{$ r $-quotient}, $ Q_r(\lam) $, and \emph{$ r $-core}, $ C_r(\lam) $, of $ \lam $ as follows. Let $ w = \va(\lam) $, and for each $ j = 1, \dots, r $, let $ w^{(j)} $ denote the subsequence of $ w $ whose positions are congruent to $ j $ modulo $ r $. Let $ \lam^{(1)} = \va^{-1}(w^{(1)}), \dots, \lam^{(r)}= \va^{-1}(w^{(r)}) $. Then, $ Q_r(\lam) $ is the result of combining $ \lam^{(1)}, \dots, \lam^{(r)} $ into a single skew shape anti-diagonally, denoted $ \lam^{(1)} \cup \dots \cup \lam^{(r)} $, with $ \lam^{(1)} $ southwest of $ \lam^{(2)} $ southwest of $ \lam^{(3)} $, etc. Finally, let $ \ti{w} $ denote the binary sequence obtained after performing all swaps of the form 
\[
	1 x_1 \dots x_{r - 1} 0 \to 0 x_1 \dots x_{r - 1} 1
\]
on $ w $ until the 0's are are far left as possible. The order of swaps is irrelevant since we only swap two numbers whose indices are congruent modulo $ r $. Then, set $ C_r(\lam) = \va^{-1}(\ti{w}) $.   

\end{Definition}

\begin{Example} \label{ex:quotients} Suppose $ \lam = (7,5,5,5,3,2,1) $ and $ r = 3 $. Tracing the boundary path of $ \lam $ and labeling vertical steps with 0's and horizontal steps with 1's,

\begin{center}
\begin{tikzpicture}[scale = 0.8]

\coordinate (0) at (0,0) {};
\coordinate (1) at (1, 0) {};
\coordinate (1') at (1,1) {};
\coordinate (2) at (2,1) {};
\coordinate (3) at ( 2,2 ) {};
\coordinate (4) at ( 3,2 ) {};
\coordinate (5) at ( 3,3 ) {};
\coordinate (6) at ( 4,3 ) {};
\coordinate (7) at ( 5,3 ) {};
\coordinate (8) at ( 5,4 ) {};
\coordinate (9) at ( 5,5 ) {};
\coordinate (10) at ( 5,6 ) {};
\coordinate (11) at ( 6,6 ) {};
\coordinate (12) at ( 7,6 ) {};
\coordinate (13) at ( 7,7 ) {};
\coordinate (14) at ( 0,7 ) {};
\draw (0) -- (1) node [midway, below] {1} ;
\draw (1) -- (1') node [midway, left] {0}  ;
\draw (1') -- (2) node [midway, below] {1}  ;
\draw (2) -- (3) node [midway, left] {0}  ;
\draw (3) -- (4) node [midway, below] {1} ;
\draw (4) -- (5) node [midway, left] {0} ;
\draw (5) -- (6) node [midway, below] {1} ;
\draw (6) -- (7) node [midway, below] {1} ;
\draw (7) -- (8) node [midway, left] {0} ;
\draw (8) -- (9) node [midway, left] {0} ;
\draw (9) -- (10) node [midway, left] {0} ;
\draw (10) -- (11) node [midway, below] {1} ;
\draw (11) -- (12) node [midway, below] {1} ;
\draw (12) -- (13) node [midway, left] {0} ;
\draw (13) -- (14) ;
\draw (14) -- (0) ;
\end{tikzpicture}
\end{center}

\noindent Therefore,
\[
	w = \va(\lam) = 10101011000110,
\]
so
\[
	w^{(1)} = 10101 = 1010, \qquad w^{(2)} = 01100 = 1100, \qquad w^{(3)} = 1001 = 100,
\]
which correspond to 
\[
	\lam^{(1)} = \va^{-1}(w^{(1)}) = (2,1), \qquad \lam^{(2)} = \va^{-1}(w^{(2)}) = (2,2), \qquad \lam^{(3)} = \va^{-1}(w^{(3)}) = (1,1).
\]
Thus, 
\[
	\Yvcentermath1 Q_3(\lam) = (2,1) \cup (2,2) \cup (1,1) =  \young(::::\hfil,::::\hfil,::\hfil\hfil,::\hfil\hfil,\hfil\hfil,\hfil) = (5,5,4,4,2,1)/(4,4,2,2).
\]
Finally, performing all possible swaps of the form $ 1 x_1 x_2 0 \to 0 x_1 x_2 1 $ on $ w $ gives
\[
	\ti{w} = 00000010111111 = 10, 
\]
which corresponds to
\[
	C_3(\lam) = \va^{-1}(\ti{w}) = (1).
\]
\end{Example}

We need to know what rectangles have empty $ r $-core and what the $ r $-quotient is in that case. Then, we can apply \Cref{cor:prfixedribbontab} to find $ \# \SYT(a^b)^{\pr^d} $, and we will know when we have constructed all of the tableaux in $ \SYT(a^b)^{\pr^d} $. These results were known by Dennis White and are implicit in \cite{Whiteprom}, an unpublished work. We include their proofs here for completeness.

\begin{Lemma} \label{lem:rarb} Suppose $ a, b, r \in \bZ_{\ge 1} $ and $ r \mid ab $. Then, $ (a^b) $ has empty $ r $-core if and only if $ r \mid a $ or $ r \mid b $.

\end{Lemma}

\begin{proof} The binary word corresponding to the rectangle $ \lam = (a^b) $ is 
\[
	\va(\lam) = 1^a 0^b.
\]
From \Cref{def:quotients}, $ \lam $ has empty $ r $-core if and only if performing all possible swaps of the form $ 1 x_1 \dots x_{r - 1} 0 \rightarrow 0 x_1 \dots x_{r - 1} 1 $ to $ 1^a 0^b $ results in $ 0^b 1^a $. Hence, $ (a^b) $ has empty $ r $-core if and only if the positions of the zeros, namely $ a + 1, \dots, a + b $ are congruent to $ 1, 2, \dots, b $ modulo $ r $, in some order. We say two multisets $ \{ \be_1, \dots, \be_m \}, \{ \g_1, \dots, \g_m \} $ are congruent modulo $ r $, denoted $ \{ \be_1, \dots, \be_m \} \equiv_r \{ \g_1, \dots, \g_m \} $, if there exists a permutation $ w \in S_b $ so that $ \be_{w_j} \equiv_r \g_j $ for all $ j = 1, \dots, m $. Then, the $ r $-core of $ (a^b) $ is empty if and only if
\begin{align}
\label{eq:emptyrcore}
	\{ a + 1, \dots, a + b \} \equiv_r \{ 1, \dots, b \}.
\end{align}
By the division algorithm, we can write $ a = kr + i $ and $ b = \ell r + j $ with $ k, \ell \in \bZ_{\ge 0 } $ and $ i, j \in [0, r - 1] $. Then,
\begin{equation}
\label{eq:reducesetmodr} 
\begin{aligned}
	\{ a + 1, \dots, a + b \} & \equiv_r \{ i + 1, \dots, i + b \} \equiv_r \{ 0, \dots, r - 1 \}^{\ell} \cup \{ i + 1, \dots, i + j \} \\
	\{ 1, \dots, b \} & \equiv_r \{ 0, \dots, r - 1 \}^\ell \cup \{ 1, \dots, j \}.
\end{aligned}
\end{equation} 
Combining \eqref{eq:emptyrcore} and \eqref{eq:reducesetmodr} gives
\begin{align}
\label{eq:emptyrcoresimp}
	\{ i + 1, \dots, i + j \} \equiv_r \{ 1, \dots, j \}.
\end{align}
Since $ i, j \in [0, r - 1] $, \eqref{eq:emptyrcoresimp} holds if and only if $ i = 0 $ or $ j = 0 $, or equivalently $ r \mid a $ or $ r \mid b $.

\end{proof}

\begin{Remark} We can also characterize a partition $ \lam $ having empty $ r $-core using $ r $-ribbons. An $ r $-ribbon is a skew shape of size $ r $ which contains no two by two boxes. The $ r $-core of a partition $ \lam $ is the smallest partition $ \mu $ that an be obtained from $ \lam $ by successively removing $ r $-ribbons. From this perspective, it is clear that if $ r \mid a $ or $ r \mid b $, then the $ r $-core of $ (a^b) $ is empty. However, it is less obvious that the $ r $-core of $ (a^b) $ is empty only if $ r \mid a $ or $ r \mid b $.

\end{Remark}

\begin{Corollary} \label{cor:prfixednonempty} Suppose $ a, b, d \in \bZ_{\ge 1} $ and $ d \mid ab $. If $ \SYT(a^b)^{\pr^d} \ne \varnothing $, then
\begin{align*}
	\SYT(a^b)^{\pr^d} = \SYT((kr)^b)^{\pr^{br}}, \quad \tx{ or } \quad \SYT(a^b)^{\pr^d} = \SYT(a^{kr})^{\pr^{ar}}
\end{align*}
for some $ k, r \in \bZ_{\ge 1} $.
\end{Corollary}

\begin{proof}
By \Cref{cor:prfixedribbontab}, if $ \SYT(a^b)^{\pr^d} \ne \varnothing $, then $ (a^b) $ has empty $ \f{ab}{d} $-core. By \Cref{lem:rarb}, this means
\begin{align} 
\label{eq:prfixednonempty}
	\f{ab}{d} \mid a, \; \tx{ or } \; \f{ab}{d} \mid b \imp b \mid d, \tx{ or } a \mid d.
\end{align}
If $ b \mid d $, let $ r \coloneqq \f{d}{b} \in \bZ_{\ge 1} $, so $ d = br $ and
\[
	d \mid ab \imp r \mid a \imp a = kr \tx{ for some } k \in \bZ_{\ge 1}.
\]
Hence, $ \SYT(a^b)^{\pr^d} = \SYT((kr)^b)^{\pr^{br}} $. Otherwise, $ a \mid d $, so similarly,  $ \SYT(a^b)^{\pr^d} = \SYT((a^{kr})^{\pr^{ar}} $ for some $ k, r \in \bZ_{\ge 1} $.
\end{proof}

Thus, $ \SYT(a^b)^{\pr^d} $ is empty unless it is of the form $ \SYT((ar)^b)^{\pr^{br}} $ or $ \SYT(b^{ar})^{\pr^{br}} $ by \Cref{cor:prfixednonempty}. However, since promotion commutes with conjugation,
\[
	\SYT(b^{ar})^{\pr^{br}} = \{ T': T \in  \SYT((ar)^b)^{\pr^{br}} \}.
\]
Therefore, it suffices to consider to describe the tableaux in $ \SYT((ar)^b)^{\pr^{br}} $ to find $ \SYT(a^b)^{\pr^d} $ in general, as in \Cref{qu:prfixedtabs}.

\begin{Lemma} \label{lem:rectquotients} For $ a, b \in \bZ_{\ge 1} $, let $ s = b \; (\Mod a) \in [0, a - 1] $. Then, for any $ r \in \bZ_{\ge 1} $,
\[
	Q_a((ar)^b) = \underbrace{ (r^{\left \lceil \f{b}{a} \right \rceil}) \cup \dots \cup (r^{\left \lceil \f{b}{a} \right \rceil}) }_{s \tx{ times}} \cup \underbrace{ (r^{\left \lfloor \f{b}{a} \right \rfloor}) \cup \dots \cup (r^{\left \lfloor \f{b}{a} \right \rfloor}) }_{a - s \tx{ times}} 
\]
\end{Lemma}

\begin{proof} Following \Cref{def:quotients}, $ \nu \coloneqq (ar)^b $ corresponds to the binary sequence
\[
	\va(\nu) = 1^{ar} 0^b.
\]
Hence, $ Q_a(\nu) = \nu^{(1)} \cup \dots \cup \nu^{(a)} $ where 
\begin{align*}
	\nu^{(1)} = \dots = \nu^{(s)} = \va^{-1}(1^r \; 0^{ \left \lceil \f{b}{a} \right \rceil} ) = (r^{\left \lceil \f{b}{a} \right \rceil}), \\
	\nu^{(s + 1)} = \dots = \nu^{(a)} = \va^{-1}(1^r \; 0^{\left \lfloor \f{b}{a} \right \rfloor }) = (r^{\left \lfloor \f{b}{a} \right \rfloor}).
\end{align*}

\end{proof}

\subsection{Generalized Sums}

Later in \Cref{sec:samesizerows}, it will be convenient to allow sums $ \sum_{j = m}^n $ where $ m > n $. For a sequence $ a_0, a_1, \dots $, we generalize the notion $ \sum_{j = m}^n a_j $ from $ m \le n $ to all $ m, n \in \bZ_{\ge 0} $ as follows, as in Section 2.6 of \cite{ConcreteMath}.


\begin{Definition} \label{def:gensum} For $ m, n \in \bZ_{\ge 0} $, define
\begin{align}
\label{eq:gensum}
	\sum_{j = m}^n a_j \coloneqq \sum_{ j = 0}^{n} a_j - \sum_{ j = 0}^{m - 1} a_j,
\end{align}
where $ \sum_{ j = 0}^{-1} a_j = 0 $. In particular, for $ m > n $, we have
\begin{align}
\label{eq:oppdirsum}
	\sum_{j = m}^n a_j = \begin{cases} 0, & \tx{ if } m = n + 1,	\\
					- \sum_{j = n + 1}^{m - 1} a_j, &\tx{ if } m > n + 1.
\end{cases}
\end{align}
\end{Definition}

We make \Cref{def:gensum} so that
\begin{align}
\label{eq:sumprop}
	\sum_{j = m}^n a_j + \sum_{j = n + 1}^{p} a_j = \sum_{j = m}^p a_j
\end{align}
for all $ m, n, p \in \bZ_{\ge 0} $, which is an immediate consequence of \eqref{eq:gensum}. We introduce this notation so we can rewrite \eqref{eq:stabshape} more elegantly as
\begin{align}
\label{eq:stabshape2}
	\lam_j = kr + \sum_{i = j}^{b - j} c_i.	
\end{align}
In addition, we will make use of the following lemma involving $ \sum_{j = m}^n a_j $ with $ m > n $ in \Cref{sec:samesizerows}. 

\begin{Lemma} \label{lem:symsum} For any positive integers $ b \ge t \ge 1 $ and any sequence $ a_1, a_2, \dots $,
\[
	\sum_{ j = 1}^t \sum_{ i = j}^{b + j -  t - 1} a_i = \sum_{ j = 1}^t \sum_{ i = j}^{b - j} a_i. 
\]
\end{Lemma}

\begin{proof} We have
\begin{align*}
	\sum_{ j = 1}^t \sum_{i = j}^{b + j  - t - 1} a_i & = \sum_{ j = 1}^t \sum_{i = j}^{b - j} a_i + \sum_{ j = 1}^t \, \sum_{i = b - j + 1}^{b - t - 1} a_i + \sum_{ j = 1}^t \, \sum_{i = b - t}^{b + j - t - 1} a_i \qquad \tx{ by } \eqref{eq:sumprop} \\
	& = \sum_{ j = 1}^t \sum_{i = j}^{b - j} a_i + \sum_{ j = 1}^t \, \lp \sum_{i = b - j + 1}^{b - t - 1} a_i + \sum_{i = b - t}^{b - j} a_i \rp
\end{align*}
by re-indexing the last sum by $ j \mapsto t + 1 - j $.  Then,
 \begin{align*}
	\sum_{ j = 1}^t \sum_{i = j}^{b + j - t - 1} a_i & = \sum_{ j = 1}^t \sum_{i = j}^{b - j} a_i + \sum_{ j = 1}^t \, \sum_{i = b - j + 1}^{b - j} a_i \qquad \tx{ by } \eqref{eq:sumprop} \\
	& = \sum_{ j = 1}^t \sum_{i = j}^{b - j} a_i \qquad \tx{ by } \eqref{eq:oppdirsum}. 
\end{align*}
\end{proof}

\section{Tableau Stabilization}
\label{sec:stab}

In this section, we recall tableau stabilization from the introduction and explore its properties. After defining row shift equivalence and equivalent definitions of tableau stabilization, we show that stabilization is invariant under dual equivalence. Next, we show that stabilization continues after its first occurrence and prove an upper bound for stabilization when it is defined.

\begin{Definition} \label{def:rse} For any skew tableau $ S $, let $ S + x $ denote the result of adding $ x $ to each cell of $ S $. We say two skew tableaux $ S_1 $ and $ S_2 $ are \emph{row shift equivalent}, denoted $ \sim $, if there exists a constant $ x $ so $ S_1 + x $ and $ S_2 $ differ only by horizontal slides. For example,
\[
	\Yvcentermath1 \begin{ytableau} 1 & 2 & 4 \\ 3 & 5 \\ 6 \end{ytableau} \sim \begin{ytableau} \none & \none & 1 & 2 & 4 \\ \none & 3 & 5 \\ 6 \end{ytableau} \sim \begin{ytableau} \none & \none & 13 & 14 & 16 \\ \none & 15 & 17 \\ 18 \end{ytableau} \, . 
\]
\end{Definition}

Row shift equivalence is a natural equivalence relation for our purposes since stabilization is unaffected by horizontal slides. Using \Cref{def:rse}, we can rephrase the definition of tableau stabilization, see \Cref{def:stab}, in various ways as follows.

\begin{Definition} \label{def:stab2} For $ c \le d $, let $ \left. S \right|_{[c,d]} $ be the restriction of $ S $ to the cells with entries in $ [c,d] $. Suppose $ S $ is a standard skew tableau with $ m $ entries and decreasing row vector. Recall that $ S^{(k)} $ is the result of attaching $ (k - 1) $ shifted copies of $ S $ to the right of $ S $ so that the result is a standard skew tableau. This means $  \left. S^{(k)} \right|_{[(j - 1)m + 1, jm]} \sim S $ for all $ j = 1, \dots, k $. We say $ S $ \emph{stabilizes at $ k $} if any of the following equivalent conditions hold:

\ms

\begin{enumerate}[(a)]
\item $ \left. \Rect(S^{(k)}) \right|_{[(k - 1)m + 1, km]} \sim S $,
\ms
\item $  \left. \Rect(S^{(k)}) \right|_{[(k - 1)m + 1, km]} $ and $ S $ have the same number of cells in each row,
\ms
\item $ \rv(\Rect(S^{(k)})) - \rv(\Rect(S^{(k - 1)})) = \rv(S) $.
\end{enumerate}

\ms

\noindent Let $ \stab(S) $ denote the minimum value at which $ S $ stabilizes. See \Cref{ex:stab} for examples.

\end{Definition}

\Cref{def:stab2}(a), (b), (c) are all equivalent to \Cref{def:stab}. First, (a) is just a rephrasing of \Cref{def:stab}. Secondly, (a) clearly implies (b). Thirdly, (b) and (c) are equivalent because
\[
	\rv ( \left. \Rect(S^{(k)}) \right|_{[(k - 1)m + 1, km]} ) = \rv(\Rect(S^{(k)})) - \rv(\Rect(S^{(k - 1)})).
\]
Finally, if (c) holds, then since $ \left. \Rect(S^{(k)}) \right|_{[(k - 1)m + 1, km]} $ is obtained by sliding the cells in $ S + (k - 1)m $, the slides from $ S + (k -1)m $ to $  \left. \Rect(S^{(k)}) \right|_{[(k - 1)m + 1, km]} $ can only be horizontal. Hence, $ \left. \Rect(S^{(k)}) \right|_{[(k - 1)m + 1, km]} \sim S $.

\begin{Remark} \label{rem:rectcont} By definition of rectification, we have
\[
	\Rect(S^{(j)}) = \left. \Rect(S^{(j + 1)}) \right|_{[1, jm]} \tx{ for all } j \ge 1,
\]
so 
\[	
	\left. \Rect(S^{(j)}) \right|_{[(i - 1)m + 1, im]} = \left. \Rect(S^{(k)}) \right|_{[(i - 1)m + 1, im]} \tx{ for all } 1 \le i \le j \le k.
\]
Thus, to determine if $ S $ stabilizes at $ j $ for any $ j \le k $, it suffices to consider $ \Rect(S^{(k)}) $. 

\end{Remark}

\begin{Remark} \label{rem:rowsizes} We only consider standard skew tableaux with decreasing row vectors so that $ S^{(k)} $ is a standard skew tableau for all $ k \in \bZ_{\ge 1} $. If the row vector of $ S $ is not decreasing, then $ S^{(k)} $ need not be a skew tableau. For example,  
\[	
	\Yvcentermath1 S = \begin{ytableau} \none & 2 \\
				1 & 3 
		\end{ytableau}
	\imp S^{(2)} = \begin{ytableau} \none & 2 & *(yellow) 5 \\
				1 & 3 & *(yellow) 4 & *(yellow) 6
		\end{ytableau} \, ,
\]
which is not a skew tableau both because of its shape and the third column not being increasing.

\end{Remark}

\begin{Lemma} \label{lem:stabrse} Suppose $ S, \, T $ are standard skew tableaux with $ S \sim T $. Then, $ \Rect(S^{(k)}) = \Rect(T^{(k)}) $ for all $ k \in \bZ_{\ge 1} $ and $ \stab(S) = \stab(T) $.

\end{Lemma}

\begin{proof} Consider any $ k \in \bZ_{\ge 1} $. By definition of row shift equivalence restricted to standard skew tableaux, $ S^{(k)} $ and $ T^{(k)} $ have the same reading word, so $ \Rect(S^{(k)}) = \Rect(T^{(k)}) $ by \Cref{lem:rectP}. Let $ m = |S| = |T| $. Then, if $ S $ stabilizes at $ k $,
\[
	\left. \Rect(T^{(k)}) \right|_{[(k - 1)m, km]} = \left. \Rect(S^{(k)}) \right|_{[(k - 1)m, km]} \sim S \sim T, 
\]
so $ T $ stabilizes at $ k $ as well. Similarly the converse holds, and $ \stab(S) = \stab(T) $. 
\end{proof}

\begin{Thm} \label{thm:dualequivstab} If $ S, T $ are dual equivalent standard skew tableaux, then $ \stab(S) = \stab(T) $.

\end{Thm}

\begin{proof}
By \Cref{def:dualequiv}, we may assume $ T = d_i(S) $ without loss of generality. Shifting the entries by a constant preserves this relationship, so $ T + x = d_{i + x}(S + x) $ for all positive integers $ x $. It follows by the definition of $ S^{(k)} $ in \Cref{def:stab} that
\begin{align}
\label{eq:dualmovescopy}
	T^{(k)} = d_i \circ d_{i + m} \circ \dots \circ d_{i + (k - 1)m}(S^{(k)}) \tx{ for all } k \ge 1,
\end{align}
where $ m = |S| = |T| $, which means $ S^{(k)}, T^{(k)} $ are dual equivalent. By \eqref{eq:DEcomm}, \eqref{eq:dualmovescopy} implies
\begin{align}
\label{eq:dualmovesstab}
	\Rect(T^{(k)}) = d_i \circ d_{i + m} \circ \dots \circ d_{i + (k - 1)m}( \Rect(S^{(k)}))  \tx{ for all } k \ge 1.
\end{align} 
Since elementary dual equivalence moves preserve the row vector,
\[
	\rv(\Rect(T^{(k)})) = \rv(\Rect(S^{(k)})) \tx{ for all } k \ge 1.
\]
Hence, $ \stab(S) = \stab(T) $ by \Cref{def:stab2}(c).
\end{proof}

\begin{Example} \label{ex:dualmovesstab} Consider
\[
	S = \begin{ytableau} \none & \none & \none & 2 \\ \none & \none & 3 \\ \none & 1 \\ 4 \end{ytableau} \, , \qquad T = \begin{ytableau} \none & \none & \none & 1 \\ \none & \none & 3 \\ \none & 2 \\ 4 \end{ytableau} \, ,
\]
which satisfy $ T = d_2(S) $. Then, 
\[
	\Rect(S^{(3)}) = \begin{ytableau} 1 & 2 & *(yellow) 6 & *(green) 10 \\ 3 & *(yellow) 5 & *(yellow) 7 & *(green) 11 \\ 4 & *(green) 9 \\ *(yellow) 8 & *(green) 12 \end{ytableau} \, , \qquad \Rect(T^{(3)}) = \begin{ytableau} 1 & 3 &  *(yellow) 5 & *(green) 9 \\ 2 &  *(yellow) 6 &  *(yellow) 7 & *(green) 11 \\ 4 & *(green) 10 \\  *(yellow) 8 & *(green) 12 \end{ytableau} \, , 
\]
which satisfy $ \Rect(T^{(3)}) = d_2 \circ d_6 \circ d_{10}( \Rect(S^{(3)}) ) $. Note that $ d_2 $ swaps 1 and 2 in $ S $ and $ T $ but swaps 2 and 3 in $ \Rect(S^{(3)}) $ and $ \Rect(T^{(3)}) $. Here, $ \stab(S) = \stab(T) = 3 $, which is consistent with \Cref{thm:dualequivstab}.
\end{Example}

In order to show any standard skew tableau $ S $ with decreasing row vector eventually stabilizes and then continues to stabilize, consider $ S^{(\8)} $, obtained by attaching infinitely many shifted copies of $ S $ to the right of $ S $ so that the result uses $ 1, 2, \dots $ exactly once. As the row vector of $ S $ is decreasing, the rows and columns of $ S^{(\8)} $ are increasing. We can rectify $ S^{(\8)} $ with inner slides just like any skew tableau using the following lemma.

\begin{Lemma} \label{lem:eventuallyhoriz} Suppose $ S $ is a standard skew tableau with $ m $ entries and decreasing row vector. When rectifying $ S^{(\8)} $, if all inner slides so far pass horizontally through the entries in $ [(k - 1)m + 1, km] $ for some $ k \in \bZ_{\ge 1} $, then they pass horizontally through all entries greater than $ km $.

\end{Lemma}
 
\begin{proof} We localize our attention to the entries in $ [(k - 1)m + 1, km] $ and suppose all inner slides so far have proceed horizontally through $ [(k - 1)m + 1, km] $. Consider one such slide:
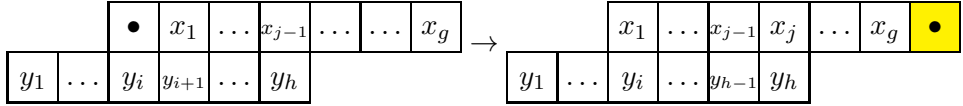
\[
	\ytableausetup{boxsize = 1.7em} \begin{ytableau} \none & \none & \bullet & x_1 & \dots & \scriptstyle{x_{j - 1}} & \dots & \dots &  x_{g} \\ y_1 & \dots & y_i & \scriptstyle{y_{i + 1}} & \dots & y_{h} \end{ytableau} 
	\rightarrow \begin{ytableau} \none & \none & x_1 & \dots & \scriptstyle{x_{j - 1}} &  x_j & \dots & x_{g} & *(yellow) \bullet \\ y_1 & \dots & y_i & \dots & \scriptstyle{y_{h - 1}} & y_{h}
\end{ytableau}
\]
where $ g \ge h $ and $ j = h  - i + 1 $. This means 
\[
	x_\ell < y_{\ell + i - 1} \qquad \tx{ for all } \ell = 1, \dots, h - i + 1. 
\]
Then, if we include the relevant entries in $ [km + 1, (k + 1)m] $, we have
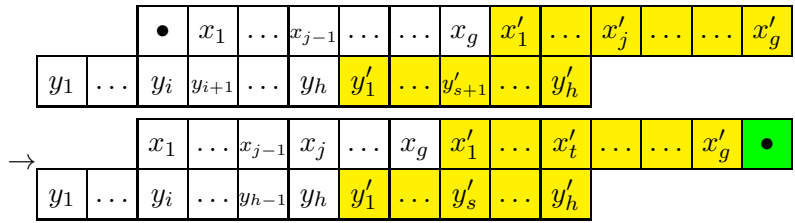
\begin{align*}
	&  \ytableausetup{boxsize = 1.7em} \begin{ytableau} \none & \none & \bullet & x_1 & \dots & \scriptstyle{x_{j -1}} & \dots & \dots & x_g &  *(yellow) x_1' &  *(yellow) \dots &  *(yellow) x_j' &  *(yellow) \dots & *(yellow) \dots &  *(yellow) x_g' \\ y_1 & \dots & y_i & \scriptstyle{y_{i + 1}} & \dots & y_h &  *(yellow) y_1' &  *(yellow) \dots &  *(yellow) \scriptstyle{y_{s + 1}'} &  *(yellow) \dots &  *(yellow) y_{h}' \end{ytableau} \\
	\rightarrow & \begin{ytableau} \none & \none & x_1 & \dots & \scriptstyle{x_{j - 1}} & x_j & \dots & x_g &  *(yellow) x_1'  &  *(yellow) \dots &  *(yellow) x_{t}' &  *(yellow) \dots & *(yellow) \dots &  *(yellow) x_g' & *(green) \bullet \\ y_1 & \dots & y_i & \dots & \scriptstyle{y_{h - 1}} & y_{h} &  *(yellow) y_1' &  *(yellow) \dots &  *(yellow) y_s' &  *(yellow) \dots &  *(yellow) y_h' \end{ytableau}
\end{align*}
where $ x_\ell' = x_\ell + m, y_\ell' = y_\ell + m, s = i + g - h \ge i , t = h - s + 1 \le j $. This slide also proceeds horizontally through $  [km + 1, (k + 1)m]  $ because
\[
	x_\ell' < y_{\ell + i - 1}' \le y_{\ell + s - 1}' \quad \tx{ for all } \; \ell = 1, \dots, h - s + 1. 
\]
Also, any cell whose diagonally southwest neighbor is in a different shifted copy of $ S $ slides horizontally. It follows inductively that this slide and thus all slides so far proceed horizontally through all entries greater than $ km $.

\end{proof} 

\begin{Lemma} \label{lem:stabcontinues} Suppose $ S $ is a standard skew tableau with decreasing row vector. Then, $ S $ must stabilize eventually, and if $ S $ stabilizes at $ k $, then $ S $ stabilizes at $ k' $ for all $ k' \ge k $.
\end{Lemma}

\begin{proof} Suppose $ S $ has $ b $ rows and skew shape $ \lam/\mu $. Then $ S^{(\8)} $ can be rectified with $ |\mu| $ inner slides, with $ \mu_i $ of them starting in row $ i $. Each inner slide starting in row $ i $ has at most $ b - i $ vertical slides. By \Cref{lem:eventuallyhoriz}, these $ b - i $ vertical slides must happen before or within the first $ b - i $ shifted copies of $ S $ that have not experienced vertical slides. Hence, all vertical slides take place in the first $ \sum_{i = 1}^{b} \mu_i (b - i) $ shifted copies of $ S $. This means $ S $ stabilizes at $ \sum_{i = 1}^{b} \mu_i (b - i)  + 1 $.

Suppose $ S $ stabilizes at $ k $. Thus, the entries in $ [(k - 1)m + 1, km] $ have not experienced vertical slides, so all entries greater than $ km $ have not experienced vertical slides by \Cref{lem:eventuallyhoriz}. Therefore, $ S $ stabilizes at $ k' $ for all $ k' \ge k $.

\end{proof}


While \Cref{lem:stabcontinues} shows that any standard skew tableau with decreasing row vector must stabilize eventually, its bound for stabilization, $ \sum_{i = 1}^{b} \mu_i (b - i)   + 1 $, is very weak. This naive bound depends on the inner shape, whereas our general bound, see \Cref{thm:stabweaklydecrows}, only depends on the number of rows $ b $, and is linear in $ b $. Any standard tableau $ T $ of straight shape has $ \stab(T) = 1 $  trivially since $ \Rect(T) = T $. Any nonempty standard skew tableau $ S $ with 1 row also has $ \stab(S) = 1 $ trivially. \Cref{thm:stabweaklydecrows} gives a general bound, $ 2b - 2 $, for tableaux with $ b $ rows ($b \ge 2 $).


\ms

\begin{Definition} \label{def:shiftcopy} For any word $ u = u_1 u_2 \dots u_n $ and $ m \in \bZ_{\ge 1} $, let
\[
	u + m = (u_1 + m) (u_2 + m) \dots (u_n + m).
\]
In addition, for $ k \in \bZ_{\ge 1} $, let
\[
	u^{(k,m)} = u (u + m) ( u + 2m) \dots (u+ (k - 1)m)
\]
If $ m $ is implicit, let $ u^{(k)} = u^{(k,m)} $. For example, if $ u = 134, k = 2, $ and $ m = 5 $, then
\[
	(134)^{(2)} = (134)^{(2,5)} = 134689.
\]
\end{Definition}

\bs

\begin{Notation} \label{not:tabstab}
Suppose $ S $ is a standard skew tableau with decreasing row vector $ (r_b, r_{b - 1}, \dots, r_2, r_1) $, so that its size is $ m \coloneqq r_1 + \dots + r_b $. Let $ w_1, w_2, \dots, w_b $ denote the entries in $ S $ read from left to right in rows $ b, (b - 1), \dots, 1 $, respectively, so that the reading word of $ S $ is $ w_1 w_2 \dots w_b $. For the rest of this paper, $ m $ will be implicit, so $ u^{(k)} = u^{(k,m)} $ for any word $ u $. Thus, the reading word of $ S^{(k)} $ is $ w_1^{(k)} w_2^{(k)} \dots w_b^{(k)} $. Hence, by \Cref{lem:rectP},
\begin{align}
	\Rect(S^{(k)}) = P(w_1^{(k)} w_2^{(k)} \dots w_b^{(k)}).
\end{align}
For $ j = 1, \dots, b $ and $ k \ge 1 $, let 
\begin{align} 
	T_j^{(k)}(S) \coloneqq P(w_1^{(k)} w_2^{(k)} \dots w_j^{(k)}).
\end{align}
so that $ T_b^{(k)}(S) = \Rect(S^{(k)}). $
\end{Notation}

\begin{Example} If
\[
	S = \byt \none & \none & \none & 4 & 5 & 6 \\ \none & \none & 3 & 7 \\ 1 & 2 \eyt \, ,
\]
then $ r_1 = 2, r_2 = 2, r_3 = 3, m = 7, w_1 = 12, w_2 = 37, w_3 = 456 $. Also,
\[
	 T_2^{(3)} = P( 1 2 8 9 (15)(16) 3 7 (10)(14)(17)(21) ) = \byt 1 & 2 & 3 & 7 & *(yellow) 10 &*(yellow) 14 & *(green) 17 & *(green) 21 \\ *(yellow) 8 & *(yellow) 9 & *(green) 15 & *(green) 16 \eyt \,. 
\]
Numbers bigger than 9 in a sequence are parenthesized to avoid confusion.
\end{Example}

\bs

Our goal is to understand the tails of the rows in $ T_b^{(k)}(S) $ for sufficiently large $ k $, which will turn out to $ k \ge 2b - 2 $, see \Cref{lem:formofTj}. In order to understand what happens to the tails of the rows under these row insertions, see \Cref{lem:consecshiftbumping}, we need to understand various comparisons between the elements we are inserting, see \Cref{lem:uvcomparisonshift}.

\begin{Lemma} \label{lem:uvcomparisonshift} Suppose $ u, v $ are increasing words on $ [m] $ with respective lengths $ r \le s $. Recall $ P(uv) $ denotes the insertion tableau of the concatenated word $ uv $, and let
\[
	c \coloneqq \rv(P(uv))_1 - s.
\]
Then, for each $ k \ge 1 $, 
\begin{align}
	v_{t}^{(k)} < u_{t + c}^{(k)} \tx{ for all } t \in [rk - c].
\end{align}
\end{Lemma}

\begin{proof} By \Cref{thm:Greene}, $ \rv(P(uv))_1 = c + s $ is the maximum length of an increasing subsequence of $ uv $. For any $ j \in [r - c] $, the subsequence
\[
	u_1 u_2 \dots u_{j + c} v_{j} v_{j + 1} \dots v_{s}
\]
of $ uv $ has size $ c + s + 1 $, so it is not increasing. As $ u, v $ are increasing, this forces
\begin{align}
\label{eq:uvineq}
	v_j < u_{j + c} \tx{ for all } j \in [r - c].
\end{align} 
Now,  by \Cref{def:shiftcopy},
\begin{align}
\label{uvshiftcopy}
	u_{ir + j}^{(k)} = u_j + (i - 1)m, \qquad v_{is + j'}^{(k)} = v_{j'} + (i - 1)m 
\end{align}
for all $ i \in [0, k - 1], j \in [r], j' \in [s] $. If $ i \in [0, k - 1], j \in [r- c] $, then $ j + c \le r $, so
\begin{align}
\label{eq:compshift1}
	v_{is + j}^{(k)} = v_{j} + (i - 1)m < u_{j + c} + (i - 1)m = u_{ir + j + c}^{(k)} \le u_{is + j + c}^{(k)}
\end{align} 
using \eqref{eq:uvineq}, $ r \le s $, and that $ u $ is increasing. On the other hand, if $ i \in [0, k - 2], j \in [r - c + 1, s] $, then $ v_j \le m $ and $ j + c \ge r + 1 $, so
\begin{align}
\label{eq:compshift2}
	v_{is + j}^{(k)} = v_{j} + (i - 1)m \le im < u_{(i + 1)r + 1}^{(k)} \le u_{ir + j + c}^{(k)} \le u_{is + j + c}^{(k)},
\end{align}
also using $ r \le s $ and that $ u $ is increasing. Combining \eqref{eq:compshift1} and \eqref{eq:compshift2} gives \Cref{lem:uvcomparisonshift}.
\end{proof}

\begin{Lemma} \label{lem:consecshiftbumping} Suppose $ i, m \in \bZ_{\ge 1} $, $ u, v $ are increasing words on $ [(i - 1)m + 1, im] $ with respective lengths $ r \le s $, and $ w, w' $ are increasing words on $ [(i  - 1)m] $. Then, for all $ k \ge 2 $,
\begin{align*}
	\boxed{ w u^{(k)} } \leftarrow w' v^{(k)} = \begin{tabular}{ | c | c | }
\hline
$ x $ & $ (v + m)^{(k - 1)} $ \\ 
\hline
$ x' $ & $ (u + m)^{(k - 1)}  $ \\
\hline
\end{tabular}
\end{align*}
where $ x, x' $ are increasing words on $ [im] $, possibly with different lengths. 


\end{Lemma}

\begin{proof} First, consider the insertion $ w u^{(k)} \leftarrow w' $. Since each element of $ w' $ is $ \le (i - 1)m $ and all elements of $ u^{(k)} $ are $ \ge (i - 1)m + 1 $, the elements of $ u^{(k)} $ that $ w' $ bumps down must form a consecutive sequence $ w'' $ starting from the beginning of $ u^{(k)} $, possibly empty. Letting $ c \coloneqq \rv(P(uv))_1 - s $, $ u \leftarrow v $ creates a tableau whose first row has size $ c + s $, meaning $ v $ bumps down $ r - c $ elements from $ u $. Hence, in $ w u^{(k)} \leftarrow w' v $, after $ w' $ bumps down $ w'' $ from $ u^{(k)} $, $ v $ bumps down at most $ r - c $ elements from $ u $. Thus, as each element of $ v $ is $ \le im $, $ v $ also bumps down at least $ c + s - r $ of the initial elements that are after both $ w'' $ and $ u $. In particular, as $ r \le s $, the $ c $ initial elements of $ (u + m)^{(k - 1)} $ are bumped down in $ w u^{(k)} \leftarrow w' v $.

Let $ c' \ge c $ be the number of initial elements of $ (u + m)^{(k - 1)} $ that are bumped down in $ w u^{(k)} \leftarrow w' v $. By \Cref{lem:uvcomparisonshift} and $ u $ being increasing, 
\[
	(v + m)^{(k - 1)}_t < (u + m)^{(k - 1)}_{t + c} \le (u + m)^{(k - 1)}_{t + c'} \tx{ for all } t \in [r(k - 1) - c']. 
\]
It follows inductively that the insertion of $ (v + m)^{(k - 1)} $ bumps down all of the remaining elements of $ (u + m)^{(k - 1)} $ in row 1. The result follows.

\end{proof}

\begin{Example} Suppose $ u = 134, v = 256 $, $ k = 3 $ and $ m = 6 $. Then, $ \boxed{u^{(3)}} \leftarrow v^{(3)} $
\begin{align*}
	& =  \ytableausetup{boxsize = 1.5em} \begin{ytableau} 1 & 3 & 4 & *(yellow) 7 & *(yellow) 9 & *(yellow) 10 & *(green) 13 & *(green) 15 & *(green) 16 \end{ytableau} \leftarrow   2 5 6 8(11)(12)(14)(17)(18) \\ 
		& = \begin{ytableau} 1 & 2 & 4 & 5 & 6 & *(yellow) 8 & *(yellow) 11 & *(yellow) 12 & *(green) 14 & *(green) 17 & *(green) 18 \\ 3 & *(yellow) 7 & *(yellow) 9 & *(yellow) 10 & *(green) 13 & *(green) 15 & *(green) 16 \end{ytableau} \\
		& = \begin{tabular}{ | c | c | }
\hline
12456 & $ (v + 6)^{(2)} $ \\ 
\hline
3 & $ (u + 6)^{(2)}  $ \\
\hline
\end{tabular}
\end{align*}
\end{Example}

\begin{Lemma} \label{lem:formofTj} For $ j = 2, \dots, b $, and $ k \ge 2b - 2 $, $ T_j ^{(k)}(S) $ is of the form

\begin{align*}
\begin{tabular}{ | c | c | }
\hline
$ x_1 $ & $ (w_j + (j  - 1)m)^{(k - j + 1)} $ \\ 
\hline
$ x_2 $ & $ (w_{j - 1} + j m)^{(k - j)}  $ \\
\hline
\vdots & \vdots \\
\hline
$ x_i $ & $ (w_{j + 1 - i} + (j -2 + i)m)^{(k - j + 2 - i)} $ \\
\hline
\vdots & \vdots \\
\hline
 $ x_{j - 1} $ & $ (w_{2} + (2j - 3)m)^{(k - 2j + 3)} $ \\
\hline
 $ x_{j} $ & $ (w_{1} + (2j - 3)m)^{(k - 2j + 3)} $ \\  
\hline
\end{tabular}
\end{align*}
where each $ x_{i} $ is an increasing word on $ [(j -2 + i)m] $, possibly with different lengths.
\end{Lemma}

\begin{proof} We proceed by induction on $ j $. The base case $ j = 1 $ holds since $ w_1^{(k)} $ is increasing. Inductively assuming the result holds for some fixed $ j \le b - 1 $, we have 
\begin{align}
	T_{j + 1}^{(k)}(S) = T_{j}^{(k)}(S) \leftarrow w_{j + 1}^{(k)}. 
\end{align}
By repeated applications of \Cref{lem:consecshiftbumping}, the bumping process in $ T_{j}^{(k)}(S) \leftarrow w_{j + 1}^{(k)} $ bumps $ (w_{j} + jm)^{(k - j)} $ from row 1 to row 2, which bumps $ (w_{j - 1} + (j + 1) m)^{(k - j - 1)} $ from row 2 to row 3, $ \dots $, which bumps $ (w_{j + 1 - i} + (j - 1 + i) m)^{(k - j + 1 - i)} $ from row $ i $ to row $ i + 1 $ for all $ i = 1, \dots, j $. This means $ (w_{j + 2 - i} + (j - 2 + i) m)^{(k - j + 2 - i)} $ bumps into row $ i $ for $ i = 1, \dots, j, j + 1 $. Thus, row $ i $ of $ T_{j + 1}^{k}(S) $ ends in $ (w_{j + 2 - i} + (j - 1 + i) m)^{(k - j + 1 - i)} $ for $ i = 1, \dots, j $ and in $ (w_{1} + (2j - 1) m)^{(k - 2j + 1)} $ in the newly created row $ j + 1 $. Smaller elements are also bumped down from each row, but they do not affect this pattern. This proves \Cref{lem:formofTj} for $ j + 1 $ and completes our induction.
\end{proof}
%
%
%

\begin{proof}[Proof of \Cref{thm:stabweaklydecrows}] Suppose $ S $ is a standard skew tableau with decreasing row vector and continue using \Cref{not:tabstab}. Plugging $ j = b $ and $ k = 2b - 2 $ into \Cref{lem:formofTj} gives 
\begin{align*}
	\Rect(S^{(2b - 2)}) = T_b^{(2b - 2)}(S) =
\begin{tabular}{ | c | c | }
\hline
$ x_1 $ & $ (w_b + (b - 1)m)^{(b - 1)} $ \\ 
\hline
$ x_2 $ & $ (w_{b - 1} + b m)^{(b - 2)}  $ \\
\hline
\vdots & \vdots \\
\hline
$ x_i $ & $ (w_{b + 1 - i} + (b -2 + i)m)^{(b - i)} $ \\
\hline
\vdots & \vdots \\
\hline
 $ x_{b - 1} $ & $ (w_{2} + (2b - 3)m)^{(1)} $ \\
\hline
 $ x_{b} $ & $ (w_{1} + (2b - 3)m)^{(1)} $ \\  
\hline
\end{tabular} \, ,
\end{align*}
where each $ x_i $ is an increasing word on $ [(b -2 + i)m] $. Hence, 
\begin{align*}
	\left. \Rect(S^{(2b - 2)}) \right|_{[(2b - 3)m + 1, (2b - 2)m]} \sim S,
\end{align*}
so $ S $ stabilizes at $ 2b - 2 $.

\end{proof}

\section{Stabilization on Skew Tableaux with Constant Row Vectors}
\label{sec:samesizerows}

In this section, we restrict our attention to skew tableaux with constant row vectors. We prove formula \eqref{eq:stabshape} for the shape of the stabilized tableau, see \Cref{thm:stabshape}. This allows us to improve our upper bound for stabilization from $ 2b - 2 $ to $ b $ in the constant row vector case, proving \Cref{thm:stab}.

\begin{Notation} \label{not:tabstab2} Suppose $ S $ is a standard skew tableau with constant row vector $ (r^b) $, so that it has $ b $ rows of size $ r $, and its size is $ m \coloneqq br $. We continue using \Cref{def:shiftcopy} and \Cref{not:tabstab} as before. Recall $ w_1, w_2, \dots, w_b $ are the reading words of the rows of $ S $ from bottom to top,
\begin{align} 
	T_j^{(k)}(S) \coloneqq P(w_1^{(k)} w_2^{(k)} \dots w_j^{(k)}) \quad \tx{ for } \quad j = 1, \dots, b \; \tx{ and } \; k \ge 1,
\end{align}
and $ T_b^{(k)}(S) = \Rect(S^{(k)}). $ In addition, let
\begin{align}
\label{eq:ci}
	c_i = c_i(S) \coloneqq \rv(P(w_i w_{i + 1}))_1 - r, \tx{ for all } i = 1, \dots, b - 1.
\end{align}
See \Cref{ex:incsubseqs} for an example computation of the $ c_i $'s.
\end{Notation}

First, we derive a lower bound on the partial sums of $ \rv(T_b^{(k)}(S)) $ using \Cref{thm:Greene}. Secondly, we will show that equality always holds for this bound using \Cref{lem:bumpingcomparisonorig} and \Cref{lem:bumpingcomparison}, which proves \Cref{thm:stabshape}. Recall from \eqref{eq:oppdirsum} that for all $ n, p \in \bZ_{\ge 0} $ with $ n > p $,
\begin{align*}
	\sum_{j = n}^p a_j = \begin{cases} 0, & \tx{ if } p = n - 1,	\\
					- \sum_{j = p + 1}^{n - 1} a_j, &\tx{ else}.
\end{cases}
\end{align*}

\begin{Lemma} \label{lem:shupperbound} Suppose $ S $ is standard skew tableau with row vector $ (r^b) $ and use \Cref{not:tabstab2}. If $ k \ge b - 1 $ and $ \rv(T_b^{(k)}(S)) = (\lam_1, \dots, \lam_b) $, then, for $ t = 1, \dots, b $,
\begin{align}
\label{eq:shupperbound}
	\sum_{j = 1}^t \lam_j \ge \sum_{j = 1}^t \lp kr +  \sum_{i = j}^{b - j} c_i \rp. 
\end{align}
\end{Lemma}

\begin{proof} Consider any positive integers $ t \le b $ and $ k \ge b - 1 $. By \Cref{thm:Greene}, we will be done if we show there exist $ t $ disjoint increasing subsequences of $ w(S,k) \coloneqq w_1^{(k)} \dots w_b^{(k)} $ whose total length is $ \sum_{j = 1}^t \lp kr +  \sum_{i = j}^{b - j} c_i \rp $. We will exhibit $ t $ such increasing subsequences of $ w(S,k) $ explicitly. Our exhibition of these subsequences uses the following $ b \x k $ matrix of words:
\[
	M_k(S) = \begin{matrix} w_b & w_b + m & \dots & w_b + (k - 1)m \\
				\vdots & \vdots & \vdots & \vdots \\
				w_2 & w_2 + m & \dots & w_2 + (k - 1)m \\
				w_1 & w_1 + m & \dots & w_1 + (k - 1)m \\
\end{matrix}
\]
Note that concatenating these words from left to right and then top to bottom gives $ w(S,k) $.



\begin{figure}[h]
\centering 
\begin{tikzpicture}[scale=1, transform shape]

\coordinate [label=left:$1$] (1a) at (0,0) {};
\coordinate (1b) at (4,0) {};
\coordinate [label=below:$t$] (1b') at (4,-0.15) {};
\coordinate (1c) at (4,0.5) {};
\coordinate (1d) at (4.5,0.5) {};
\coordinate (1e) at (4.5,1) {};
\coordinate (1f) at (5,1) {};
\coordinate (1g) at (6,2) {};
\coordinate (1h) at (6,2.5) {};
\coordinate (1i) at (8,2.5) {};
\coordinate (1j) at (8,2.5) {};

\coordinate [label=left:$2$] (2a) at (0,0.5) {};
\coordinate (2b) at (3.5,0.5) {};
\coordinate (2c) at (3.5,1) {};
\coordinate (2d) at (4,1) {};
\coordinate (2e) at (4,1.5) {};
\coordinate (2f) at (4.5,1.5) {};
\coordinate (2g) at (5.5, 2.5) {};
\coordinate (2h) at (5.5, 3) {};
\coordinate (2i) at (6, 3) {};
\coordinate (2j) at (8, 3) {};

\coordinate [label=left:$3$] (3a) at (0,1) {};
\coordinate (3b) at (3,1) {};
\coordinate (3c) at (3,1.5) {};
\coordinate (3d) at (3.5,1.5) {};
\coordinate (3e) at (3.5,2) {};
\coordinate (3f) at (4,2) {};
\coordinate (3g) at (5, 3) {};
\coordinate (3h) at (5, 3.5) {};
\coordinate (3i) at (6, 3.5) {};
\coordinate (3j) at (8,3.5) {};

\coordinate [label=left:$(t - 1)$] (4a) at (0,3.5) {};
\coordinate (4b) at (0.5,3.5) {};
\coordinate (4c) at (0.5, 4) {};
\coordinate (4d) at (1, 4) {};
\coordinate (4e) at (1,4.5) {};
\coordinate (4f) at (1.5, 4.5) {};
\coordinate (4g) at (2.5, 5.5) {};
\coordinate (4h) at (2.5, 6) {};
\coordinate (4i) at (6, 6) {};
\coordinate (4j) at (8, 6) {};

\coordinate [label=left:$t$] (5a) at (0,4) {};
\coordinate (5b) at (0,4) {};
\coordinate (5c) at (0,4.5) {};
\coordinate (5d) at (0.5, 4.5) {};
\coordinate (5e) at (0.5,5) {};
\coordinate (5f) at (1, 5) {};
\coordinate (5g) at (2, 6) {};
\coordinate [label=above:$b - t$] (5h) at (2, 6.5) {};
\coordinate [label=above:$b - 1$] (5i) at (6, 6.5) {};
\coordinate [label=above:$k$] (5j) at (8, 6.5) {};

\coordinate [label=right: $ b $] (6) at (8,6.5) {};
\coordinate [label=below: $1$] (7) at (0,-0.15) {};
\coordinate [label=right: $ (b - t + 1)$] (8) at (8,2.5) {};
\coordinate [label=right: $ (b - t + 2)$] (9) at (8,3) {};
\coordinate [label=right: $ (b - t + 3)$] (11) at (8,3.5) {};
\coordinate [label=right: $ (b - 1) $] (10) at (8,6) {};

\draw (1a) -- (1b);
\draw (1b) -- (1c);
\draw (1c) -- (1d);
\draw (1d) -- (1e);
\draw (1e) -- (1f);
\draw[style = dashed] (1f) -- (1g);
\draw (1g) -- (1h);
\draw (1h) -- (1i);
\draw (1i) -- (1j);

\draw (2a) -- (2b);
\draw (2b) -- (2c);
\draw (2c) -- (2d);
\draw (2d) -- (2e);
\draw (2e) -- (2f);
\draw[style = dashed] (2f) -- (2g);
\draw (2g) -- (2h);
\draw (2h) -- (2i);
\draw (2i) -- (2j);

\draw (3a) -- (3b);
\draw (3b) -- (3c);
\draw (3c) -- (3d);
\draw (3d) -- (3e);
\draw (3e) -- (3f);
\draw[style = dashed] (3f) -- (3g);
\draw (3g) -- (3h);
\draw (3h) -- (3i);
\draw (3i) -- (3j);

\draw (4a) -- (4b);
\draw (4b) -- (4c);
\draw (4c) -- (4d);
\draw (4d) -- (4e);
\draw (4e) -- (4f);
\draw[style = dashed] (4f) -- (4g);
\draw (4g) -- (4h);
\draw (4h) -- (4i);
\draw (4i) -- (4j);

\draw (5a) -- (5b);
\draw (5b) -- (5c);
\draw (5c) -- (5d);
\draw (5d) -- (5e);
\draw (5e) -- (5f);
\draw[style = dashed] (5f) -- (5g);
\draw (5g) -- (5h);
\draw (5h) -- (5i);
\draw (5i) -- (5j);

\draw [fill = black] (1a) circle (0.050);
\draw [fill = black] (1b) circle (0.050);
\draw [fill = black] (1c) circle (0.050);
\draw [fill = black] (1d) circle (0.050);
\draw [fill = black] (1e) circle (0.050);
\draw [fill = black] (1f) circle (0.050);
\draw [fill = black] (1g) circle (0.050);
\draw [fill = black] (1h) circle (0.050);
\draw [fill = black] (1i) circle (0.050);
\draw [fill = black] (1j) circle (0.050);

\draw [fill = black] (2a) circle (0.050);
\draw [fill = black] (2b) circle (0.050);
\draw [fill = black] (2c) circle (0.050);
\draw [fill = black] (2d) circle (0.050);
\draw [fill = black] (2e) circle (0.050);
\draw [fill = black] (2f) circle (0.050);
\draw [fill = black] (2g) circle (0.050);
\draw [fill = black] (2h) circle (0.050);
\draw [fill = black] (2i) circle (0.050);
\draw [fill = black] (2j) circle (0.050);

\draw [fill = black] (3a) circle (0.050);
\draw [fill = black] (3b) circle (0.050);
\draw [fill = black] (3c) circle (0.050);
\draw [fill = black] (3d) circle (0.050);
\draw [fill = black] (3e) circle (0.050);
\draw [fill = black] (3f) circle (0.050);
\draw [fill = black] (3g) circle (0.050);
\draw [fill = black] (3h) circle (0.050);
\draw [fill = black] (3i) circle (0.050);
\draw [fill = black] (3j) circle (0.050);

\draw [fill = black] (4a) circle (0.050);
\draw [fill = black] (4b) circle (0.050);
\draw [fill = black] (4c) circle (0.050);
\draw [fill = black] (4d) circle (0.050);
\draw [fill = black] (4e) circle (0.050);
\draw [fill = black] (4f) circle (0.050);
\draw [fill = black] (4g) circle (0.050);
\draw [fill = black] (4h) circle (0.050);
\draw [fill = black] (4i) circle (0.050);
\draw [fill = black] (4j) circle (0.050);

\draw [fill = black] (5a) circle (0.050);
\draw [fill = black] (5b) circle (0.050);
\draw [fill = black] (5c) circle (0.050);
\draw [fill = black] (5d) circle (0.050);
\draw [fill = black] (5e) circle (0.050);
\draw [fill = black] (5f) circle (0.050);
\draw [fill = black] (5g) circle (0.050);
\draw [fill = black] (5h) circle (0.050);
\draw [fill = black] (5i) circle (0.050);
\draw [fill = black] (5j) circle (0.050);

\draw [fill = black] (2,2) circle (0.015);
\draw [fill = black] (2.2,1.8) circle (0.015);
\draw [fill = black] (1.8,2.2) circle (0.015);
\draw [fill = black] (4,4) circle (0.015);
\draw [fill = black] (4.2,3.8) circle (0.015);
\draw [fill = black] (3.8,4.2) circle (0.015);

\end{tikzpicture}
\caption{The lattice paths used to construct increasing subsequences of $ w(S,k) $}
      \label{fig:LPs}
\end{figure}
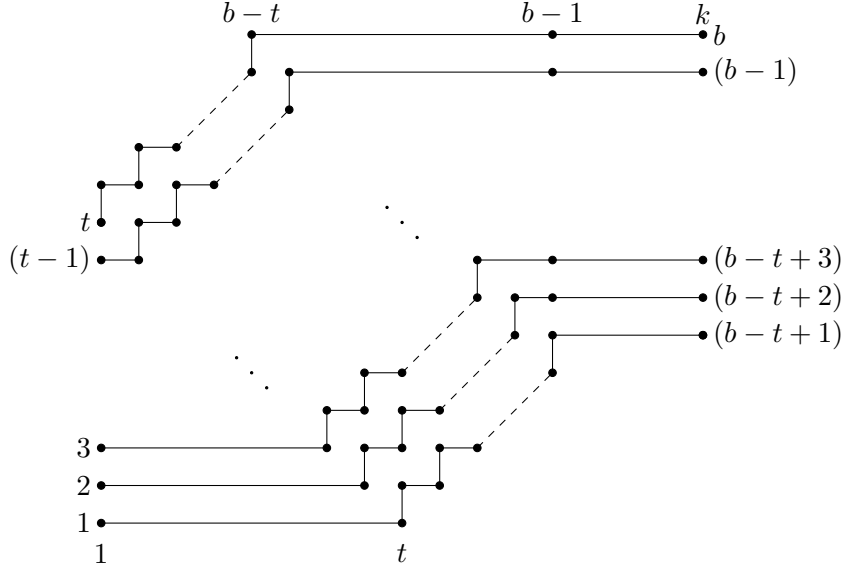

In addition, we will use the lattice paths on the Cartesian plane in \Cref{fig:LPs}. Label the paths $ L_1, \dots, L_t $ in \Cref{fig:LPs} from bottom to top. For $ j = 1, \dots, t $, the path $ L_j $ proceeds as follows:
\begin{enumerate}[(a)]
\item Walk horizontally from $ (1, j) $ to $ (t + 1 - j, j) $, then
\item Alternate between $ (0,1) $ and $ (1,0) $ steps, starting with $ (0,1) $, until reaching $ (b - j, b + j - t) $, then
 \item Walk horizontally to $ (k, b + j - t) $.
\end{enumerate}
Now, overlay the $ b \x k $ grid containing these lattice paths on $ M_k(S) $. We can then convert each lattice path $ L_j $ into a word $ y_j $ as follows. First concatenate the words in $ M_k(S) $ that correspond to lattice points in $ L_j $ in order starting from $ (1,j) $ and ending at $ (k, b + j - t) $ to obtain $ x_j $. Note $ x_j $ need not be increasing because it contains $ (w_i + sm) (w_{i + 1} + sm) $ in each column $ i $ where there is a vertical step. To obtain an increasing subsequence associated to $ L_j $, we will replace each $ (w_i + sm) (w_{i + 1} + sm) $ by a maximum length increasing subsequence of $ (w_i + sm) (w_{i + 1} + sm) $. By \Cref{thm:Greene} and the definition of $ c_i $, \eqref{eq:ci}, there exists a maximum length increasing subsequence $ u_i $ of $ w_i w_{i + 1} $ with length $ c_i + r $. Thus, replace each instance of $ (w_i + sm) (w_{i + 1} + sm) $ in $ x_j $ with $ u_i + sm $ to obtain the increasing subsequence $ y_j $. Hence,
\begin{align*}
	y_j  = w_j^{(t - j)} & (u_{j} + (t - j)m) (u_{j + 1} + (t + 1 - j)m) \dots (u_{b + j - t - 1} + (b - 1 - j)m) \\
	& (w_{b + j - t} + (b + 1 - j)m)^{(k - b + j)}.
\end{align*}

Since the lattice paths $ L_1, \dots, L_t $ are disjoint, $ y_1, \dots, y_t $ are disjoint as well. For $ j = 1, \dots, t $, as $ \ell(w_i) = r $ and $ \ell(u_i) = c_i + r $ for all $ i $, we have
\begin{align*}
	\ell(y_j) & = (t - j)r + \lp \sum_{i = j}^{b + j - t - 1} (c_i + r) \rp + (k - b + j)r = kr + \sum_{i = j}^{b + j - t - 1} c_i.
\end{align*}
Hence,
\begin{align*}
	\sum_{j = 1}^t \ell(y_j) & = tkr + \sum_{ j = 1}^t \sum_{i = j}^{b + j - t - 1} c_i = \sum_{ j = 1}^t \lp kr + \sum_{i = j}^{b - j} c_i \rp
\end{align*}
by \Cref{lem:symsum}. Thus, $ y_1, \dots, y_t $ are $ t $ disjoint increasing subsequences of $ w(S,k) $ with total length $ \sum_{j = 1}^t \lp kr + \sum_{i = j}^{b - j} c_i \rp $, completing the proof.
\end{proof}


\begin{Example} \label{ex:incsubseqs} Suppose
\[
	\Yvcentermath1 S = \begin{ytableau} \none & \none & \none & \none & \none & \none & 7 & 9 & 10 \\
				    \none & \none & \none & \none & 2 & 4 & 12 \\
				    \none & \none & \none & 6 & 8 & 11 \\
				    1 & 3 & 5
\end{ytableau} \, .
\]
Thus, $ b = 4, r = 3, m = 12 $, $ w_1 = 135, w_2 = 6 8 (11), w_3 = 24(12), w_4 = 79(10)$, 
\begin{align*}
	P(w_1 w_2) & = \begin{ytableau}1 & 3 & 5 & 6 & 8 & 11 \end{ytableau} \, , 
	P(w_2 w_3) = \begin{ytableau} 2 & 4 & 11 & 12 \\ 6 & 8 \end{ytableau} \, , \\
	P(w_3 w_4) & = \begin{ytableau} 2 & 4 & 7 & 9 & 10 \\ 12 \end{ytableau} \, ,
\end{align*}
so $ c_1 = 3, c_2 = 1, c_3 = 2 $. Maximum length increasing subsequences of $ w_1 w_2, w_2 w_3 $ and $ w_3 w_4 $ are
\[
	u_1 = 13568(11), \quad u_2 = 68(11)(12), \quad u_3 = 2 4 7 9 (10),
\]
respectively. The matrix of words for this example is
\[
	M_3(S) = \begin{matrix} w_4 & w_4 + 12 & w_4 + 24 \\
			w_3 & w_3 + 12 & w_3 + 24 \\
			w_2 & w_2 + 12 & w_2 + 24 \\
			w_1 & w_1 + 12 & w_1 + 24 \end{matrix} \quad
	= \quad \begin{matrix} 79(10) & (19)(21)(22) & (31)(33)(34) \\
			24(12) & (14)(16)(24) & (26)(28)(36) \\
			68(11) & (18)(20)(23) & (30)(32)(35) \\
			135 & (13)(15)(17) & (25)(27)(29) \end{matrix} \, .
\]
The lattice paths from the proof of \Cref{lem:shupperbound} for $ b = 4, k = 3 $ and $ t = 1, 2, 3, 4 $ are shown in \Cref{fig:exLPs} from left to right, respectively. 

\begin{figure}[h]
\centering
\begin{tikzpicture}[scale=1, transform shape]

\coordinate (1a) at (0,0) {};
\coordinate (1b) at (1,0) {};
\coordinate (1c) at (2,0) {};
\coordinate (2a) at (0,1) {};
\coordinate (2b) at (1,1) {};
\coordinate (2c) at (2,1) {};
\coordinate (3a) at (0,2) {};
\coordinate (3b) at (1,2) {};
\coordinate (3c) at (2,2) {};
\coordinate (4a) at (0,3) {};
\coordinate (4b) at (1,3) {};
\coordinate (4c) at (2,3) {};

\draw (1a)--(2a)--(2b)--(3b)--(3c)--(4c);

\draw [fill = black] (1a) circle (0.050);
\draw [fill = black] (1b) circle (0.050);
\draw [fill = black] (1c) circle (0.050);
\draw [fill = black] (2a) circle (0.050);
\draw [fill = black] (2b) circle (0.050);
\draw [fill = black] (2c) circle (0.050);
\draw [fill = black] (3a) circle (0.050);
\draw [fill = black] (3b) circle (0.050);
\draw [fill = black] (3c) circle (0.050);
\draw [fill = black] (4a) circle (0.050);
\draw [fill = black] (4b) circle (0.050);
\draw [fill = black] (4c) circle (0.050);

\coordinate (1a') at (4,0) {};
\coordinate (1b') at (5,0) {};
\coordinate (1c') at (6,0) {};
\coordinate (2a') at (4,1) {};
\coordinate (2b') at (5,1) {};
\coordinate (2c') at (6,1) {};
\coordinate (3a') at (4,2) {};
\coordinate (3b') at (5,2) {};
\coordinate (3c') at (6,2) {};
\coordinate (4a') at (4,3) {};
\coordinate (4b') at (5,3) {};
\coordinate (4c') at (6,3) {};

\draw (1a')--(1b')--(2b')--(2c')--(3c');
\draw (2a')--(3a')--(3b')--(4b')--(4c');

\draw [fill = black] (1a') circle (0.050);
\draw [fill = black] (1b') circle (0.050);
\draw [fill = black] (1c') circle (0.050);
\draw [fill = black] (2a') circle (0.050);
\draw [fill = black] (2b') circle (0.050);
\draw [fill = black] (2c') circle (0.050);
\draw [fill = black] (3a') circle (0.050);
\draw [fill = black] (3b') circle (0.050);
\draw [fill = black] (3c') circle (0.050);
\draw [fill = black] (4a') circle (0.050);
\draw [fill = black] (4b') circle (0.050);
\draw [fill = black] (4c') circle (0.050);

\coordinate (1a'') at (8,0) {};
\coordinate (1b'') at (9,0) {};
\coordinate (1c'') at (10,0) {};
\coordinate (2a'') at (8,1) {};
\coordinate (2b'') at (9,1) {};
\coordinate (2c'') at (10,1) {};
\coordinate (3a'') at (8,2) {};
\coordinate (3b'') at (9,2) {};
\coordinate (3c'') at (10,2) {};
\coordinate (4a'') at (8,3) {};
\coordinate (4b'') at (9,3) {};
\coordinate (4c'') at (10,3) {};

\draw (1a'')--(1b'')--(1c'')--(2c'');
\draw (2a'')--(2b'')--(3b'')--(3c'');
\draw (3a'')--(4a'')--(4b'')--(4c'');

\draw [fill = black] (1a'') circle (0.050);
\draw [fill = black] (1b'') circle (0.050);
\draw [fill = black] (1c'') circle (0.050);
\draw [fill = black] (2a'') circle (0.050);
\draw [fill = black] (2b'') circle (0.050);
\draw [fill = black] (2c'') circle (0.050);
\draw [fill = black] (3a'') circle (0.050);
\draw [fill = black] (3b'') circle (0.050);
\draw [fill = black] (3c'') circle (0.050);
\draw [fill = black] (4a'') circle (0.050);
\draw [fill = black] (4b'') circle (0.050);
\draw [fill = black] (4c'') circle (0.050);

\coordinate (1a''') at (12,0) {};
\coordinate (1b''') at (13,0) {};
\coordinate (1c''') at (14,0) {};
\coordinate (2a''') at (12,1) {};
\coordinate (2b''') at (13,1) {};
\coordinate (2c''') at (14,1) {};
\coordinate (3a''') at (12,2) {};
\coordinate (3b''') at (13,2) {};
\coordinate (3c''') at (14,2) {};
\coordinate (4a''') at (12,3) {};
\coordinate (4b''') at (13,3) {};
\coordinate (4c''') at (14,3) {};

\draw (1a''')--(1b''')--(1c''');
\draw (2a''')--(2b''')--(2c''');
\draw (3a''')--(3b''')--(3c''');
\draw (4a''')--(4b''')--(4c''');

\draw [fill = black] (1a''') circle (0.050);
\draw [fill = black] (1b''') circle (0.050);
\draw [fill = black] (1c''') circle (0.050);
\draw [fill = black] (2a''') circle (0.050);
\draw [fill = black] (2b''') circle (0.050);
\draw [fill = black] (2c''') circle (0.050);
\draw [fill = black] (3a''') circle (0.050);
\draw [fill = black] (3b''') circle (0.050);
\draw [fill = black] (3c''') circle (0.050);
\draw [fill = black] (4a''') circle (0.050);
\draw [fill = black] (4b''') circle (0.050);
\draw [fill = black] (4c''') circle (0.050);

\end{tikzpicture}
\caption{The lattice paths used to construct increasing subsequences of $ w(S,k) $ for $ b = 4, k = 3, t = 1, 2, 3, 4 $ }
\label{fig:exLPs}
\end{figure}
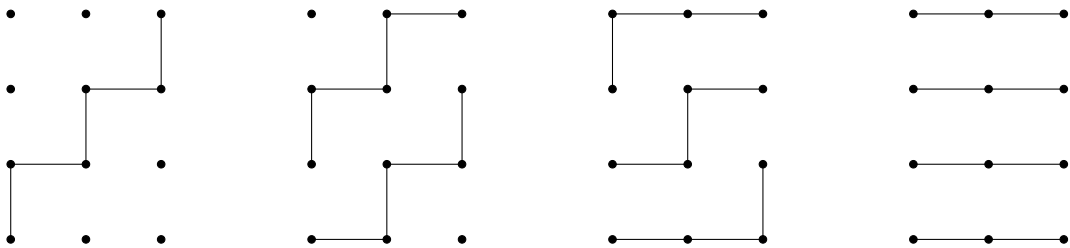

\Cref{fig:exLPs} corresponds to the following lists of disjoint increasing subsequences of $ w_1^{(3)} w_2^{(3)} w_3^{(3)} w_4^{(3)} $:
\begin{enumerate}[1.]
\item $ u_1 (u_2 + 12) (u_3 + 24) $ with size 15,
\item $ w_1(u_1 + 12)(u_2 + 24), u_2(u_3 + 12)(w_4 + 24) $ with total size 25,
\item $ w_1(w_1 + 12)(u_1 + 12), w_2(u_2 + 12)(w_3 + 24), u_3 (w_4 + 12)(w_4 + 24) $ with total size 33,
\item $ w_1^{(3)}, w_2^{(3)}, w_3^{(3)}, w_4^{(3)} $ with total size $ 36 $.
\end{enumerate}
Letting $ (\lam_1, \lam_2, \lam_3, \lam_4) = \rv(T_4^{(3)}(S)) $, \Cref{lem:shupperbound} says
\begin{align*}
	\lam_1 & \ge 3r + \sum_{i = 1}^3 c_i = 15, \\
	\lam_1 + \lam_2 & \ge 15 + 3r + \sum_{i = 2}^2 c_i = 25, \\ 
	\lam_1 + \lam_2 + \lam_3 & \ge 25 + 3r + \sum_{i = 3}^1 c_i = 33, \\
	 \lam_1 + \lam_2 + \lam_3 + \lam_4 & \ge 33 + 3r + \sum_{i = 4}^0 c_i = 36,
\end{align*}
which is consistent with \Cref{thm:Greene} and the total sizes of our increasing subsequences. Now, not only do these inequalities hold, but they are all equalities. We have $ \rv(T_4^{(3)}(S)) = (15, 10, 8, 3) $, as shown below:
\[
	\Yvcentermath1 T_4^{(3)}(S) = \Rect(S^{(3)}) = \begin{ytableau} 1 & 2 & 4 & 6 & 7 & 9 & 10 &  *(yellow) 14 &  *(yellow) 16 &  *(yellow) 19 &  *(yellow) 21 &  *(yellow) 22 & *(green) 31 & *(green) 33 & *(green) 34 \\
	3 & 5 & 8 & 11 & 12 &  *(yellow) 23 &  *(yellow) 24 & *(green) 26 & *(green) 28 & *(green) 36 \\
	*(yellow) 13 & *(yellow) 15 &  *(yellow) 17 &  *(yellow) 18 &  *(yellow) 20 & *(green) 30 & *(green) 32 & *(green) 35 \\
	*(green) 25 & *(green) 27 & *(green) 29
\end{ytableau} \, .
\]
\end{Example}

\begin{proof}[Proof of \Cref{thm:stabshape}] Suppose $ S $ is a standard skew tableau with row vector $ (r^b) $, and continue using \Cref{not:tabstab2}. Fix $ k \ge b - 1 $. Recall that
\[
	T_i^{(k)}(S) = P(w_1^{(k)} w_2^{(k)} \dots w_i^{(k)} ) \quad \tx{ for } \quad i = 1, \dots, b.
\]
We need to show 
\begin{align}
\label{eq:rvformula}
	\rv(T_b^{(k)}(S))_j = kr + \sum_{i = j}^{b - j} c_i \quad \tx{ for all } \quad j = 1, \dots, b \; \tx{ and } \; k \ge b - 1.
\end{align}
We proceed by induction on $ b $. If $ b = 1 $, $ P(w_1^{(k)}) $ is a tableau with 1 row of size $ kr $, so \Cref{thm:stabshape} holds for $ b = 1 $. Inductively assume that \Cref{thm:stabshape} holds for all standard skew tableaux with less than $ b $ rows and all $ k \ge b - 2 $. 

Let $ B_1, \dots, B_{kr} $ denote the bumping chains of the row insertions 
\[
	T_{b - 1}^{(k)}(S) = T_{b - 2}^{(k)}(S) \leftarrow w_{b - 1}^{(k)}.
\]
Each letter of $ w_{b - 1}^{(k)} $ has its own bumping chain. Since $ w_{b -1}^{(k)} $ is increasing, $ B_1, \dots, B_{kr} $ proceed strictly left to right by \Cref{lem:bumpingcomparisonorig}. By the induction hypothesis, we have
\begin{align*}
	\rv(T_{b - 2}^{(k)}(S))_j & = kr + \sum_{i = j}^{b - 2 - j} c_i, \tx{ for all } j = 1, \dots, b - 2, \\
	\rv(T_{b - 1}^{(k)}(S))_j & = kr + \sum_{i = j}^{b - 1- j} c_i, \tx{ for all } j = 1, \dots, b - 1.    
\end{align*}
Therefore, for $ j = 1, \dots, b - 2 $, 
\[
	\rv(T_{b - 1}^{(k)}(S))_j - \rv(T_{b - 2}^{(k)}(S))_j = c_{b - 1 - j},
\] 
which implies that $ c_{b - 1 - j} $ of the bumping chains among $ B_1, \dots, B_{kr} $ end in row $ j $. Furthermore, for $ t = 1, \dots, b - 2 $, exactly $ c_{b - 1- t} + \dots + c_{b - 2} $ of $ B_1, \dots, B_{kr} $ end in or above row $ t $.

Next, let $ B_1', \dots, B_{kr}' $ denote the bumping chains of the row insertions 
\[
	T_b^{(k)}(S) = T_{b - 1}^{(k)}(S) \leftarrow w_{b}^{(k)}.
\]
Since $ w_{b}^{(k)} $ is increasing, $ B_1', \dots, B_{kr}' $ again proceed strictly left to right by \Cref{lem:bumpingcomparisonorig}. We claim that $ B_j' $ is weakly to the left of $ B_{j + c_{ b - 1}} $ for all $ j = 1, \dots, kr - c_{b - 1} $, where
\[
	c_{b - 1} = \rv(P(w_{b - 1} w_{b}))_1 - r,
\]
from \eqref{eq:ci}. By \Cref{lem:uvcomparisonshift},
\begin{align}
\label{eq:lastcomparison}
	(w_{b}^{(k)})_j < (w_{b - 1}^{(k)})_{j + c_{b - 1}} \tx{ for all }  j = 1, \dots, kr - c_{b - 1}.
\end{align}
Since $ B_2, \dots, B_{kr} $ are disjoint from $ B_1 $, and $ (w_{b}^{(k)})_1 < (w_{b - 1}^{(k)})_{1 + c_{b - 1}} $, \Cref{lem:bumpingcomparison} implies that $ B_1' $ is weakly left of $ B_{1 + c_{b - 1}} $. Inductively assuming that $ B_1', \dots, B_j' $ are weakly left of $ B_{1 + c_{b - 1}}, \dots, B_{j + c_{b - 1}} $, respectively, $ B_{j + 2 + c_{b - 1}}, \dots, B_m, B_1', \dots B_j' $ are all disjoint from $ B_{ j + 1 + c_{ b - 1}} $ and $ (w_{b}^{(k)})_{j + 1} < (w_{b - 1}^{(k)})_{j + c_{b - 1}} $. Thus, by \Cref{lem:bumpingcomparison}, $ B_{j + 1}' $ is weakly left of $ B_{j + 1 + c_{b - 1}} $, which proves our claim. Hence, $ B_j' $ must end in a strictly lower row than $ B_{j + c_{ b - 1}} $ for all $ j = 1, \dots, kr - c_{b - 1} $.

Combining this with exactly $ c_{b - t} + \dots + c_{b - 2} $ of $ B_1, \dots, B_{kr} $ ending in or above row $ t - 1 $, at most $ c_{b - t} + \dots + c_{b - 2} + c_{b - 1} $ of $ B_1', \dots, B_{kr}' $ end in or above row $ t $ for $ t = 1, \dots, b - 1 $. Thus, using our original induction hypothesis,
\begin{equation}
\label{eq:shlowerbound}
\begin{aligned}
	\sum_{ j = 1}^t \rv(T_b^{(k)}(S))_j & \le \lp \sum_{ j = 1}^t \rv(T_{b - 1}^{(k)}(S))_j \rp + c_{b - t} + \dots + c_{b - 1} \\
			& = \sum_{ j = 1}^t \lp kr + \sum_{i = j}^{b - 1- j} c_i \rp + \sum_{ j = 1}^t c_{b - j} \\
			 & = \sum_{ j = 1}^t \lp kr + \sum_{i = j}^{b - j} c_i \rp.
\end{aligned}
\end{equation}
for $ t = 1, \dots, b - 1 $. Combining \eqref{eq:shlowerbound} with \Cref{lem:shupperbound},
\begin{align}
\label{eq:t<b}
	\sum_{ j = 1}^t \rv(T_b^{(k)}(S))_j = \sum_{ j = 1}^t \lp kr + \sum_{i = j}^{b - j} c_i \rp \quad \tx{ for all } \quad t = 1, \dots, b - 1. 
\end{align}

For the case $ t = b $,
\begin{equation}
\label{eq:t=b}
\begin{aligned}
	\sum_{ j = 1}^b \lp kr + \sum_{i = j}^{b - j} c_i \rp = bkr + \sum_{j = 1}^b \sum_{i = j}^{j - 1} c_i = bkr = \sum_{ j = 1}^b \rv(T_b^{k}(S))_j
\end{aligned}
\end{equation}
by \Cref{lem:symsum} with $ t = b $ and $ \ell(w_1^{(k)} w_2^{(k)} \dots w_b^{(k)}) = bkr $. Combining \eqref{eq:t<b} and \eqref{eq:t=b} yields 
\begin{align*}
	\sum_{ j = 1}^t \rv(T_b^{(k)}(S))_j  = \sum_{ j = 1}^t \lp kr + \sum_{i = j}^{b - j} c_i \rp \tx{ for all } t = 1, \dots, b.
\end{align*}
This shows that the partial sums of the two sequences $ \{ \rv(T_b^{(k)}(S))_j \}_{j = 1}^b $ and $ \lb kr + \sum_{i = j}^{b - j} c_i \rb_{j = 1}^b $ both agree, so the sequences agree as well, proving \eqref{eq:rvformula}.

\end{proof}




\begin{proof}[Proof of \Cref{thm:stab}] Suppose $ S $ is a standard skew tableau with row vector $ (r^b) $. In order to show $ S $ stabilizes at $ b $, it suffices to show that
\begin{align}
	\rv(\Rect(S^{(b)}))_j - \rv(\Rect(S^{(b - 1)}))_j = \rv(S)_j \quad\tx{ for all } \quad j = 1, \dots, b.
\end{align}
by \Cref{def:stab2}(c). In fact, using \Cref{thm:stabshape},
\begin{align*}
	\rv(\Rect(S^{(b)}))_j - \rv(\Rect(S^{(b - 1)}))_j = br + \sum_{i = j}^{b - j} c_i - (b - 1)r - \sum_{i = j}^{b - j} c_i = r = \rv(S)_j
\end{align*}
for all $ j = 1, \dots, b $.
\end{proof}

\begin{Corollary} \label{cor:stabshapeearly} Suppose $ S $ is a standard skew tableau with row vector $ (r^b) $. Let $ w_1, \dots, w_b $ denote the entries in each row read from left to right, starting from the bottom. For $ k \ge \stab(S) - 1 $, $ \Rect(S^{(k)}) $ has shape $ (\lam_1, \dots, \lam_b) $, where
\begin{align*}
	\lam_j = kr + \sum_{i = j}^{b - j} c_i \quad \tx{ for all } j = 1, \dots, b,
\end{align*}
where 
\[
	c_i = \rv( P(w_i w_{i + 1}))_1 - r \quad \tx{ for all } i = 1, \dots, b - 1.
\]
\end{Corollary}

\begin{proof} This is exactly \Cref{thm:stabshape} with the bound $ k \ge b - 1 $ changed to $ k \ge \stab(S) - 1 $. Let $ m = |S| $ and suppose $ \stab(S) - 1 \le k \le b - 1 $. Because $ k \ge \stab(S) - 1 $, 
\[
	\left. \Rect(S^{(b - 1)}) \right|_{[km + 1, (b - 1)m]} \sim S^{(b - 1 - k)}
\] 
by \Cref{lem:stabcontinues}. Since $ \Rect(S^{(k)}) \sube \Rect(S^{(b - 1)}) $ as well, we have
\begin{align*}
	\rv( \Rect(S^{(k)}) )_j & = \rv( \Rect(S^{(b - 1)}) )j - \rv( S^{(b - 1 - k)} )j \\
				& = (b - 1)r + \sum_{i = j}^{b - j} c_i - (b - 1 - k)r \\
				& = kr + \sum_{i = j}^{b - j} c_i
\end{align*}
for all $ j = 1, \dots, b $.
\end{proof}

\begin{Example} \label{ex:equalrowsstabbound} Continuing with the same $ S $ as in \Cref{ex:incsubseqs},
\[
	\Yvcentermath1 \Rect(S^{(4)}) = \begin{ytableau} 1 & 2 & 4 & 6 & 7 & 9 & 10 & *(yellow) 14 & *(yellow) 16 & *(yellow) 19 & *(yellow) 21 & *(yellow) 22 & *(green) 31 & *(green) 33 & *(green) 34 & *(orange) 43 & *(orange) 45 & *(orange) 46 \\
	3 & 5 & 8 & 11 & 12 & *(yellow) 23 & *(yellow) 24 & *(green) 26 & *(green) 28 & *(green) 36 & *(orange) 38 & *(orange) 40 & *(orange) 48 \\
	*(yellow) 13 & *(yellow) 15 & *(yellow) 17 & *(yellow) 18 & *(yellow) 20 & *(green) 30 & *(green) 32 & *(green) 35 & *(orange) 42 & *(orange) 44 & *(orange) 47 \\
	*(green) 25 & *(green) 27 & *(green) 29 & *(orange) 37 & *(orange) 39 & *(orange) 41
\end{ytableau} \, ,
\]
so $ \left. \Rect(S^{(4)}) \right|_{[25,36]} \sim S, \left. \Rect(S^{(4)}) \right|_{[13,24]} \not\sim S $, and $ \stab(S) = 3 $. If
\[
	T = \begin{ytableau} \none & \none & \none & \none & \none & 5 & 10 & 11 \\
	\none & \none & \none & \none & 8 & 9 & 12 \\
	\none & 2 & 4 & 6 \\
	1 & 3 & 7
\end{ytableau} \, ,
\] 
then
\[
	\Rect(T^{(4)}) = \begin{ytableau} 1 & 2 & 4 & 5 & 8 & 9 & 10 & 11 & *(yellow) 17 & *(yellow) 22 & *(yellow) 23 & *(green) 29 & *(green) 34 & *(green) 35 &  *(orange) 41 &  *(orange) 46 &  *(orange) 47 \\
	3 & 6 & 12 & *(yellow) 14 & *(yellow) 16 & *(yellow) 18 & *(yellow) 20 & *(yellow) 21 & *(yellow) 24 & *(green) 32 & *(green) 33 & *(green) 36 & *(orange) 44 & *(orange) 45 & *(orange) 48 \\
	7 & *(yellow) 13 & *(green) 25 & *(green) 26 & *(green) 28 & *(green) 30 & *(orange) 38 & *(orange) 40 & *(orange) 42 \\
	*(yellow) 15 & *(yellow) 19 & *(green) 27 & *(green) 31 & *(orange) 37 & *(orange) 39 & *(orange) 43
\end{ytableau} \, ,
\]
so $ \left. \Rect(T^{(4)}) \right|_{[37,48]} \sim T $, $ \left. \Rect(T^{(4)}) \right|_{[25,36]} \not\sim T $, and $ \stab(S) = 4 $.

\end{Example}

As seen in \Cref{ex:equalrowsstabbound}, some, but not all, skew tableaux with $ b $ rows of the same size stabilize before $ b $. To prove the bound in \Cref{thm:stab} is tight, it suffices to find a standard tableau $ I_{b,r} $ with row vector $ (r^b) $ and $ \stab(I_{b,r}) = b $ for each $ b, r \in \bZ_{\ge 1} $. We can form such an $ I_{b,r} $ by making its reading word increasing, as in \Cref{ex:stabatb}. 

\begin{Example} \label{ex:stabatb} Consider
\[
	\ytableausetup{boxsize = 1.8em} S = \begin{ytableau}
	\none & \none & \none & \none & \none & \none & \none & \none & \none & \none & & \dots & br \\
	\none & \none & \none & \none & \none & \none & \none & \none & \none & \reflectbox{$\ddots$} \\
	\none & \none & \none & \none & \none & \none & \scriptstyle{2r + 1} & \dots & 3r \\
	\none & \none & \none & \scriptstyle{r + 1} & \dots & 2r \\
	1 & \dots & r \\
	\end{ytableau}
	\ytableausetup{boxsize = 1.5em}
\]
Then, we have $ c_1 = c_2 = \dots = c_{b - 1} = r $, so in particular, setting $ k = b -1 $ and $ j = b $  in \Cref{thm:stabshape} gives
\[
	\rv( \Rect(I_{b,r}^{(b - 1)}) )_b = (b - 1)r + \sum_{i = b}^{0} c_i = (b -1)r - \sum_{i = 1}^{b -1} r = 0. 
\]
Hence, $ \Rect(I_{b,r}^{(b - 1)}) $ consists of at most $ b -1 $ rows. Therefore, $ \left. \Rect(I_{b,r}^{(b - 1)})  \right|_{[b(b - 2) + 1, b(b - 1)]} \not \sim I_{b,r} $, so $ I_{b,r} $ does not stabilize at $ b - 1 $. But, by \Cref{thm:stab}, $ I_{b,r} $ stabilizes at $ b $, so $ \stab(I_{b,r}) = b $. Note that $ I_{b,r} $ is not the only skew tableaux with $ b $ rows of size $ r $ and $ \stab(I_{b,r}) = b $, as illustrated by $ T $ in \Cref{ex:equalrowsstabbound}.

\end{Example}

\begin{Remark} We believe any standard tableau $ S $ with decreasing row vector stabilizes at $ b $, see \Cref{conj:stab}, though we are currently unable to prove it. We cannot directly generalize our proof of \Cref{thm:stab} to prove \Cref{conj:stab} because \Cref{thm:stabshape}, a key ingredient in our proof of \Cref{thm:stab}, does not generalize to this setting. Let $ w_1, \dots, w_b $ denote the reading words of the rows of $ S $, from bottom to top, $ r_i = \ell(w_i) $ for all $ i $, and
\[
	c_i \coloneqq \rv(P(w_i w_{i + 1}))_1 - r_{i + 1} \tx{ for all } i = 1, \dots, (b - 1). 
\] 
Unlike \Cref{thm:stabshape}, in this more general setting, $ \rv(\Rect(S^{(b)})) $ depends on more than just $ b, r_1, \dots, r_b, c_1, \dots, c_{b - 1} $, as demonstrated in \Cref{ex:samecjsdiffhape}. Hence, some other technique is required to prove \Cref{conj:stab}, if it is indeed true. 

\end{Remark}

\begin{Example} \label{ex:samecjsdiffhape} Consider 
\[
	S = \begin{ytableau}
\none & \none & 2 & 4  \\
 \none & 3\\
 1 \\
\end{ytableau} \,, \quad
	T = \begin{ytableau}
\none & \none & 1 & 3  \\
 \none & 4 \\
 2 \\
\end{ytableau} \, ,
\]
which both have $ b = 3, r_1 = 1, r_2 = 1, r_3 = 2, c_1 = 1, c_2 = 0 $. Then,
\begin{align*}
	\Rect(S^{(3)}) & = \begin{ytableau}
1 & 2 & 4 & *(yellow) 6 & *(yellow) 8 &  *(green)  10 &  *(green) 12  \\
 3 & *(yellow) 7 &  *(green)  11 \\
*(yellow) 5 &  *(green)  9 \\
\end{ytableau} \, , \quad
	\Rect(T^{(3)}) = \begin{ytableau}
1 & 3 & *(yellow) 5 & *(yellow) 7 &  *(green) 9 &  *(green) 11 \\
 2 & 4 & *(yellow) 8 &  *(green)  12 \\
 *(yellow) 6 &  *(green)  10 \\
\end{ytableau} \, .
\end{align*}
Thus, $ \rv(\Rect(S^{(3)}) ) = (7,3,2) \ne (6,4,2) = \rv(\Rect(T^{(3)})) $.  
\end{Example}


\section{Anti-Stabilization}
\label{sec:antistab}

In this section, we discuss anti-stabilization, which involves rectifying toward a southeast corner rather than the northwest corner. We show that the stabilized and anti-stabilized skew tableaux have reflected shapes and that the stabilization and anti-stabilization statistics coincide when both are defined.

\begin{Definition} \label{def:antirect} For any skew tableau $ S $ with $ b $ rows, fix $ a $ so that $ S \sub (a^b) $, and let $ \Rect^\ast(S) $, the \emph{anti-rectification} of $ S $, be the skew tableau obtained by continually performing outer slides on $ S $ within $ (a^b) $ until the southeast corner of $ (a^b) $ is filled. $ \Rect^\ast(S) $ is independent of $ a $ up to row shift equivalence. Also, for a standard skew tableau $ S $ of size $ m $, let $ S^\dagger $ denote the tableau obtained from $ S $ by rotating $ 180^\circ $ and then flipping the entries by $ x \mapsto m + 1 - x $, and let $ \lp \lam/\mu \rp^{\dagger} $ denote the skew shape obtained by rotating $ \lam/\mu $ by $ 180^{\circ} $.

	Suppose $ S $ is a standard skew tableau with $ m $ entries and increasing row vector. Define $ S^{(\ast k)} $ to be the standard skew tableau obtained by attaching $ k -1 $ shifted copies of $ S $ to the left of $ S + (k - 1)m $ so that $ \left. S^{(\ast k)} \right|_{[(j - 1)m + 1, jm]} \sim S $ for all $ j = 1, \dots, k $. Then, we say $ S $ \emph{anti-stabilizes} at $ k $ if $ \left. \Rect^\ast(S^{(\ast k)}) \right|_{[1,m]} \sim S $. Let $ \stab^\ast(S) $ denote the minimum value at which $ S $ anti-stabilizes.

\end{Definition}

\begin{Example} Consider 
\[
	\Yvcentermath1 S = \begin{ytableau} \none & \none & 4 \\
				      \none & 2 & 5 \\
				    1 & 3
		\end{ytableau} \, . \\ 
\]
Then,
\begin{align*}
	 \Yvcentermath1 S^{(\ast 3)} & = \begin{ytableau} \none & \none & \none & \none & 4 & *(yellow) 9 &  *(green)  14 \\
				      \none & 2 & 5 & *(yellow) 7 & *(yellow) 10 &  *(green)  12 &  *(green)  15 \\
				    1 & 3 & *(yellow) 6 & *(yellow) 8 &  *(green)  11 &  *(green)  13
	\end{ytableau} \, , \\ 
	\Yvcentermath1 \Rect^{\ast}(S^{(\ast 3)}) & = \begin{ytableau} \none & \none & \none & \none & \none & 4 & *(yellow) 9 \\
				      \none & 2 & 5 & *(yellow) 7 & *(yellow) 10 &  *(green)  12 &  *(green) 14 \\
				    1 & 3 & *(yellow) 6 & *(yellow) 8 &  *(green) 11 &  *(green)  13 &  *(green) 15
		\end{ytableau} \, , \\ 
	\Yvcentermath1 \Rect^{\ast}(S^{(\ast 3)})^{\dagger} & = \begin{ytableau} *(green) 1 & *(green) 3 & *(green) 5 & *(yellow) 8 & *(yellow) 10 & 13 & 15 \\
				      *(green) 2 & *(green) 4 & *(yellow) 6 & *(yellow) 9 & 11 & 14 \\
				    *(yellow) 7 & 12
		\end{ytableau} \, . 	
\end{align*}

\noindent Since $ \left. \Rect^{\ast}(S^{(\ast 3)}) \right|_{[11,15]} \not \sim S $, but $ \left. \Rect^{\ast}(S^{(\ast 3)}) \right|_{[6,10]} \sim S $, we have $ \stab^\ast(S) = 2 $.

\end{Example}

\begin{Remark} Similarly to \Cref{rem:rowsizes}, if $ S $ does not have increasing row vector, $ S^{(\ast k)} $ need not be a standard skew tableau. Hence, the notions of stabilization and anti-stabilization only make sense simultaneously when the row vector is constant.

\end{Remark}

Anti-stabilization has many of the same properties as stabilization. If we apply the same reasoning using anti-rectification instead of rectification and $ S^{(\ast k)} $ instead of $ S^{(k)} $, we get the following analogues of \Cref{lem:stabrse}, \Cref{thm:dualequivstab}, \Cref{lem:stabcontinues}, \Cref{thm:stabweaklydecrows}, and \Cref{thm:stab} for anti-stabilization. 


\begin{Corollary} \label{cor:antistabfacts}
Suppose $ S $ is a standard skew tableau with $ b $ rows and increasing row vector. Then,
\begin{enumerate}[(a)]
\item Anti-stabilization is well-defined up to row shift equivalence.
\item $ \stab^\ast $ is constant on dual equivalence classes of standard skew tableaux.
\item $ S $ anti-stabilizes eventually, and if $ S $ anti-stabilizes at $ k $, it anti-stabilizes at any $ k' \ge k $.
\item If $ b \ge 2 $, then $ S $ anti-stabilizes at $ 2b - 2 $.
\item If $ S $ has constant row vector, then $ S $ anti-stabilizes at $ b $. 
\end{enumerate}

\end{Corollary}

%
%

If $ S $ is a standard skew tableau with constant row vector, then $ S^{(k)} $ and $ S^{(\ast k)} $ are both standard tableaux with the same reading word. Thus, by \Cref{lem:rectP},
\begin{align}
\label{eq:leftrightsim}
	S^{(\ast k)} \sim S^{(k)}, \; \Rect(S^{(\ast k)}) = \Rect(S^{(k)}), \; \tx{ and } \Rect^\ast(S^{(\ast k)}) = \Rect^\ast(S^{(k)}).
\end{align}

\begin{Lemma} \label{lem:evacflip} For any standard skew tableau $ S $,
\[
	\Rect^\ast(S) = \Rect(S^\dagger)^\dagger = ( e( \Rect(S) ) )^\dagger,
\]
where $ e $ is Sch\"utzenberger's evacuation operator.
\end{Lemma}

\begin{proof} Let $ w = w_1 \dots w_m $ be the reading word of $ S $, so $ w $ is a permutation of size $ m $. Anti-rectifying $ S $ and rectifying $ S^\dagger $ differ by a rotation of $ 180^\circ $, implying
\begin{align}
\label{eq:rectvsantirect}
	\Rect^\ast(S) = \Rect(S^\dagger)^\dagger.
\end{align}
As $ \dagger $ reverses the order of the cells in the reading word and reverses the entry values, the reading word of $ S^\dagger $ is $ w_0 w w_0 $, where $ w_0 = [m, m - 1, \dots, 2, 1] $ is the reversal permutation of size $ m $. Thus, by \Cref{lem:rectP},
\begin{align}
\label{eq:conjrw}
	\Rect(S^\dagger)^\dagger = P(w_0 w w_0)^{\dagger}. 
\end{align}
Moreover, $ P(w_0 w w_0) = e(P(w)) $ by \cite[Theorem 3.9.4]{MR1824028}, so
\begin{align}
\label{eq:rectevac}
	 P(w_0 w w_0)^{\dagger} = e(P(w))^{\dagger} = ( e( \Rect(S) ) )^\dagger
\end{align}
by \Cref{lem:rectP}. Combining \eqref{eq:rectvsantirect}, \eqref{eq:conjrw}, and \eqref{eq:rectevac} gives \Cref{lem:evacflip}.


\end{proof}

\begin{Lemma} \label{lem:antistab} Suppose that $ S $ is a standard skew tableau with row vector $ (r^b) $. Then, for all $ k \ge 1 $, and $ j = 1, \dots, b $,
\begin{align}
\label{eq:antistabshape}
	\rv( \Rect^\ast ( S^{(k)} ) )_j = \rv( \Rect ( S^{(k)} ) )_{b + 1 - j}.
\end{align}
In addition, 
\begin{align}
\label{eq:antistab}
	\stab^\ast(S) = \stab(S). 
\end{align}
\end{Lemma}

\begin{proof} For any $ k \ge 1 $, we have $ S^{(k)} \sim S^{(\ast k)} $ by \eqref{eq:leftrightsim}. Thus, by \Cref{lem:evacflip}, 
\begin{align*}
	\rv( \Rect^\ast ( S^{(k)} ) )_j =  \rv ( e(\Rect(S^{(k) } ) )^\dagger )_j = \rv ( \Rect(S^{(k) })^\dagger )_j = \rv ( \Rect(S^{(k) }) )_{b + 1 - j}
\end{align*}
for all $ k \ge 1 $ and $ j = 1, \dots, b $, proving \eqref{eq:antistabshape}. By \Cref{def:stab2}(c), $ S $ stabilizes at $ k $ if and only if
\begin{align}
\label{eq:stabchar}
	\rv(\Rect( S^{(k)} ))_j - \rv( \Rect(S^{(k - 1)}) )_j = r \quad \tx{ for all } j = 1, \dots, b.
\end{align}
Similarly, $ S $ anti-stabilizes at $ k $ if and only if
\begin{align}
\label{eq:antistabchar}
	\rv ( \Rect^\ast( S^{(k)} ))_{j} - \rv( \Rect^\ast(S^{(k - 1)} ) )_{j} = r  \quad \tx{ for all } j = 1, \dots, b.	
\end{align}
By \eqref{eq:antistabshape}, \eqref{eq:stabchar} holds if and only if \eqref{eq:antistabchar} holds. Thus, $ S $ stabilizes at $ k $ if and only if $ S $ anti-stabilizes at $ k $, and \eqref{eq:antistab} follows.

\end{proof}

\section{Sufficiently Large Tableaux Fixed by Powers of Promotion}
\label{sec:prfixedtabs}

In this section, we first give an alternative method for doing multiple promotions at once. Recall $ \pr: \SYT(a^b) \to \SYT(a^b) $ denotes Sch\"utzenberger's promotion operator, see \Cref{def:pr}. We then construct the sufficiently large rectangular standard tableaux fixed by promotion powers. In particular, we construct the tableaux in $ \SYT((ar)^b)^{\pr^{br}} $ for $ a \ge 2b - 1 $ and prove \Cref{thm:prfixedtabs}. \Cref{thm:stabshape} plays a key role in showing these tableaux are rectangular, tableau stabilization is central in showing these tableaux are fixed by $ \pr^{br} $, and the bound in \Cref{thm:stab} lets us control the size of these rectangular tableaux. We also prove \Cref{cor:prfixedpr}, which describes the action of promotion on $ \SYT((ar)^b)^{\pr^{br}} $.

\begin{Lemma} \label{lem:multprom} Suppose $ a, b \in \bZ_{\ge 1} $ and $ n \coloneqq ab $. For any tableau $ T \in \SYT(a^b) $ and $ k = 1, \dots, n $, we have
\begin{align}
\label{eq:multprom1}
	\left. \pr^k(T) \right|_{[1,n - k]} & = \Rect( \left. T \right|_{[k + 1,n]} ) - k, \\
\label{eq:multprom2}
	\left. \pr^k(T) \right|_{[n - k + 1, n]} & = \Rect^\ast( \left. T \right|_{[1, k]}) + (n - k).
\end{align}

\end{Lemma}

\begin{proof} Equation \eqref{eq:multprom1} follows from \Cref{def:pr} and rectification being well-defined. By \Cref{thm:prnid},  $ \pr^n = \id $ on $ \SYT(a^b) $, so $ \pr^k = (\pr^{-1})^{n- k} $. Then, \eqref{eq:multprom2} follows from the fact that $ T $ is rectangular, meaning demotion and promotion are inverses, \Cref{def:pr}, and anti-rectification being well-defined.

\end{proof}


\begin{Definition} \label{def:Raconstruct} Suppose $ S $ is a standard skew tableau with row vector $ (r^b) $, and let $ k \coloneqq \stab(S) $. For any $ a \ge 2k - 1 $, let $ R_a(S) $ denote the filling of rectangular shape formed by row-concatenating $ \Rect(S^{(k - 1)}), S^{(a - 2k + 2)} + (k - 1)br, $ and $ \Rect^\ast(S^{(k - 1)}) + (a - k + 1)br $ together from left to right.

\end{Definition}

\begin{Example} Let 
\[
	S = \begin{ytableau} \none & \none & \none & \none & 2 & 6 \\
				      \none & \none & 4 & 5 \\
				    1 & 3
		\end{ytableau} \,,
\]
which has $ k = \stab(S) = 3 $, and let $ a = 6 $. Observe that
\begin{align*}
	\Rect(S^{(2)}) =  \begin{ytableau} 1 & 2 & 4 & 5 & 6 &  *(yellow) 8 & *(yellow) 12 \\
				     3 & *(yellow) 9 & *(yellow) 10 & *(yellow) 11 \\
				    *(yellow) 7
		\end{ytableau} \, , \\
	S^{(2)} + 12 = \begin{ytableau} \none & \none & \none & \none & *(green) 14 & *(green) 18 & *(orange) 20 & *(orange) 24 \\
				     \none & \none & *(green) 16 & *(green) 17 & *(orange) 22 & *(orange) 23 \\
				    *(green) 13 & *(green) 15 & *(orange) 19 & *(orange) 21
		\end{ytableau} \, , \\
	\Rect^\ast(S^{(2)}) + 24 = \begin{ytableau} \none & \none & \none & \none & \none & \none & *(cyan) 26 \\
				     \none & \none & \none & *(cyan) 28 & *(cyan) 29 & *(cyan) 30 & *(red) 32 \\
				    *(cyan) 25 & *(cyan) 27 & *(red) 31 & *(red) 33 & *(red) 34 & *(red) 35 & *(red) 36 \\
		\end{ytableau} \, ,
\end{align*}
so
\[
	R_6(S) = \begin{ytableau} 1 & 2 & 4 & 5 & 6 & *(yellow) 8 & *(yellow) 12 & *(green) 14 & *(green) 18 & *(orange) 20 & *(orange) 24 & *(cyan) 26 \\
				     3 & *(yellow) 9 & *(yellow) 10 & *(yellow) 11 & *(green) 16 & *(green) 17 & *(orange) 22 & *(orange) 23 & *(cyan) 28 & *(cyan) 29 & *(cyan) 30 & *(red) 32 \\
				    *(yellow) 7 & *(green) 13 & *(green) 15 & *(orange) 19 & *(orange) 21 & *(cyan) 25 & *(cyan) 27 & *(red) 31 & *(red) 33 & *(red) 34 & *(red) 35 & *(red) 36
\end{ytableau} \, .
\]
\end{Example}

\bs

\begin{Thm} \label{thm:fixedbyprom} Suppose that $ S $ is a standard skew tableau with row vector $ (r^b) $. Then, for any integer $ a \ge 2\stab(S) - 1$,
\[
	R_a(S) \in \SYT((ar)^b)^{\pr^{br}}. 
\]
\end{Thm}

\begin{proof} Let $ k \coloneqq \stab(S) $ and $ R \coloneqq R_a(S) $. First, we check that $ R $ has rectangular shape $ ((ar)^b) $. By \Cref{cor:stabshapeearly} and \eqref{eq:antistabshape},
\begin{align}
\label{eq:innershape*}
	\rv(\Rect(S^{(k - 1)}))_j & = (k - 1)r + \sum_{i = j}^{b - j} c_i, \\ 
	\rv(\Rect^\ast(S^{(k - 1)}))_j & = \rv(\Rect(S^{(k - 1)}))_{b + 1 - j} =  (k - 1)r + \sum_{i = b + 1 - j}^{j - 1} c_i
\end{align}
for all $ j = 1, \dots, b $. Hence, for $ j = 1, \dots, b $,
\begin{align*}
	\rv(R)_j & = \rv(\Rect(S^{k - 1}))_j  + \rv(S^{(a - 2k + 2)})_j + \rv(\Rect^\ast(S_{k - 1}))_j \\  
		& = (k - 1)r + \sum_{i = j}^{b - j} c_i + (a - 2k + 2)r + (k - 1)r + \sum_{i = b + 1 - j}^{j - 1} c_i \\
		& = ar + \sum_{i = j}^{j - 1} c_i \\
		& = ar,
\end{align*}
which shows $ \rv(R) = ((ar)^b) $. Since $ R $ is left-justified, it is a filling of shape $ ((ar)^b) $.

Secondly, we check that $ R $ is a standard tableau. We break $ R $ into three pieces: 
\begin{itemize}
\item $ R_1 \coloneqq \left. R \right|_{[1,kbr]} = \Rect(S^{(k)}) $, \tx{ because $ S $ stabilizes at $ k $ }
\item $ R_2 \coloneqq \left. R \right|_{[kbr + 1, (a - k + 1)br]} \sim S^{(a - 2k + 1)} $,
\item $ R_3 \coloneqq \left. R \right|_{[(a - k + 1)br + 1, abr]} = \Rect^\ast(S^{(k - 1)}) + (a - k + 1)br $. 
\end{itemize}
The fillings $ R_1 $ and $ R_3 $ are skew tableaux because they are the rectification and anti-rectification of skew tableaux, respectively. The filling $ R_2 $ is a skew tableau since it is formed by row-concatenating shifted copies of $ \left. R_1 \right|_{[(k - 1)br + 1, kbr]} $. Since $ R_1, R_2, R_3 $ use the entries $ [1, kbr] $, $ [kbr + 1, (a - k + 1)br] $, and $ [(a - k + 1)br + 1, abr] $ respectively, $ R $ uses each of $ 1, 2, \dots, abr $ exactly once. Because $ R_2 $ adjoins to the right of $ R_1 $ and $ R_3 $ adjoins to the right of $ R_2 $, the rows and columns of $ R $ are increasing. Hence, $ R $ is a standard tableau.

\ms

Thirdly, we check that $ R $ is fixed by $ br $ promotions. Since $ S $ stabilizes at $ k $,
\begin{align}
\label{eq:stabstart2}
	\left. R \right|_{[1, jbr]} = \Rect(S^{(j)}) \quad \tx{ for all } j \in [0, a - k + 1]. 
\end{align}
By \Cref{lem:antistab}, $ S $ also anti-stabilizes at $ k $, so similarly,
\begin{align}
\label{eq:stabstart1}
	\left. R \right|_{[jbr + 1,abr]} = \Rect^\ast(S^{(a - j)}) + jbr \quad \tx{ for all } j \in [k - 1, a].
\end{align}
Note that both sides of \eqref{eq:stabstart1} are empty when $ j = a $. Then, for all $ j \in [k - 1, a] $,
\begin{equation}
\label{eq:prfixedrest1}	
\begin{aligned}
	\left. \pr^{jbr}(R) \right|_{[1,(a - j)br]} & = \Rect( \left. R \right|_{[jbr + 1, abr]} ) - jbr \qquad \tx{ by } \eqref{eq:multprom1} & \\
			& = \Rect (\Rect^\ast(S^{ (a - j) } ) + jbr) - jbr \qquad \tx{ by } \eqref{eq:stabstart1} & \\
			& = \Rect(S^{ (a - j) }) \qquad \tx{ as rectification is well-defined } & \\
			& = \left. R \right|_{[1,(a- j)br]} \qquad \tx{ by } \eqref{eq:stabstart2}. &
\end{aligned}
\end{equation}
Similarly, for all $ j \in [1,a - k + 1] $,
\begin{equation}
\label{eq:prfixedrest2}
\begin{aligned}
	\left. \pr^{jbr}(R) \right|_{[(a - j)br + 1, abr]} & = \Rect^\ast( \left. R \right|_{[1,jbr]} ) + (a - j)br \qquad \tx{ by } \eqref{eq:multprom2} & \\
					& = \Rect^\ast( \Rect(S^{(j)}) ) + (a - j)br \qquad \tx{ by } \eqref{eq:stabstart2} & \\ 
					& = \Rect^\ast(S^{(j)}) + (a - j)br \quad \tx{ as anti-rectification is well-defined } & \\
					& = \left. R \right|_{[(a - j)br + 1, abr]} \qquad \tx{ by } \eqref{eq:stabstart1}. &
\end{aligned}
\end{equation}
Putting \eqref{eq:prfixedrest1} and \eqref{eq:prfixedrest2} together,
\[
	\pr^{jbr}(R) = R \; \tx{ for all } \; j \in [k - 1, a - k + 1].
\] 
In particular, because $ a \ge 2k -1 $, we have $ k - 1, k \in [k - 1, a - k + 1] $, so
\[
	\pr^{br}(R) = \pr^{kbr - (k - 1)br}(R) = \pr^{kbr}(\pr^{-(k - 1)br}(R)) = \pr^{kbr}(R) = R.
\]
\end{proof}

%

\begin{Notation} Fix $ b, r \in \bZ_{\ge 1} $, and let $ \be $ be the block anti-diagonal skew shape
\[
	\be \coloneqq \underbrace{(r) \cup (r) \cup \dots \cup (r)}_{b \tx{ times}}. 
\]
\end{Notation}

Although \Cref{thm:fixedbyprom} holds for any standard skew tableau $ S $ with row vector $ (r^b) $, $ S $ is row shift equivalent to a standard skew tableau $ S' $ of shape $ \be $ by performing horizontal slides. Thus, by \Cref{lem:stabrse}, $ \Rect(S^{(k)}) = \Rect((S')^{(k)}) $ and $ \stab(S) = \stab(S') $, which means $ R_a(S) = R_a(S') $ for all $ a \ge 2 \stab(S) - 1 $. Hence, any rectangular tableau of the form $ R_a(S) $ for some standard skew tableau $ S $ with row vector $ (r^b) $ is also of the form $ R_a(S') $ for some $ S' \in \SYT(\be) $.


\begin{Corollary} \label{cor:prfixedsub} For all $ a \in \bZ_{\ge 1} $,

\begin{align}
	\lb R_a(S) : S \in \SYT(\be), \, \stab(S) \le \f{a + 1}{2}  \rb \subseteq \SYT((ar)^b)^{\pr^{br}}.
\end{align}

\end{Corollary}

\begin{proof} For fixed $ S \in \SYT(\be) $, $ R_a(S) $ is defined for $ a \ge 2 \stab(S) - 1 $. Thus, for fixed $ a \in \bZ_{\ge 1} $ and all $ S \in \SYT(\be) $ with $ \stab(S) \le \f{a + 1}{2} $, $ R_a(S) $ is defined and lies in $  \SYT((ar)^b)^{\pr^{br}} $ by \Cref{thm:fixedbyprom}. 

It remains to show the elements of $ \lb R_a(S) : S \in \SYT(\be) \rb $ are distinct. Suppose that $ R_a(S) = R_a(S') $ for $ S, S' \in \SYT(\be) $. Letting $ k = \stab(S) $, we have $ \rv( \left. R_a(S) \right|_{[(k - 1)br + 1, k br]} ) = (r^b) $ and hence $ \left. \rv( R_a(S') \right|_{[(k - 1)br + 1, k br]} ) = (r^b) $ as well. This and \Cref{def:Raconstruct} mean
\[
	S' \sim  \left. R_a(S') \right|_{[(k - 1)br + 1,  k br]} = \left. R_a(S) \right|_{[(k - 1)br + 1, k br]} \sim S,
\] 
forcing $ S = S' $ because $ S, S' \in \SYT(\be) $.

\end{proof}

\begin{Example} \label{ex:notconstructed} Consider 
\[
	\Yvcentermath1 S = \begin{ytableau} \none & \none & \none & 2 \\ \none & \none & 1 \\ \none & 3 \\  4 \end{ytableau} \, , \qquad T = \begin{ytableau} \none & \none & \none & 2 \\ \none & \none & 1 \\ \none & 4 \\  3 \end{ytableau} \, ,
\qquad U = \begin{ytableau} \none & \none & \none & 4 \\ \none & \none & 3 \\ \none & 2 \\  1 \end{ytableau} \, , 
\]
which have
\[
	\stab(S) = 2, \qquad \stab(T) = 3, \qquad \stab(U) = 4.
\]
Thus, the smallest $ R_a(S), R_{a'}(T), R_{a''}(U) $ we can construct with \Cref{def:Raconstruct} are
\begin{align*}
	\Yvcentermath1 R_3(S) & = \begin{ytableau} 1 & 2 & *(yellow) 6 \\ 3 & *(yellow) 5 & *(green) 10 \\ 4 & *(yellow) 7 & *(green) 11 \\ *(yellow) 8 & *(green) 9 & *(green) 12 \end{ytableau} \in \SYT(3^4)^{\pr^4}, \\
	R_5(T) & = \begin{ytableau} 1 & 2 & *(yellow) 5 & *(yellow) 6 & *(green) 10 \\ 3 & 4 & *(green) 9 & *(orange) 13 &  *(orange) 14 \\ *(yellow) 7 & *(yellow) 8 & *(green) 12 &  *(cyan) 17 & *(cyan) 18 \\ *(green) 11 &  *(orange) 15 &  *(orange) 16 & *(cyan) 19 &*(cyan)  20 \end{ytableau} \in \SYT(5^4)^{\pr^4}, \\
	R_7(U) & = \begin{ytableau} 1 & 2 & 3 & 4 & *(yellow) 8 & *(green) 12 &  *(orange) 16 \\ *(yellow) 5 & *(yellow) 6 & *(yellow) 7 & *(green) 11 &  *(orange) 15 & *(cyan) 19 & *(cyan) 20 \\ *(green) 9 & *(green) 10 &  *(orange) 14 & *(cyan) 18 & *(red) 22 & *(red) 23 & *(red) 24 \\  *(orange) 13 & *(cyan) 17 & *(red) 21 & *(pink) 25 &  *(pink) 26 &  *(pink) 27 &  *(pink) 28 \end{ytableau} \in \SYT(7^4)^{\pr^4}.
\end{align*}
On the other hand, neither
\[
	\Yvcentermath1 V = \begin{ytableau} 1 & 3 & 4 \\ 2 & *(yellow) 5 & *(yellow) 8 \\ *(yellow) 6 & *(yellow) 7 & *(green) 9 \\ *(green) 10 & *(green) 11 & *(green) 12 \end{ytableau} \in \SYT(3^4)^{\pr^4}, \quad \tx{ nor } \quad W = \begin{ytableau} 1 & 2 & 3 & 4 & *(yellow) 8 \\ *(yellow) 5 & *(yellow) 6 & *(yellow) 7 & *(green) 11 & *(green) 12 \\ *(green) 9 & *(green) 10 & *(orange) 14 & *(orange) 15 & *(orange) 16 \\ *(orange) 13 & *(cyan) 17 & *(cyan) 18 & *(cyan) 19 & *(cyan) 20 \end{ytableau} \in \SYT(5^4)^{\pr^4}
\]
is of the form $ R _a(S) $ for $ S \in \SYT((1) \cup (1) \cup (1) \cup (1)) $ because $ \left. V \right|_{[4(j - 1) + 1,4j]} $ and $ \left. W \right|_{[4(j - 1) + 1,4j]} $ never use all 4 rows for any value of $ j $. 
\end{Example}

The containment in \Cref{cor:prfixedsub} can be strict in general, as illustrated by $ V, W $ in \Cref{ex:notconstructed}. However, equality in \Cref{cor:prfixedsub} holds for $ a \ge 2b -1 $, when $ \stab(S) \le \f{a + 1}{2} $ holds for all $ S \in \SYT(\be) $ by \Cref{thm:stab}. Thus, \Cref{cor:prfixedsub} specializes to
\begin{align}
\label{eq:prfixedsub2}
	\lb R_a(S) : S \in \SYT(\be) \rb \subseteq \SYT((ar)^b)^{\pr^{br}} \qquad \tx{ when } a \ge 2b -1.
\end{align}
The fact that equality holds in \eqref{eq:prfixedsub2} is the content of \Cref{thm:prfixedtabs}.

\begin{Remark} In \Cref{thm:prfixedtabs}, $ R_a(S) $ is defined by row concatenating shifted copies of $ \Rect(S^{(b - 1)}) $, $ S^{(a - 2b + 2)} $, and $ \Rect^\ast(S^{(b - 1)}) $ instead of shifted copies of $ \Rect(S^{(k - 1)}), $ $ S^{(a - 2k + 2)}, $ and $ \Rect^\ast(S^{(k - 1)}) $ where $ k = \stab(S) $ as in \Cref{def:Raconstruct}. These 2 constructions agree by the definition of tableau stabilization, see \Cref{def:stab}, and because $ \stab(S) \le b $ for $ S \in \SYT(\be) $, see \Cref{thm:stab}.

\end{Remark}

\begin{proof}[Proof of \Cref{thm:prfixedtabs}] Fix a positive integer $ a \ge 2b - 1 $. Equality in \eqref{eq:prfixedsub2} will follow from the fact that both sides have the same size. By \Cref{lem:rectquotients}, we have the quotient
\begin{align}
\label{eq:aquotientsufflarge}
	Q_e((er)^b) = \underbrace{(r) \cup (r) \cup \dots \cup (r)}_{b \tx{ times}} = \be \tx{ whenever } e \ge b.
\end{align} 
Thus, $ Q_a((ar)^b) = \be $ because $ a \ge 2b - 1 \ge b $. By \Cref{cor:prfixedribbontab} and \eqref{eq:aquotientsufflarge},
\begin{align*}
	\# \SYT((ar)^b)^{\pr^{br}} & = \# \SYT(Q_a((ar)^b) ) = \# \SYT \lp \be \rp = \# \lb R_a(S) : S \in \SYT(\be) \rb,
\end{align*}
which shows both sides of \eqref{eq:prfixedsub2} have the same size. Therefore,
\[
	\lb R_a(S) : S \in \SYT(\be) \rb = \SYT((ar)^b)^{\pr^{br}}.
\]

Finally, we can choose $ S \in \SYT(\be) $ by choosing the sets of $ r $ elements that go in each of $ b $ rows, which determines $ S $ uniquely since each row of $ S $ is increasing. Thus,
\[
	\# \SYT((ar)^b)^{\pr^{br}} = \ch{br}{r, \dots, r}.
\]
\end{proof}

Not only can we describe the elements of $ \SYT((ar)^b)^{\pr^{br}} $ for $ a \ge 2b -1 $, but we can also describe the action of promotion on $ \SYT((ar)^b)^{\pr^{br}} $, which is closed under promotion. In fact, using the definition of promotion for skew shapes, see \Cref{def:pr}, promotion commutes with the $ R_a $ operator:
\[
	\pr(R_a(S)) = R_a(\pr(S)) \tx{ for all } S \in \SYT(\be), 
\]
as in \Cref{cor:prfixedpr}.

\begin{proof}[Proof of \Cref{cor:prfixedpr}] For any $ T \in \SYT((ar)^b)^{\pr^{br}} $, we have
\[
	\pr^{br}(\pr(T)) = \pr(\pr^{br}(T)) = \pr(T), \quad \tx{ implying } \quad \pr(T) \in \SYT((ar)^b)^{\pr^{br}}. 
\] 
Thus, $ \SYT((ar)^b)^{\pr^{br}} $ is closed under promotion. So, fixing $ S \in \SYT(\be) $,
\[
	\pr(R_a(S)) = R_a(S') \qquad \tx{ for some } S' \in \SYT(\be)
\]
by \Cref{thm:prfixedtabs}. It suffices to show $ S' = \pr(S) $. As rectification is well-defined, we get the same result whether we perform an inner slide at 1's cell before or after rectification, so 
\begin{align}
\label{eq:1slide}
	\left. \pr(S^{(b)}) \right|_{[1,b^2 r - 1]} = \left. \Rect(S^{(b)} \right|_{[2,b^2 r]}) - 1.
\end{align}
Since $ S \in \SYT(\be) $, the inner slide started by removing 1's cell in $ S^{(b)} $ consists of only horizontal slides, so
\begin{align}
\label{eq:horizslide}
	\pr(S^{(b)}) = \pr(S)^{(b)}.
\end{align}
Then,
\begin{equation}
\label{eq:promprfixedrest}
\begin{aligned}
	\left. R_a(S') \right|_{[1,b^2 r - 1]} & = \left. \pr(R_a(S)) \right|_{[1,b^2 r - 1]} \\
				& = \Rect ( \left. R_a(S) \right|_{[2,b^2 r]} )  - 1 \qquad \tx{ by  \Cref{def:pr} } & \\
				& = \Rect ( \left. S^{(b)} \right|_{[2,b^2 r]} )  - 1 \qquad \tx{ by  \Cref{def:Raconstruct} } & \\
				& = \Rect ( \left. \pr(S^{(b)})  \right|_{[1,b^2 r - 1]} ) \qquad \tx{ by } \eqref{eq:1slide} & \\
				& = \Rect ( \left. \pr(S)^{(b)} \right|_{[1,b^2 r - 1]} ) \qquad \tx{  by } \eqref{eq:horizslide} & \\
				& = \left. R_a(\pr(S)) \right|_{[1,b^2 r - 1]} \qquad  \tx{ by \Cref{def:Raconstruct}.} &
\end{aligned}
\end{equation}
By \Cref{thm:stab}, $ \left. R_a(S') \right|_{[(b - 1) br + 1, b^2 r]} $ and $ \left. R_a(\pr(S)) \right|_{[(b - 1) br + 1, b^2 r]} $ each have row vector $ (r^b) $, which together with \eqref{eq:promprfixedrest}, forces $ b^2 r $ to be in the same cell in $ R_a(S') $ and $ R_a(\pr(S)) $. Combining this with \eqref{eq:promprfixedrest} again,
\[
	\left. R_a(S') \right|_{[1,b^2 r]} = \left. R_a(\pr(S)) \right|_{[1,b^2 r]}.  
\] 
In particular,
\[
	S' \sim \left. R_a(S') \right|_{[(b - 1)br + 1, b^2r]} = \left. R_a(\pr(S)) \right|_{[(b - 1)br + 1, b^2 r]} \sim \pr(S),
\]
forcing $ S' = \pr(S) $ since $ S', \pr(S) \in \SYT(\be) $.  
\end{proof}

\begin{Corollary} \label{cor:LRCcount} Let $ \lam = ((b + 1)r, br, \dots, 2r, r) $ and $ \mu = (br, (b - 1)r, \dots, r, 0) $. Then, for all integers $ a \ge 2b - 1 $, $ \nu \vdash br $, and $ T_0 \in \SYT(\nu) $,
\[
	\# \{ T \in \SYT((ar)^b)^{\pr^{br}} : \left. T \right|_{[1,br]} = T_0 \} = c^{\lam}_{\mu, \nu},
\]
the Littlewood-Richardson coefficient from \Cref{def:LRCs}.
\end{Corollary}

\begin{proof} Note that $ \be = \lam / \mu $. By \Cref{thm:prfixedtabs} and $ \left. R_a(S) \right|_{[1,br]} = \Rect(S) $, we have
\begin{align*}
	\# \{ T \in \SYT((ar)^b)^{\pr^{br}} : \left. T \right|_{[1,br]} = T_0 \} & = \# \left\{ R_a(S) : S \in \SYT(\be), \left. R_a(S) \right|_{[1,br]} = T_0 \right\} \\
			& = \# \left\{ R_a(S) : S \in \SYT(\be), \Rect(S) = T_0 \right\} \\
			& = \# \left\{ S \in \SYT(\lam / \mu) : \Rect(S) = T_0 \right\} \\
			& = c^{\lam}_{\mu, \nu}
\end{align*}
by \Cref{def:LRCs}.
\end{proof}

\section{Other Tableaux Fixed by Promotion Powers}
\label{sec:otherprfixed}

In this section, we present the construction of $ \SYT((2r)^b)^{\pr^{br}} $ by White and Rhee \cite{Whiteprom} \cite{RheeThesis}, see \Cref{thm:fixedbyab/2prom}, using a similar construction to \Cref{thm:prfixedtabs}. Then, we describe the action of promotion on $ \SYT((2r)^b)^{\pr^{br}} $ in \Cref{cor:prfixedpr2}. Next, we characterize the block diagonal skew tableaux fixed by a power of promotion in terms of tableaux of straight shape fixed by certain powers of promotion, see \Cref{thm:prfixedpieces}, inspired by White \cite{Whiteprom}. This has consequences for describing the tableaux in $ \SYT((2r)^b)^{\pr^{br/k}} $ in terms of rectangular tableaux fixed by smaller promotion powers, see \Cref{cor:R2prfixedwithin}.

We next describe the tableaux in $ \SYT((ar)^b)^{\pr^{br}} $ when $ a = 2 $. Using \Cref{lem:rectquotients} and that $ \left \lceil \f{b}{2} \right \rceil = \left \lfloor \f{b}{2} \right \rfloor $ when $ b $ is even, we have the quotient
\begin{align}
\label{eq:rect2quotient}
	Q_2((2r)^b) = ( r^{\left \lceil \f{b}{2} \right \rceil} ) \cup ( r^{\left \lfloor \f{b}{2} \right \rfloor} ).
\end{align}
Therefore, by \Cref{cor:prfixedribbontab},
\begin{align}
\label{eq:2quotientofrectangle}
	\# \SYT((2r)^b)^{\pr^{br}} = \# \SYT \lp ( r^{\left \lceil \f{b}{2} \right \rceil} ) \cup ( r^{\left \lfloor \f{b}{2} \right \rfloor} ) \rp.
\end{align} 

\begin{Definition} \label{def:R2construct} Fix $ b,r \in \bZ_{\ge 1} $ and let $ \g \coloneqq ( r^{\left \lceil \f{b}{2} \right \rceil} ) \cup ( r^{\left \lfloor \f{b}{2} \right \rfloor} ) $. For $ S \in \SYT(\g) $, let $ R_2(S) $ be the filling formed by row-concatenating $ \Rect(S) $ and $ \Rect^\ast(S) + br $ together from left to right.

\end{Definition}

\begin{Example} \label{ex:a=2} Suppose $ b = 4, r = 2 $, and choose
\[
	S = \begin{ytableau} \none & \none & 2 & 5 \\ \none & \none & 6 & 8 \\ 1 & 3 \\ 4 & 7 \end{ytableau} \in \SYT((2,2) \cup (2,2)).
\]	 
Then,
\begin{align*}
	\Rect(S) & = \begin{ytableau} 1 & 2 & 5 & 8 \\ 3 & 6 \\ 4 & 7 \end{ytableau} \, , \qquad \Rect^\ast(S) + 8 = \begin{ytableau} \none & \none & \none & \none \\ \none & \none & *(yellow) 10 & *(yellow) 13 \\ \none & \none & *(yellow) 11 & *(yellow) 14 \\ *(yellow) 9 & *(yellow) 12 & *(yellow) 15 & *(yellow) 16 \end{ytableau} \, , \\
	R_2(S) & = \begin{ytableau} 1 & 2 & 5 & 8 \\ 3 & 6 & *(yellow) 10 & *(yellow) 13 \\ 4 & 7 & *(yellow) 11 & *(yellow) 14 \\ *(yellow) 9 & *(yellow) 12 & *(yellow) 15 & *(yellow) 16 \end{ytableau} \, .
\end{align*}

\end{Example}

%

We use the name $ R_2(S) $ because the construction is similar to the definition of $ R_a(S) $ if we set $ a = 2 $. Note that we cannot just set $ a =2 $ in \Cref{def:Raconstruct} because \Cref{def:Raconstruct} requires $ a \ge 2 \stab(S) - 1 $. Similar to \Cref{thm:fixedbyprom}, we will show $ R_2(S) $ has shape $ ((2r)^b) $, is a standard tableau, and is fixed by $ \pr^{br} $.

\begin{Lemma} \label{lem:symrect} The map 
\[
	\Rect: \SYT(\g) \to \bigcup_{ \substack{ \nu = ( \nu_1, \dots, \nu_b) \, \vdash \, br, \\ \tx{s.t.} \, \nu_j + \nu_{b + 1 - j} = 2r } } \SYT(\nu) 
\] 
is bijective.

\end{Lemma}

\begin{proof} Say $ \g = \lam/\mu $ as a skew shape. By \cite[Lemma 3.3]{MR859302}, we have the Littlewood-Richardson coefficient
\begin{align}
\label{eq:symLRCs}
	c^{\lam}_{\mu, \nu} = \begin{cases} 1, \quad &\tx{ if } \nu_j + \nu_{b + 1 - j} = 2r \tx{ for all } j \in [1,b], \\
							0, \quad & \tx{ else,} \end{cases}
\end{align}
for all $ \nu =  (\nu_1, \dots, \nu_b)  \vdash br $. \Cref{lem:symrect} follows immediately from \eqref{eq:symLRCs} and \Cref{def:LRCs}.

\end{proof}

\begin{proof}[Proof of \Cref{thm:fixedbyab/2prom}] \cite{Whiteprom} \cite{RheeThesis} Fix $ S \in \SYT(\g) $. For all $ j = 1, \dots, b $ and 
\begin{align*}
	\rv(R_2(S)) & = \rv(\Rect(S))_j + \rv(\Rect^\ast(S))_{j} \\
			& = \rv(\Rect(S))_j + \rv(e(\Rect(S))^\dagger)_{j} \qquad \tx{ by \Cref{lem:evacflip} } & \\
			& = \rv(\Rect(S))_j + \rv(\Rect(S))_{b + 1 - j} \\
			& = 2r  \qquad \tx{ by \Cref{lem:symrect} } &.
\end{align*}
Indeed, $ R_2(S) $ is a filling of shape $ ((2r)^b) $. The fillings $ \left. R_2(S) \right|_{[1,br]} = \Rect(S) $ and $ \left. R_2(S) \right|_{[br + 1, 2br]} = \Rect^\ast(S) + br $ are skew tableaux using the entries $ [1,br] $ and $ [br + 1, 2br] $, repsectively. Thus, their row-concatenation from left to right, $ R_2(S) $, is a standard tableau. 

Similar to the proof of \Cref{thm:fixedbyprom},
\begin{equation}
\begin{aligned}
\label{eq:prfixedab/2rest1}
	\left. \pr^{br}(R_2(S)) \right|_{[1, br]} & = \Rect(\Rect^\ast(S) + br) - br \qquad \tx{ by } \eqref{eq:multprom1} \\
						   & = \Rect(S) \qquad \tx{ as rectification is well-defined} \\
						   & = \left. R_2(S) \right|_{[1,br]} \qquad \tx{ by \Cref{def:R2construct} },
\end{aligned}
\end{equation}
and
\begin{equation}
\begin{aligned}
\label{eq:prfixedab/2rest2}
	\left. \pr^{br}(R_2(S)) \right|_{[br + 1, 2br]} & = \Rect^\ast(\Rect(S)) + br \qquad \tx{ by } \eqref{eq:multprom2} \\
						   & = \Rect^\ast(S) + br \qquad \tx{ as anti-rectification is well-defined} \\
						   & = \left. R_2(S) \right|_{[br + 1, 2br]} \qquad \tx{ by \Cref{def:R2construct} }. 
\end{aligned}
\end{equation}
Putting \eqref{eq:prfixedab/2rest1} and \eqref{eq:prfixedab/2rest2} together gives $ \pr^{br}(R_2(S)) = R_2(S) $. Thus,
\begin{align}
\label{eq:prab/2sub}
	\{ R_2(S) : S \in \SYT(\g) \} \sube \SYT((2r)^b)^{\pr^{br}}.
\end{align}

The fact that $ \{ R_2(S) : S \in \SYT(\g) \} $ is a set and not a multiset is a consequence of \Cref{lem:symrect}. By \eqref{eq:2quotientofrectangle}, both sides of \eqref{eq:prab/2sub} have the same size, so 
\[
	 \SYT((2r)^b)^{\pr^{br}} = \{ R_2(S) : S \in \SYT(\g) \}.
\] 

Finally, we count these tableaux by first partitioning $ \{ 1, 2, \dots, br \} $ into two blocks of size $ r \dd \left \lceil \f{b}{2} \right \rceil $ and $ r \dd \left \lfloor \f{b}{2} \right \rfloor $. Then, we choose a filling of $ ( r^{\left \lceil \f{b}{2} \right \rceil} ) $ with the first block and a filling of $ ( r^{\left \lfloor \f{b}{2} \right \rfloor} ) $ with the second block so that both fillings have increasing rows and columns. This yields
\[
	\# \SYT((2r)^b)^{\pr^{br}} = \ch{br}{ \left \lfloor \f{b}{2} \right \rfloor r} \# \SYT( r^{\left \lceil \f{b}{2} \right \rceil}) \dd \# \SYT( r^{\left \lfloor \f{b}{2} \right \rfloor}).
\] 
\end{proof}

We also describe the action of promotion on $ \SYT((2r)^b)^{\pr^{br}} $, which is closed under promotion. Like $ R_a $ in \Cref{cor:prfixedpr}, the $ R_2 $ operator also commutes with promotion.


\begin{Corollary} \cite{Whiteprom} \cite{RheeThesis} \label{cor:prfixedpr2} For any $ S \in \SYT(\g) $,
\[
	\pr(R_2(S)) = R_2(\pr(S)).
\]
\end{Corollary}

\begin{proof} Similar to \Cref{cor:prfixedpr}, $ \SYT((2r)^b)^{\pr^{br}} $ is closed under promotion. Thus, fixing $ S \in \SYT(\g) $, 
\[
	p(R_2(S)) = R_2(S') \qquad \tx{ for some } S' \in \SYT(\g).
\]
Our reasoning from the proof of \Cref{cor:prfixedpr} still holds through \eqref{eq:promprfixedrest}, so
\begin{align}
\label{eq:promprfixedrest2}
	\left. R_2(S') \right|_{[1,br - 1]} = \left. R_2(p(S)) \right|_{[1, br - 1]}.
\end{align}
By \Cref{lem:symrect}, $ \left. R_2(S') \right|_{[1,br]} $ and $ \left. R_2(\pr(S)) \right|_{[1,br]} $ each have shapes $ \nu = (\nu_1, \dots, \nu_b) $ satisfying $ \nu_j + \nu_{b + 1 - j} = 2r $, which together with \eqref{eq:promprfixedrest2}, forces $ br $ to be in the same cell in $ R_2(S') $ and $ R_2(\pr(S)) $. Thus, 
\begin{align}
\label{eq:promprfixedrest3}
	\left. \Rect(S') = R_2(S') \right|_{[1,b r]} = \left. R_2(\pr(S)) \right|_{[1,b r]} = \Rect(\pr(S)).
\end{align}
Finally, \Cref{lem:symrect} forces $ S' = \pr(S) $.  
\end{proof}

For any $ k \mid br $, 
\[
	\SYT((ar)^b)^{\pr^{br/k}} \sube \SYT((ar)^b)^{\pr^{br}}.
\]
By \Cref{cor:prfixedpr} and \Cref{cor:prfixedpr2}, we have, for $ a = 2 $ or $ a \ge 2b - 1 $,
\begin{align}
\label{eq:prfixedwithin}
	\SYT((ar)^b)^{\pr^{br/k}} = \{ R_a(S) : S \in \SYT(Q_a((ar)^b)), \pr^{br/k}(S) = S \}.
\end{align}
For $ a = 2 $, we will describe $ \{ S \in \SYT(\g): \pr^{br/k}(S) = S \} $ in terms of rectangular tableaux fixed by smaller promotion powers. 

%
%
%
%
%

\begin{Definition} \label{def:nonstandardprom} For a standard tableau $ T $ of size $ n $, $ m \in \bZ_{\ge 1} $, and $ A = \{ a_1, \dots, a_k \} \sube [m] $ with $ a_1 < \dots < a_k $ and $ k \mid n $, let $ I_{m,A}(T) $ denote the skew tableau obtained from $ T $ by replacing $ ik + j $ by $ im + a_j $ for all $ i = 0,1, \dots, n/k - 1 $ and $ j = 1, \dots, k $. For example, if 
\[
	T = \begin{ytableau} 1 & 2 & *(yellow) 5 & *(green) 8 \\ 3 & *(yellow) 6 & *(green) 9 & *(orange) 11 \\ *(yellow) 4 & *(green) 7 & *(orange) 10 & *(orange) 12 \end{ytableau} \,, 
\]
then
\[
	I_{6, \{2, 4, 5 \}} (T) = \begin{ytableau} 2 & 4 & *(yellow) 10 & *(green) 16 \\ 5 & *(yellow) 11 & *(green) 17 & *(orange) 22 \\ *(yellow) 8 & *(green) 14 & *(orange) 20 & *(orange) 23 \end{ytableau} \,.
\]

\end{Definition}

\begin{Definition} On a tableau with entries $ i_1 < \dots < i_n $, let promotion act on the indices as it would on a standard tableau. For example, since 
\begin{align*}
	\pr & : \begin{ytableau} 1 & 3 & 5 & 6 \\ 2 & 4 & 8 & 11 \\ 7 & 9 & 10 & 12 \end{ytableau} \mapsto \begin{ytableau} 1 & 2 & 4 & 5 \\ 3 & 7 & 9 & 10 \\ 6 & 8 & 11 & 12 \end{ytableau} \, , \\
	\pr & : \begin{ytableau} i_1 & i_3 & i_5 & i_6 \\ i_2 & i_4 & i_8 & i_{11} \\ i_{7} & i_9 & i_{10} & i_{12} \end{ytableau} \mapsto \begin{ytableau} i_1 & i_2 & i_4 & i_5 \\ i_3 & i_7 & i_9 & i_{10} \\ i_6 & i_8 & i_{11} & i_{12} \end{ytableau} \quad \tx{ for any } i_1 < \dots < i_{12}. 
\end{align*}
We realize there are other variations on how promotion acts on tableaux with distinct non-standard entries, such as in \cite{MR2557880}, but this will prove convenient for our purposes.
\end{Definition}

\begin{Thm} \label{thm:prfixedpieces} (inspired by \cite{Whiteprom}) For any partitions $ \lam^{(1)}, \dots, \lam^{(b)} $ with total size $ n $, and $ k \mid n $,
\begin{equation}
\begin{aligned}
\label{eq:prfixedpieces}
	\SYT & (\lam^{(1)} \cup \dots \cup \lam^{(b)})^{\pr^{n/k}} = \lb I_{n/k, M_1}(T_1) \cup \dots \cup I_{n/k, M_b}(T_b) : \right. \\ 
	& \left. T_j \in \SYT(\lam^{(j)})^{\pr^{ |\lam^{(j)}|/k }}, (M_1, \dots, M_b) \in \ch{[ n/k]}{ |\lam^{(1)}|/k , \dots, |\lam^{(b)}|/k} \rb
\end{aligned}
\end{equation}
if $ k \mid \gcd( |\lam^{(1)}|, \dots, |\lam^{(b)}|) $. Else, 
\[
	\SYT(\lam^{(1)} \cup \dots \cup \lam^{(b)})^{\pr^{n/k}} = \varnothing.
\]
\end{Thm}

\begin{Example} \label{ex:prfixedinpieces} Suppose $ b = 2, \lam^{(1)} = (3,3), \lam^{(2)} = (2,2), k = 2 $. Let
\begin{align*}
	T_1 = \begin{ytableau} 1 & 2 & *(yellow) 5 \\ 3 & *(yellow) 4 & *(yellow) 6 \end{ytableau} \in \SYT(3,3)^{\pr^3}, \; & \; T_2  = \begin{ytableau} 1 & *(yellow) 3 \\ 2 & *(yellow) 4 \end{ytableau} \in \SYT(2,2)^{\pr^2}, \\
	(M_1, M_2) = (\{ 1, 4, 5 \} &, \{ 2, 3 \}) \in \ch{[5]}{3,2}.
\end{align*}
Then,
\begin{align*}
	I_{5,M_1}(T_1) & = \begin{ytableau} 1 & 4 & *(yellow) 9 \\ 5 & *(yellow) 6 & *(yellow) 10 \end{ytableau} \, , \quad I_{5,M_2}(T_2) = \begin{ytableau} 2 & *(yellow) 7 \\ 3 & *(yellow) 8 \end{ytableau} \, , \\
	I_{5, M_1}(T_1) \cup I_{5, M_2}(T_2) & = \begin{ytableau} \none & \none & \none & 2 & *(yellow) 7 \\ \none & \none & \none & 3 & *(yellow) 8 \\ 1 & 4 & *(yellow) 9 \\ 5 & *(yellow) 6 & *(yellow) 10 \end{ytableau} \in \SYT((3,3) \cup (2,2))^{\pr^{5}}.
\end{align*}
\end{Example}

\begin{proof} Suppose $ \lam^{(1)}, \dots, \lam^{(b)} $ are partitions with total size $ n $ and $ k \mid \gcd( |\lam^{(1)}|, \dots, |\lam^{(b)}|) $. Let
\[
	m \coloneqq \f{n}{k}, \qquad m_j \coloneqq \f{ |\lam^{(j)}| }{k} \quad \tx{ for all } j = 1, \dots, b.
\]
Consider $ T_j \in \SYT(\lam^{(j)})^{\pr^{m_j}} $ for $ j = 1, \dots, b $ and $ (M_1, \dots, M_b) \in \ch{[m]}{m_1, \dots, m_b} $. Let 
\begin{align}
\label{eq:defSj}
	S_j \coloneqq I_{m,M_j}(T_j) \qquad \tx{ for all } j = 1, \dots, b,
\end{align}
and
\[
	S \coloneqq S_1 \cup \dots \cup S_b.
\]
Also, write
\begin{align}
\label{eq:prbypieces}
	\pr^m(S) = \pr^m(S)_1 \cup \dots \cup \pr^m(S)_b.
\end{align}
For each $ j $, the entries in $ S_j $ are $ \cup_{i = 1}^{k} ( M_j + (i - 1)m ) $. Thus, since $ \pr^m $ performs a total of $ m $ decrements modulo $ n $, the entries in $ \pr^m(S)_j $ are
\[
	\bigcup_{\ell = 1}^{k - 1} ( M_j + (\ell - 1)m) \bigcup (M_j + n - m) = \bigcup_{\ell = 1}^{k} ( M_j + (\ell - 1)m ).  
\] 
Thus, the same set of entries appears in $ S_j $ and $ \pr^m(S)_j $. By this and the fact that $ m_j $ of the $ m $ promotions on $ S $ slide through $ S_j $,
\begin{align}
\label{eq:prrestrict}
	\pr^m(S)_j = \pr^{m_j}(S_j).
\end{align}
Now,
\begin{equation}
\label{eq:pronpiece}
\begin{aligned}
	 \pr^{m_j}(S_j) & = \pr^{m_j}(I_{m,M_j} (T_j)) \qquad \tx{ by } \eqref{eq:defSj} \\
		 & = I_{m,M_j}(\pr^{m_j}(T_j)) \qquad \tx{ by \Cref{def:nonstandardprom} } \\
		& = I_{m, M_j}(T_j) \qquad \tx{ because $ T_j \in \SYT(\lam^{(j)})^{\pr^{m_j}} $} \\
		& = S_j \qquad \tx{ by \eqref{eq:defSj} }
\end{aligned}
\end{equation}
for all $ j = 1, \dots, b $. Combining \eqref{eq:prbypieces}, \eqref{eq:prrestrict}, and \eqref{eq:pronpiece} yields $ \pr^m(S) = S $, so $ \supseteq $ holds in \eqref{eq:prfixedpieces}.

\ms

On the other hand, suppose $ S \in \SYT (\lam^{(1)} \cup \dots \cup \lam^{(b)})^{\pr^m} $. Let $ S = S_1 \cup \dots \cup S_b $, $ M_j $ denote the entries in $ S_j $ that lie in $ [m] $ and $ m_j \coloneqq \# M_j $. Since $ \pr^m(S) = S $, the set of entries in $ S_j $ must be fixed by $ m $ decrements modulo $ n $. When applying $ \pr^m $ to $ S $, we apply $ m_j $ inner slides and $ m $ decrements modulo $ n $ to $ S_j $, so letting $ \{ i_1, \dots, i_{ |\lam^{(j)}| } \} $ denote the set of entries of $ S_j $ in increasing order, we must have
\begin{align}
\label{eq:listequality}
	i_{t + m_j \, (\Mod |\lam^{(j)}|)} \equiv_n i_{t} + m \qquad \tx{ for all } t = 1, \dots, |\lam^{(j)}|.
\end{align}
Equation \eqref{eq:listequality} means that shifting the indices by $ m_j $ modulo $ |\lam^{(j)}| $ corresponds to shifting the entries by $ m $ modulo $ n $. This means $ |\lam^{(j)}| = k_j m_j $ for some $ k_j \in \bZ_{\ge 1} $, $ M_j = \{ i_1, \dots, i_{m_j} \} $, and that the set of entries in $ S_j $ is
\[
	\bigcup_{\ell = 1}^{k_j} (M_j + (\ell - 1)m).
\] 
As $ S \in \SYT(\lam^{(1)} \cup \dots \cup \lam^{(b)}), $ $ S $ must use exactly the entries $ 1, 2, \dots, n $, forcing 
\[
	k_1 = \dots = k_b \coloneqq k.
\]
Now we have
\[
	k \mid \gcd( |\lam^{(1)}|, \dots, |\lam^{(b)}|), \quad m = \f{n}{k}, \quad m_j = \f{ |\lam^{(j)}| }{k} \tx{ for all } j = 1, \dots b.
\]
Hence, $ \SYT(\lam^{(1)} \cup \dots \cup \lam^{(b)})^{\pr^{n/k}} = \varnothing $ unless $ k \mid \gcd( |\lam^{(1)}|, \dots, |\lam^{(b)}|) $.

\ms

Now, assume that $ k \mid \gcd( |\lam^{(1)}|, \dots, |\lam^{(b)}|) $. Since the set of entries in $ S_j $ is $ \bigcup_{\ell = 1}^{k_j} (M_j + (\ell - 1)m) $, we can write $ S_j = I_{m,M_j}(T_j) $ for some $ T_j \in \SYT(\lam^{(j)}) $. Then,
\begin{align*}
	\pr^m(S) & = S \\
	\imp \pr^{m_j}(S_j) & = S_j \qquad \tx{ by \eqref{eq:prbypieces} and \eqref{eq:prrestrict} } \\
	\imp \pr^{m_j}( I_{m,M_j}(T_j)) & = I_{m,M_j}(T_j) \\
	\imp I_{m,M_j}( \pr^{m_j}(T_j) ) & = I_{m,M_j}(T_j) \qquad \tx{ by \Cref{def:nonstandardprom} } \\
	\imp \pr^{m_j}(T_j) & = T_j
\end{align*}
for all $ j = 1, \dots, b $. This proves that $ \sube $ holds in \eqref{eq:prfixedpieces} and hence proves \Cref{thm:prfixedpieces}.
 
\end{proof}

\begin{Corollary} \label{cor:R2prfixedwithin} For any $ r, b \in \bZ_{\ge 1} $, and $ k \mid br $,
\begin{align*}
	\SYT ((2r)^b)^{\pr^{br/k}} & = \left \{ R_2(I_{br/k, M_1}(T_1) \cup I_{br/k, M_2}(T_2)) : \right. \left. T_1 \in \SYT \lp r^{\left \lceil \f{b}{2} \right \rceil} \rp^{\pr^{r \left \lfloor \f{b}{2} \right \rfloor / k}} \right. , \\
	& \left. T_2 \in \SYT \lp r^{\left \lceil \f{b}{2} \right \rceil} \rp^{\pr^{r \left \lfloor \f{b}{2} \right \rfloor /k}}, (M_1, M_2) \in \ch{[br/k]} {r \left \lceil \f{b}{2} \right \rceil / k, r \left \lfloor \f{b}{2} \right \rfloor / k} \rb
\end{align*}
if $ k \mid \gcd(r^{\left \lfloor \f{b}{2} \right \rfloor}, r^{\left \lceil \f{b}{2} \right \rceil}) $. Else, 
\[
	\SYT ((2r)^b)^{\pr^{br/k}} = \varnothing.
\]
\end{Corollary}

\begin{proof} By \Cref{thm:fixedbyab/2prom} and \eqref{eq:prfixedwithin}, 
\begin{align}
\label{eq:R2prfixed}
	\SYT( (2r)^{b})^{\pr^{br/k}} & = \lb R_2(S) : S \in \SYT(r^{\left \lfloor \f{b}{2} \right \rfloor} \cup r^{\left \lceil \f{b}{2} \right \rceil})^{\pr^{br/k}} \rb.
\end{align}
\Cref{cor:R2prfixedwithin} follows from \eqref{eq:R2prfixed} and \Cref{thm:prfixedpieces}.
\end{proof}

\begin{Example} Continuing \Cref{ex:prfixedinpieces},
\[
	R_2  \begin{ytableau} \none & \none & \none & 2 & *(yellow) 7 \\ \none & \none & \none & 3 & *(yellow) 8 \\ 1 & 4 & *(yellow) 9 \\ 5 & *(yellow) 6 & *(yellow) 10 \end{ytableau} = \begin{ytableau} 1 & 2 & *(yellow) 7 & *(green) 12 & *(orange) 17 \\ 3 & *(yellow) 6 & *(yellow) 8 & *(green) 13 & *(orange) 18 \\ 4 & *(yellow) 9 & *(green) 11 & *(green) 14 & *(orange) 19 \\ 5 & *(yellow) 10 & *(green) 15 & *(orange) 16 & *(orange) 20 \end{ytableau} \in \SYT(4^5)^{\pr^5}.
\]
\end{Example}

\begin{Remark} \label{rem:prfixedwithin} One can generalize \Cref{cor:R2prfixedwithin} to give a similar description of $ \SYT((ar)^{b})^{\pr^{br/k}} $ for $ a \ge 2b - 1, k \ge 2 $. However, such a generalization would not give us anything new.
In order for $ \SYT((ar)^b)^{\pr^{br/k}} $ to be nonempty, we must have $ b \mid \f{br}{k} $ or $ ar \mid \f{br}{k} $ by \eqref{eq:prfixednonempty}. Since $ ar \ge (2b - 1)r > \f{br}{k} $, $  ar \not| \f{br}{k} $, so
\[
	b \mid \f{br}{k} \imp k \mid r.
\]
Note also that $ ak > a \ge 2b - 1 $. Then, by \Cref{thm:prfixedtabs},
\begin{align*}
	\SYT((ar)^b)^{\pr^{br/k} } & = \SYT \lp \lp ak \dd \f{r}{k} \rp^b \rp^{\pr^{b \dd \f{r}{k}} }
	= \lb R_{ak}(\ti{S}) : \ti{S} \in \SYT ( \underbrace{(r/k) \cup \dots \cup (r/k) }_{b \tx{ times} } ) \rb.
\end{align*}
These tableaux were already constructed in \Cref{thm:prfixedtabs}.
\end{Remark}

\section{Stabilization as a Permutation Statistic}
\label{sec:stabonperms}

We can identify each permutation $ w \in S_n $ with the unique skew tableau with shape $ (1) \cup \dots \cup (1) $ whose reading word is $ w $. Under this identification, we translate the stabilization statistic to the symmetric group. We realize stabilization is invariant on dual equivalence classes and is bounded below strictly by the number of ascents of $ w $. We characterize the permutations with stabilization statistic 1 and 2 in terms of their recording tableaux. 

\begin{Definition} For $ w = w_1 \dots w_n \in S_n $, let $ T(w) $ be the unique standard skew tableau with shape $ (1) \cup \dots \cup (1) $ with reading word $ w $, and define $ \stab(w) \coloneqq \stab(T(w)) $. Let 
\[
	\asc(w) \coloneqq \# \{ i : w_i < w_{ i + 1} \} 
\]
denote the number of ascents of $ w $.

\begin{Example}
\begin{align*}
	\stab(4231) & = \stab \begin{ytableau} \none & \none & \none & 1 \\ \none & \none & 3 \\ \none & 2 \\ 4 \end{ytableau} =  3, \qquad \stab(4123) = \begin{ytableau} \none & \none & \none & 3 \\ \none & \none & 2 \\ \none & 1 \\ 4 \end{ytableau} = 3, \\
	\asc(4231) & = 1, \qquad \qquad \qquad \qquad \qquad \qquad  \asc(4123) = 2. 
\end{align*}
\end{Example}

\end{Definition}

%

Recall from \eqref{eq:DEQ} that $ v, w \in S_n $ are dual equivalent if and only if $ Q(v) = Q(w) $. As in \Cref{thm:dualequivstab}, dual equivalence plays well with stabilization.

\begin{Lemma} \label{lem:dualpermstab} If $ v, w \in S_n $ are dual equivalent, then $ \stab(v) = \stab(w) $.

\end{Lemma}

\begin{proof} Since $ v, w $ are dual equivalent, $ T(v) $ and $ T(w) $ are dual equivalent, see \cite[Lemma~2.11]{MR1158783}. It follows by \Cref{thm:dualequivstab} that 
\[
	\stab(v) = \stab(T(v)) = \stab(T(w)) = \stab(w).
\]
\end{proof}

\begin{Lemma} \label{lem:ascentbound} For all $ w \in S_n $, $ \stab(w) > \asc(w) $. 

\end{Lemma}

\begin{proof} Suppose $ w = w_1 \dots w_n \in S_n $, and let $ r = \stab(w) $, so $ T(w) $ stabilizes at $ r $. In particular, this means $ \Rect(T(w)^{(r)}) $ has $ n $ rows. The reading word of $ T(w)^{(r)} $ is
\begin{align*}
	w_1^{(r)} \dots w_n^{(r)} = w_1 \, (w_1 + n) \, \dots \, (w_1 + (r - 1)n) \, \dots \, w_n \, (w_n + n) \, \dots \, (w_n + (r - 1)n),
\end{align*}
so
\begin{align}
\label{eq:Pwr}	
	\Rect(T(w)^{(r)}) = P(w_1^{(r)} \dots w_n^{(r)})
\end{align}
by \Cref{lem:rectP}. By \Cref{thm:Greene}, $ w_1^{(r)} \dots w_n^{(r)} $ must have a decreasing subsequence of size $ n $. Yet, at most one of $ w_j, (w_j + n), \dots, (w_j + (r - 1)n) $ can be in such a decreasing subsequence. Thus, $ w_1^{(r)} \dots w_n^{(r)} $ must have a decreasing subsequence of the form
\[
	(w_1 + (r - 1)n) \dots (w_{i_1} + (r - 1)n) \, (w_{i_1 + 1} + (r - 2)n) \dots (w_{i_{2}} + (r - 2)n) \dots (w_{i_{r -1} + 1}) \dots (w_n)
\]
for some $ 1 \le i_1 \le i_2 \dots \le i_{r - 1} \le n $. Since this subsequence is decreasing,
\[
	w_1 > \dots > w_{i_1}, \quad w_{i_1 + 1} > \dots > w_{i_2}, \quad \dots, \quad w_{i_{r - 1} + 1} > \dots > w_n.
\]
Therefore, $ w $ can only have ascents at possibly $ i_1, \dots, i_{r - 1} $, so $ w $ has at most $ r - 1 $ ascents. Hence, $ \stab(w) = r > r - 1 \ge \asc(w) $.

\end{proof}

\Cref{lem:dualpermstab} and \Cref{lem:ascentbound} give us some information about the distribution of $ \stab $ on $ S_n $, but what else can we say about this distribution? By \Cref{thm:stab}, we know $ 1 \le \stab(w) \le n $ for $ w \in S_n $. In \Cref{fig:stabtri}, we give the distribution of $ \stab $ on $ S_1, \dots, S_8 $. We have a formula for the number of permutations in $ S_n $ with $ \stab $ 1 or $ \stab $ 2. In fact, we characterize exactly which permutation has $ \stab $ 1 and which permutations have $ \stab 2 $ in terms of their recording tableaux.

\begin{figure}[h]
      \centering
\begin{align*}
	1& \\
1 \quad & \quad 1 \\
1 \quad \; 4 & \quad \quad 1 \\
1 \quad 8 \quad & \; \, 14 \quad \; 1 \\
1 \quad 18 \quad 6&3 \quad 37 \quad 1 \\
1 \quad \, 33 \quad \, 175 \; & \; 434 \quad 76 \quad 1 \\
1 \quad \,  68 \quad 549 \; \; 23&45 \; \, 1927 \; \, 149 \; 1 \\
1 \quad 124 \; \; 1787 \; \; 7807 \, & \, 23760 \; \, 6552 \; \; 288 \; \; 1
\end{align*}
\caption{The distribution of $ \stab $ on $ S_1, \dots, S_8 $. The $ k $-th entry from the left in row $ n $ is $ \# \{ w \in S_n : \stab(w) = k \} $.}
\label{fig:stabtri}
\end{figure}

\begin{Lemma} The permutation $ w = n \, (n - 1) \, \dots \, 2 \, 1 $ is the only permutation with $ \stab(w) = 1 $.

\end{Lemma}

\begin{proof} If $ w = n \, (n - 1) \, \dots \, 2 \, 1 $, then $ \Rect(T(w)) = P(w) $ has a single column, so $ \stab(w) = 1 $. If $ w \in S_n $ with $ \stab(w) = 1 $, then $ P(w) $ consists of $ n $ rows, forcing $ w = n \, (n - 1) \, \dots \, 2 \, 1 $.

\end{proof}

\begin{Notation} Fix $ n \in \bZ_{\ge 1} $. For $ k = 1, \dots, n $, let
\[
	\Yvcentermath1 T_k = \ytableausetup{boxsize = 1.8em} \begin{ytableau} 1 &  \scriptstyle{k + 1} \\ \vdots & \vdots \\ k & 2k \\ x \\ \vdots \\ n \end{ytableau} \, , \qquad \tx{ for } k < \left\lfloor \f{n}{2} \right\rfloor, \qquad T_k = \ytableausetup{boxsize = 1.8em} \begin{ytableau} 1 & \scriptstyle{k + 1} \\ \vdots & \vdots \\ \scriptstyle{n - k} & n \\ y \\ \vdots \\ k \end{ytableau} \, , \qquad \tx{ for } k \ge  \left\lfloor \f{n}{2} \right\rfloor.
\]
where $ x = 2k + 1, y = n - k + 1 $. This also means
\begin{align}
\label{eq:PQstab2}
	T_k = P( k \, \dots \, 1 \, n \, \dots \, (k + 1)) = Q ( k \, \dots \, 1 \, n \, \dots \, (k + 1)).
\end{align}
\end{Notation}

\begin{Thm} \label{thm:stab2} For all $ n \in \bZ_{\ge 1} $,
\begin{align}
\label{eq:stab2perms}
	\{ w \in S_n : \stab(w) = 2 \} = \{ w \in S_n : Q(w) = T_k \tx{ for some k} \}.
\end{align}
Consequently,
\begin{align}
\label{eq:numstab2perms}
	\# \{ w \in S_n : \stab(w) = 2 \} = \ch{n + 1}{ \lfloor \f{n + 1}{2} \rfloor} - 2.
\end{align}
See OEIS entry \href{https://oeis.org/A201686}{A201686}.
\end{Thm}

We break the proof of \Cref{thm:stab2} into 3 steps. First, we use \Cref{lem:ascentbound} to prove \eqref{eq:stab2perms}. Secondly, we calculate the number of standard tableaux of a given size with 1 or 2 columns. Thirdly, we use this result to prove \eqref{eq:numstab2perms}. 

\begin{proof}[Proof of \eqref{eq:stab2perms}] Suppose $ w \in S_n $ has $ \stab(w) = 2 $. By \Cref{lem:ascentbound}, $ \asc(w) \le 1 $, so we can write $ w = w_1 \dots w_k w_{k + 1} \dots w_n $ where
\[
	w_1 > \dots > w_k, \quad w_{ k + 1} > \dots > w_n
\]
for some $ k = 1, \dots, n $. We must also have $ w_k < w_{k + 1} $ or else $ w = n \, (n - 1) \, \dots \, 2 \, 1 $, which has $ \stab(w) = 1 $. Thus, $ w_k < w_{ k + 1} $. Since $ \stab(w) = 2 $, $ \Rect(T(w)^{(2)}) $ must have $ n + 1, \dots, 2n $ in distinct rows. Note that we can perform inner slides to $ T(w)^{(2)} $ to get
\[
	\ytableausetup{boxsize = 1.8em} \begin{ytableau} \none & w_n & w_n' \\ \none & \vdots & \vdots \\ \none & \scriptstyle{w_{ k + 1}} & \scriptstyle{w_{ k + 1}'} \\ w_k & w_k' \\ \vdots & \vdots \\ w_1 & w_1' \end{ytableau} \, ,
\]
where $ x' = x + n $ for all $ x $. In order for $ T(w) $ to stabilize at 2, $ w_1', \dots, w_n' $ can only slide horizontally. Since $ w_ k < w_{ k + 1} $ and thus $ w_k' < w_{ k + 1}' $, $ w_k' $ will slide vertically if the cell above it is vacated while it is in the second column. Since any entry above $ w_k' $ sliding into column 1 forces $ w_k' $ to slide up, if $ w_k' $ slides into column 1, it must do so before any of the entries above it. Hence, $ w_k, \dots, w_1 $ must slide up either above $ w_k' $ or to the top before $ w_n, \dots, w_{k + 1} $ experience any horizontal slides.

If $ k \ge\left \lfloor \f{n}{2} \right \rfloor $, then $ w_k, \dots, w_1 $ slide to the top and $ w_n, \dots, w_{k + 1} $ stay still, so
\[
	P(w^{(2)}) = \Rect(T(w)^{(2)}) = \ytableausetup{boxsize = 2em} \begin{ytableau} w_k & w_n & *(yellow) w_n' \\ \vdots & \vdots & *(yellow) \vdots \\ x & \scriptstyle{w_{k + 1}} & *(yellow) \scriptstyle{w_{ k + 1}'} \\ \scriptstyle{w_{2k - n}} & *(yellow) w_k' \\ \vdots & *(yellow) \vdots \\ w_1 & *(yellow) y' \\ *(yellow) \scriptstyle{w_{n - k}'} \\ *(yellow) \vdots \\ *(yellow) w_1' \end{ytableau} \quad \imp \quad
	P(w) = \begin{ytableau} w_k & w_n \\ \vdots & \vdots \\ x & \scriptstyle{w_{k + 1}} \\ \scriptstyle{w_{2k - n}} \\ \vdots \\ w_1 \\ \end{ytableau} \, .
\]
where $ x = w_{2k - n + 1}, y' = w_{n - k + 1}' $. Since $ w = w_1 \dots w_k w_{k + 1} \dots w_n $, this means during RSK on $ w $, the whole first column of $ P(w) $ was created first. Thus,
\[
	\Yvcentermath1 Q(w) = \ytableausetup{boxsize = 1.8em} \begin{ytableau} 1 & \scriptstyle{k + 1} \\ \vdots & \vdots \\ \scriptstyle{n - k} & n \\ \vdots \\ k \end{ytableau} = T_k.
\]

If $ k < \left \lfloor \f{n}{2} \right \rfloor$, then $ w_k, \dots, w_1 $ slide up above $ w_k' $ before $ w_n, \dots, w_{k + 1} $ experience any horizontal slides, so
\[
	P(w^{(2)}) = \Rect(T(w)^{(2)}) = \Rect \ytableausetup{boxsize = 2em} \begin{ytableau} \none & w_n & *(yellow) w_n' \\ \none & \vdots & *(yellow) \vdots \\ \none & \scriptstyle{w_{2k + 1}} & *(yellow) \scriptstyle{w_{2k + 1}'} \\ w_k & w_{2k} & *(yellow) w_{2k}' \\ \vdots & \vdots & *(yellow) \vdots \\ w_1 & \scriptstyle{w_{k + 1}} & *(yellow)  \scriptstyle{w_{k + 1}'} \\ *(yellow) w_{k}' \\ *(yellow) \vdots \\ *(yellow) w_1' \end{ytableau}
\]
In particular, based on the position of $ w_1 \dots w_{2k} $, $ \RSK(w_1 \dots w_{2k}) $ began with a column of size k followed by a second column of size $ k $. Thus,
\[
	P(w_1 \dots w_{2k}) = \ytableausetup{boxsize = 2em} \begin{ytableau} w_k & w_{2k} \\ \vdots & \vdots \\ w_1 & \scriptstyle{w_{k + 1}} \end{ytableau} \, , \qquad Q(w_1 \dots w_{2k}) = \begin{ytableau} 1 & \scriptstyle{k + 1} \\ \vdots & \vdots \\ k & 2k \end{ytableau} \, .
\]
As $ w $ has a decreasing subsequence of size $ n - k $, the first column of $ Q(w) $ must have size at least $ n - k $ by \Cref{thm:Greene}. As $ Q(w) $ contains $ Q(w_1 \dots w_{2k}) $, we must have
\[
	Q(w) = \begin{ytableau} 1 & \scriptstyle{k + 1} \\ \vdots & \vdots \\ k & 2k \\ \scriptstyle{2k + 1} \\ \dots \\ n \end{ytableau} = T_k.
\]
This shows that if $ \stab(w) = 2 $, then $ Q(w) = T_k $ for some $ k $, and hence
\begin{align}
\label{eq:stab2permscont}  
	\{ w \in S_n : \stab(w) = 2 \} \sube \{ w \in S_n : Q(w) = T_k \tx{ for some k} \}.
\end{align}

In order to show the reverse containment of \eqref{eq:stab2permscont} and thus complete the proof of \eqref{eq:stab2perms}, suppose $ Q(w) = T_k $ for some $ k $. To show $ \stab(w) = 2 $, it suffices to verify $ \stab(w) = 2 $ for only 1 such word $ w $ in each dual equivalence class by \Cref{lem:dualpermstab}. As the recording tableau characterizes dual equivalence classes, it suffices to check that
\[
	\stab( k  \dots 1 \, n \dots (k + 1) ) = 2 \quad \tx{ for all } k = 1, \dots n.
\]
by \eqref{eq:PQstab2}. Let $ w = k \dots 1 \, n \dots (k + 1) $ so that by \eqref{eq:Pwr}, 
\begin{align*}
	\Rect(T(w)^{(2)}) = P ( k \, (n + k) \dots \, 1 \, (n + 1) \, n \, 2n \, \dots \, (k + 1) \, (n + k + 1) ) \\
 = \ytableausetup{boxsize = 1.8em} \begin{ytableau} 1 & \scriptstyle{k + 1} & *(yellow) y \\ \vdots & \vdots & *(yellow) \vdots \\ k & 2k & *(yellow) \scriptstyle{n + 2k} \\ \scriptstyle{2k + 1} & *(yellow) z \\ \vdots & *(yellow) \vdots \\ n & *(yellow) 2n \\ *(yellow) \scriptstyle{n + 1} \\ *(yellow) \vdots \\ *(yellow) \scriptstyle{n + k} \end{ytableau} \, \tx{ if } k < \left \lfloor \f{n}{2} \right \rfloor, \qquad 
	\begin{ytableau} 1 & \scriptstyle{k + 1} & *(yellow) y \\ \vdots & \vdots & *(yellow) \vdots \\  \scriptstyle{n - k} & n & *(yellow) 2n \\ x & *(yellow) \scriptstyle{n + 1} \\ \vdots & *(yellow) \vdots \\ k & *(yellow) 2k \\ *(yellow) \scriptstyle{2k + 1} \\ *(yellow) \vdots \\ *(yellow) \scriptstyle{n + k} \end{ytableau} \, \tx{ if } k \ge \left \lfloor \f{n}{2} \right \rfloor.
\end{align*}
where $ x = n - k + 1, y = n + k + 1, z = n + 2k + 1 $. Either way, $ \stab(w) = 2 $.

\end{proof}

\begin{Lemma} \label{lem:le2coltabs} Let $ \SYT^{(2)}(n) $ denote the set of size $ n $ standard tableaux with at most 2 columns. Then,
\[
	\# \SYT^{(2)}(n) = \ch{n}{\left \lfloor n/2 \right \rfloor}.
\]

\end{Lemma}

\begin{proof} Because RSK: $ S_n \to \cup_{\lam \vdash n} \SYT(\lam) \x \SYT(\lam) $  is a bijection, 
\begin{align}
\label{eq:numpermswithQ}
	\# \{ w \in S_n : Q(w) = T \} = \# \SYT(\rv(T)).
\end{align}
For all $ k = 0, 1, \dots, \left \lfloor \f{n}{2} \right \rfloor $, let 
\[
	Q_k = \ytableausetup{boxsize = 1.8em} \begin{ytableau} 1 & \scriptstyle{m + 1} \\ \vdots & \vdots \\ k & \scriptstyle{m + k} \\ \vdots \\ m \\ x \\ \vdots \\ n  \end{ytableau} \,  \tx{ where } m = \left \lfloor \f{n}{2} \right \rfloor, \, x = m + k + 1
\]
which has shape $ (2^k, 1^{n - 2k}) $. Now, $ \{ Q_k : k = 0, 1, \dots, \left \lfloor \f{n}{2} \right \rfloor \} $ is the set of size $ n $ standard tableaux with descents at $ 1, \dots, \left \lfloor \f{n}{2} \right \rfloor - 1, \left \lfloor \f{n}{2} \right \rfloor + 1, \dots, n - 1 $. Thus, 
\begin{align*}
	\# \SYT^{(2)}(n) & = \sum_{k = 0}^{\left \lfloor \f{n}{2} \right \rfloor} \# \SYT(2^k, 1^{n - 2k}) \\
	& = \# \{ w \in S_n : Q(w) = Q_k \tx{ for some } k \} \qquad \tx{ by } \eqref{eq:numpermswithQ} \\
	& = \# \lb w \in S_n : 1, \dots, \left \lfloor \f{n}{2} \right \rfloor - 1, \left \lfloor \f{n}{2} \right \rfloor + 1, \dots, n - 1 \in \Des(Q(w)) \rb \\
	& = \# \lb w \in S_n : 1, \dots, \left \lfloor \f{n}{2} \right \rfloor - 1, \left \lfloor \f{n}{2} \right \rfloor + 1, \dots, n - 1 \in \Des(w) \rb \qquad \tx{ by } \eqref{eq:DesQ} \\
	& = \# \{ w \in S_n : w_1 > \dots > w_{ \left \lfloor \f{n}{2} \right \rfloor }, w_{\left \lfloor \f{n}{2} \right \rfloor + 1} > \dots > w_n \} \\
	& = \ch{n}{ \left \lfloor \f{n}{2} \right \rfloor }
\end{align*}
by choosing the set $ \{ w_1, \dots, w_{\left \lfloor \f{n}{2} \right \rfloor} \} $ from $ [n] $, which determines $ w $ uniquely.

\end{proof}

\begin{proof}[Proof of \eqref{eq:numstab2perms}] We can add $ (n + 1) $ to 2 positions to tableaux in $ \SYT(2^r, 1^{n - 2r}) $ for $ r = 0, \dots, \left \lfloor \f{n - 1}{2} \right \rfloor $, but only 1 position to tableaux in $ \SYT \lp 2^{n/2} \rp $, meaning
\begin{align}
\label{eq:addnplus1}
	\# \SYT^{(2)}(n+ 1)  = 2 \sum_{ r = 0}^{ \left \lfloor \f{n - 1}{2} \right \rfloor} \# \SYT(2^r, 1^{n - 2r}) + \# \SYT \lp 2^{n/2} \rp	
\end{align}
Note $ \SYT \lp 2^{n/2} \rp = \varnothing $ if $ n $ is odd. Therefore,
\begin{align*}
	\# \{ w \in S_n : \stab(w) = 2 \} & = \# \{ w \in S_n : Q(w) = T_k \tx{ for some } k \} \\
						& = \sum_{k = 1}^n \# \SYT(\rv(T_k)) \qquad  \tx{ by } \eqref{eq:numpermswithQ} \\
						& = 2 \sum_{r = 1}^{ \left \lfloor \f{n - 1}{2} \right \rfloor} \# \SYT(2^r, 1^{n - 2r}) + \# \SYT \lp 2^{n/2} \rp \\
	 & = 2 \sum_{r = 0}^{ \left \lfloor \f{n - 1}{2} \right \rfloor} \# \SYT (2^r, 1^{ n - 2r}) + \# \SYT \lp 2^{n/2} \rp - 2 \qquad \tx{ as } \# \SYT(1^n) = 1, \\
	& = \# \SYT^{(2)}(n + 1)  - 2 \qquad \tx{ by } \eqref{eq:addnplus1} \\
	&= \ch{n + 1}{\left \lfloor \f{n + 1}{2} \right \rfloor} - 2 \qquad \tx{ by \Cref{lem:le2coltabs}}.
\end{align*}	

 \end{proof}


\section{Open Problems}
\label{sec:open}

We finish by discussing related open questions about tableau stabilization and promotion. While we have proven some of the important properties of tableau stabilization, much remains unknown.
The most glaring open problem is \Cref{conj:stab}. 

\ms
\noindent \textbf{\Cref{conj:stab}}: Any standard skew tableau with $ b $ rows and decreasing row vector stabilizes at $ b $.

\ms

\noindent We have proven this bound is tight for skew tableaux with constant row vectors, see \Cref{thm:stab} and \Cref{ex:stabatb}. Due to \Cref{ex:samecjsdiffhape}, our approach to proving \Cref{thm:stab} does not readily generalize to proving \Cref{conj:stab}.



The distribution of the stabilization statistic remains to be explored as well. 

\begin{Problem} \label{prob:stabdist} What is the distribution of $ \stab $ on tableaux of a fixed skew shape?

\end{Problem}

\noindent We expect \Cref{prob:stabdist} to be especially difficult. The permutation case, i.e. shape $ (1) \cup \dots \cup (1) $, could be more tractable. The triangular array in \Cref{fig:stabtri} describes stabilization's distribution on permutations in $ S_1, \dots, S_8 $. The distribution is unimodal and log-concave for these small cases. While we know the leftmost entry is always 1, is the rightmost entry is always 1 as well? In addition, since $ \stab $ is invariant on dual equivalence classes, $ Q(w) $ determines $ \stab(w) $.

\begin{Conjecture} The permutation $ w = 1 \, 2 \, \dots \, (n - 1) \, n $ is the only permutation with $ \stab(w) = n $.

\end{Conjecture}

\begin{Problem} What is the distribution of $ \stab $ on $ S_n $? 

\end{Problem}

\begin{Problem} Is the distribution of $ \stab $ on $ S_n $ unimodal? Is it log-concave? 

\end{Problem}

\begin{Problem} Characterize $ \stab(w) $ directly in terms of $ Q(w)$ for all $ w \in S_n $.

\end{Problem}

We have made substantial progress on the problem of specifying the fixed points of the powers of promotion on rectangular tableaux, but some cases remain open. We constructed the tableaux in $ \SYT((ar)^b)^{\pr^{br}} $ for $ a \ge 2b - 1 $ in \Cref{sec:prfixedtabs} and for $ a = 2 $ in \Cref{sec:otherprfixed}. The case $ a = 1 $ is trivial, and the complete case $ r = 1 $ was solved by Purbhoo and Rhee, see \cite{MR3625918}. To completely solve the problem of specifying the fixed points of the powers of promotion on rectangular tableaux, \Cref{prob:otherprfixed} remains.


\begin{Problem} \label{prob:otherprfixed} Fix $ b, r \in \bZ_{\ge 1} $. For each $ a \in [3, 2b - 2] $ ( $a \in [3, b - 1] $ for $ r = 1 $), describe the tableaux in $ \SYT((ar)^b)^{\pr^{br}} $.

\end{Problem}
 
Generalizing Purbhoo and Rhee's construction for $ r = 1, a \ge b $ to $ r \ge 2, a \ge b $ is a nontrivial potential future task. By \Cref{cor:prfixedribbontab} and \eqref{eq:aquotientsufflarge}, we have
\[
	\# \SYT((ar)^b)^{\pr^{br}} = \ch{br}{r, r, \dots, r}, \qquad \tx{ for all } a \ge b.
\]
Purbhoo and Rhee's proof relies on the fact that for all $ T \in \SYT(b^b)^{\pr^b} $, the entries in the anti-diagonal cells $ (b, 1), (b - 1, 2), \dots, (1,b) $ of $ T $ form a permutation in $ S_b $ when taken modulo $ b $. For example,
\[
	\ytableausetup{boxsize = 1.5em} \Yvcentermath1 \begin{ytableau} 1 & 2 & 3 & *(yellow) 7 & *(yellow) 8 \\ 4 & 5 & *(yellow) 10 & *(green) 12 & *(green) 13 \\ *(yellow) 6 & *(yellow) 9 & *(green) 15 & *(orange) 17 & *(orange) 18 \\ *(green) 11 & *(green) 14 & *(orange) 20 & *(cyan) 22 & *(cyan) 23 \\ *(orange) 16 & *(orange) 19 & *(cyan) 21 & *(cyan) 24 & *(cyan) 25 \end{ytableau} \in \SYT(5^5)^{\pr^5}
\]
has the anti-diagonal entries $ [16, 14, 15, 12, 8] \equiv_5 [1,4,5,2,3] \in S_5 $. One might hope that for any $ T \in \SYT((b \dd r)^b)^{\pr^{br}} $, the cells $ (b,1), \dots, (b,r), (b - 1, r + 1), \dots, (b - 1, 2r), \dots, (1, (b - 1)r + 1), \dots (1, br) $ have entries that form a permutation in $ S_{br} $ when taken modulo $ br $. But, considering such examples as
\[
	\begin{ytableau} 1 & 2 & 4 & 5 & 6 & *(yellow) 8 \\ 3 & *(yellow) 7 & *(yellow) 9 & *(yellow) 10 & *(yellow) 11 & *(yellow) 12 \end{ytableau} \in \SYT(6^2)^{\pr^6}, \quad \begin{ytableau} 1 & 2 & 5 & 6 & *(yellow) 8 & *(green) 14 \\ 3 & *(yellow) 7 & *(yellow) 9 & *(yellow) 11 & *(yellow) 12 & *(green) 15 \\ 4 & *(yellow) 10 & *(green) 13 & *(green) 16 & *(green) 17 & *(green) 18 \end{ytableau} \in \SYT(6^3)^{\pr^{6}},
\]
this is not the case, since $ [3,7,9,5,6,8] \equiv_6 [3,1,3,5,6,2] \notin S_6 $ and $ [4, 10, 9, 11, 8, 14] \equiv_6 [4,4,3,5,2,2] \notin S_6 $. In the first case, there is not even a symmetric choice of 3 consecutive entries in row 1 and 3 consecutive entries in row 2 which reduces to a permutation modulo $ 6 $. We want a symmetric choice with respect to reflection so we can attach 2 symmetric shapes together to make a rectangle like we did for tableau stabilization.  

The case $ a < b $ is even more challenging since $ \# \SYT((ar)^b)^{\pr^{br}} = \# \SYT(Q_a((ar)^b)) $ changes according to \Cref{lem:rectquotients}, using \Cref{cor:prfixedribbontab}. Since now the pieces of the quotient are rectangles, one might have to do some tableau stabilization-like procedure with rectangles to find $ \SYT((ar)^b)^{\pr^{br}} $, and it is so far unclear how this would work without just reducing to tableau stabilization. When $ a = 2 $, row-concatenating the rectification to the anti-rectification works, but this does not easily generalize to $ a \in [3, 2b - 1] $. Row-concatenating $ \Rect(S^{\lp \lceil \f{a}{2} \rceil \rp}) $ and $ \Rect^\ast(S^{\lp \lfloor \f{a}{2} \rfloor \rp}) + \lp \lceil \f{a}{2} \rceil \rp br $ will not produce a rectangular tableaux in general.

For example, consider $ a = 3, r = 2, b = 6 $, whence
\[
	\# \SYT(6^6)^{\pr^{12}} = \# \SYT(Q_3(6^6)) = \# \SYT((2,2) \cup (2,2) \cup (2,2)). 
\]
by \Cref{cor:prfixedribbontab} and \Cref{lem:rectquotients}. Choose
\[
	S = \begin{ytableau} \none & \none & \none & \none & 5 & 9 \\ \none & \none & \none & \none & 6 & 11 \\ \none & \none & 2 & 4 \\ \none & \none & 8 & 10 \\ 1 & 3 \\ 7 & 12 \end{ytableau} \, \in \SYT((2,2) \cup (2,2) \cup (2,2)). 
\]
Then, we have
\begin{align*}
	\Rect(S^{(2)}) = \begin{ytableau} 1 & 2 & 4 & 5 & 6 & *(yellow) 16 & *(yellow) 17 & *(yellow) 18 \\ 3 & 8 & 9 & 11 & *(yellow) 20 & *(yellow) 21 & *(yellow) 23 \\ 7 & 10 & *(yellow) 13 & *(yellow) 14 & *(yellow) 22 \\ 12 & *(yellow) 15 \\ *(yellow) 19 & *(yellow) 24 \end{ytableau} \, , \qquad
	\Rect^{\ast}(S) + 24 = \begin{ytableau}  \none \\ \none \\  \none & \none & \none & \none & *(green) 30 \\ \none & \none & \none & *(green) 26 & *(green) 33 \\ \none & *(green) 27 & *(green) 28 & *(green) 29 & *(green) 35 \\ *(green) 25 & *(green) 31 & *(green) 32 & *(green) 34 & *(green) 36 \end{ytableau}. 
\end{align*}
Row-concatenating these along 6 rows gives
\[
	\begin{ytableau} 1 & 2 & 4 & 5 & 6 & *(yellow) 16 & *(yellow) 17 & *(yellow) 18 \\ 3 & 8 & 9 & 11 & *(yellow) 20 & *(yellow) 21 & *(yellow) 23 \\ 7 & 10 & *(yellow) 13 & *(yellow) 14 & *(yellow) 22 & *(green) 30 \\ 12 & *(yellow) 15 & *(green) 26 & *(green) 33 \\ *(yellow) 19 & *(yellow) 24 & *(green) 27 & *(green) 28 & *(green) 29 & *(green) 35 \\ *(green) 25 & *(green) 31 & *(green) 32 & *(green) 34 & *(green) 36 \end{ytableau} \, , 
\]
which is not even partitioned-shaped, let alone rectangular.

In \Cref{thm:prfixedtabs} and \Cref{thm:fixedbyab/2prom}, our construction of $ \SYT((ar)^b)^{\pr^{br}} $ for $ a \ge 2b - 1 $ and White and Rhee's construction for $ a = 2 $ included a natural bijection $ R_a: \SYT(Q_a((ar)^b)) \to \SYT((ar)^b)^{\pr^{br}} $. Moreover, promotion commuted with $ R_a $ as in \Cref{cor:prfixedpr} and \Cref{cor:prfixedpr2}. One would hope these properties extended to a construction of $ \SYT((ar)^b)^{\pr^{br}} $ for $ a \in [3,2b - 2] $ as well.

\begin{Problem} Fix $ b, r \in \bZ_{\ge 1} $ and $ a \in [3, 2b - 2] $. Find a bijection 
\[
	R_a: \SYT((ar)^b)^{\pr^{br}} \to \SYT(Q_a((ar)^b))
\]
satisfying $ \pr(R_a(S)) = R_a(\pr(S)) $ for all $ S \in \SYT(Q_a((ar)^b)) $, or show no such bijection exists.
\end{Problem}

\section*{Acknowledgments}

I sincerely thank my advisor, Sara Billey, for helpful discussions and extensive comments on the manuscript. I thank Dennis White very much for his work on the problem of rectangular tableaux fixed by promotion powers, notes that inspired \Cref{thm:prfixedpieces}, helpful notes on proper citation, and for conjecturing \Cref{thm:rhoadesCSP} in the first place. We would also like to thank Brendon Rhoades for proving \Cref{thm:rhoadesCSP} and posing \cite[Problem~9.4]{MR2557880}, which inspired this work. I would also like to thank Kevin Purbhoo for making and pointing out to me previous progress on this problem in \cite{MR3625918}.

\bibliography{refs}{}
\bibliographystyle{alpha}

\end{document}